\documentclass{article}

\usepackage{algpseudocode}
\usepackage{alltt}
\usepackage{amsbsy}
\usepackage{amscd}
\usepackage{amsfonts}
\usepackage{amsmath}
\usepackage{amssymb}
\usepackage{amsthm}
\usepackage{boldtensors}
\usepackage{booktabs}
\usepackage{color}
\usepackage{enumitem}
\usepackage{fancyhdr}
\usepackage{float}
\usepackage{geometry}
\usepackage{graphicx}
\usepackage{microtype}
\usepackage{multirow}
\usepackage{siunitx}
\usepackage{supertabular}
\usepackage{times}
\usepackage{bm}
\usepackage[ruled,vlined]{algorithm2e}
\usepackage[dvipsnames]{xcolor}
\usepackage{tikz}
\usepackage{tikz-3dplot}
\usepackage{url}
\usepackage{verbatim}
\usepackage{xspace}
\usepackage{soul}

\usepackage[nomessages]{fp}
\usepackage[english]{babel}
\usepackage{subcaption} 
\usepackage[T1]{fontenc}
\usepackage{ae,aecompl}
\usepackage[abs]{overpic}
\usepackage[normalem]{ulem}
\usepackage[font=footnotesize]{caption}

\usepackage{marginnote} 

\usepackage[colorlinks]{hyperref}
\usepackage{cleveref}

\usetikzlibrary{arrows,calc,decorations.markings,math,arrows.meta}
\usetikzlibrary{shapes.geometric,patterns}

\crefname{paragraph}{section}{sections}
\Crefname{paragraph}{Section}{Sections}

\newcommand{\uvec}{\mathbf{u}}

\newcommand{\M}{\boldsymbol M}

\newcommand{\evec}{\mathbf{e}}

\newcommand{\I}{\boldsymbol{I}}

\newcommand{\bvec}{\mathbf{b}}

\newcommand {\norm}[1]{\left\lVert#1\right\rVert}

\newcommand{\entryvec}[2]{\left[#1_1\ \cdots\ #1_{#2}\right]^T}

\newcommand{\PetrovGalerkin}[1]{{#1}_\text{PG}}
\newcommand{\PrecondPetrovGalerkin}[1]{{#1}_\text{PPG}}

\newcommand{\solsymb}{w}
\newcommand{\sol}{\boldsymbol \solsymb}
\newcommand{\res}{\boldsymbol r}
\newcommand{\dummySol}{\boldsymbol y}

\newcommand{\approximate}{\tilde}
\newcommand{\reduced}{\hat}
\newcommand{\reference}[1]{{#1}_0}

\newcommand{\newton}[1]{^{(#1)}}

\newcommand{\solredi}[1]{\hat\solsymb_{#1}}

\newcommand{\RR}[1]{\mathbb{R}^{#1}}
\newcommand{\ndof}{N}
\newcommand{\nrb}{M}

\newcommand{\rb}{\bm{\Phi}} 
\newcommand{\lb}{\bm{\Psi}} 

\newcommand{\Mk}{\M^{(k)}}

\newcommand{\rk}{\res^{(k)}}
\newcommand{\jac}{\boldsymbol J}
\newcommand{\jk}{\jac\newton{k}}
\newcommand{\jkT}{{\jk}^T}
\newcommand{\jkinv}{{\jk}^{-1}}

\newcommand{\lengthk}{\alpha_k}
\newcommand{\defeq}{:=}

\newcommand{\trilinos}{\textsc{Trilinos}}

\newcommand{\aeras}{\textsc{Aeras}}
\newcommand{\albany}{\textsc{Albany}}
\newcommand{\ace}{\textsc{ACE}}

\newcommand{\aztec}{\textsc{AztecOO}}

\newcommand{\felix}{\textsc{ALI}}

\newcommand{\ifpack}{\textsc{Ifpack}}

\newcommand{\lcm}{\textsc{LCM}}

\newcommand{\qcad}{\textsc{QCAD}}
\newcommand{\sacado}{\textsc{Sacado}}

\newcommand{\anasazi}{\textsc{Anasazi}}
\newcommand{\loca}{\textsc{LOCA}}
\newcommand{\rbgen}{\textsc{RBGen}}
\newcommand{\exodus}{\textsc{Exodus}}
\newcommand{\op}[1]{\operatorname{\text{#1}}}
\newcommand{\Grad}{\op{Grad}}

\newcommand{\Div}{\op{Div}}

\newcommand{\e}[1]{\times10^{#1}} 

\newcommand{\mass}[1][\mbox{ }]{#1\mbox{kg}}
\newcommand{\len}[1][\mbox{ }]{#1\mbox{m}}
\renewcommand{\temp}[1][\mbox{ }]{#1\mbox{K}}
\renewcommand{\time}[1][\mbox{ }]{#1\mbox{s}}
\newcommand{\press}[1][\mbox{ }]{#1\mbox{Pa}}
\newcommand{\dens}[1][\mbox{ }]{#1\mass[]/\len[]^3}

\newcommand{\units}[1]{\mbox{ }[#1[]]}

\providecommand{\keywords}[1]{\textbf{Keywords:} #1}

\newcommand{\subfig}[3]{ 
    \begin{subfigure}[c]{1.0\textwidth}
        \centering
        \includegraphics[width=\textwidth]{#1}
        \caption{#2}
        \label{#3}
    \end{subfigure}
} 

\newenvironment{varray}{\left(\begin{array}{c}}{\end{array}\right)} 
\newenvironment{marray}{\left(\begin{array}{cc}}{\end{array}\right)} 

\pdfoutput=1

\newif\ifFiguresInGrid
\FiguresInGridtrue
\newif\ifLegendInFig
\LegendInFigtrue

\begin{document}

\pagestyle{fancyplain}

\lhead [\fancyplain{}{\emph{P. Lindsay, J. Fike, I. Tezaur, K. Carlberg}}]
{\fancyplain{}{\emph{P. Lindsay, J. Fike, I. Tezaur, K. Carlberg}}}

\rhead
[\fancyplain{}{\emph{Preconditioned LSPG ROMs}}]
{\fancyplain{}{\emph{Preconditioned LSPG ROMs}}}

\title{Preconditioned Least-Squares Petrov--Galerkin Reduced Order Models}

\author{\Large Payton Lindsay$^1$,
  Jeffrey Fike$^2$, Irina Tezaur$^3$\thanks{Email: ikalash@sandia.gov},
  Kevin Carlberg$^4$
  \\
  \\
  $^1$Computational Solid Mechanics and Structural Dynamics Department\\
  Sandia National Laboratories\\
  Albuquerque, NM 87185, USA\\
  \\
  $^2$Aerosciences Department\\
  Sandia National Laboratories\\
  Albuquerque, NM 87185, USA\\
  \\
  $^3$Quantitative Modeling and Analysis Department\\
  Sandia National Laboratories\\
  Livermore, CA 94551, USA\\
  \\
  $^4$University of Washington\\
  Seattle, WA 98195, USA\\
}

\date{\today}

\maketitle

\begin{abstract}

In this paper, we introduce a methodology for
improving the accuracy and efficiency 
of reduced-order models (ROMs)
 constructed using the least-squares Petrov--Galerkin (LSPG) projection
method 
through the introduction of
preconditioning.  Unlike prior related work, 
which focuses on preconditioning 
the linear systems arising within the ROM numerical solution procedure
to improve linear solver performance, 
our approach
leverages a preconditioning matrix directly within the minimization problem
underlying the LSPG formulation.  
Applying preconditioning in this way has the potential to improve
 ROM accuracy for several reasons.
First, preconditioning the LSPG formulation 
changes the norm defining the residual minimization,
which can improve the
residual-based stability constant bounding the ROM solution's error.
The incorporation of a preconditioner into the LSPG formulation can have the additional effect of scaling
the
components of the residual being minimized to make them roughly of the same magnitude,
which can be beneficial when applying 
the LSPG method to problems with disparate scales
(e.g., dimensional equations, multi-physics problems).
Importantly, we demonstrate that an `ideal preconditioned' LSPG ROM (a ROM in which the
preconditioner is the inverse of the Jacobian of its corresponding full order model)
emulates projection of the full order model solution increment onto
the reduced basis.  This quantity defines a lower bound on the error of a ROM
solution for a given
reduced basis.
By designing preconditioners that approximate the Jacobian inverse---as is
common in designing preconditioners for solving linear systems---it is possible to obtain a ROM whose error approaches
this lower bound.
The proposed approach is evaluated on several mechanical and thermo-mechanical 
problems implemented within the \albany{} HPC code and run in the predictive regime,
with prediction across material parameter space.  We demonstrate numerically 
that the introduction
of simple Jacobi, Gauss-Seidel and ILU preconditioners into the Proper Orthogonal 
Decomposition (POD)/LSPG formulation
reduces significantly the ROM solution error, the reduced Jacobian condition number, 
the number of nonlinear iterations required to reach convergence, and the wall time
(thereby improving efficiency).  
Moreover, our numerical results reveal that the introduction 
of preconditioning can deliver a robust and accurate solution for 
test cases in which the unpreconditioned LSPG method fails to converge.

\end{abstract}

\keywords{Reduced order model (ROM), proper orthogonal decomposition (POD), 
Least-squares Petrov--Galerkin (LSPG) projection, preconditioners, computational
solid mechanics, mechanical, thermo-mechanical.}

\section{Introduction}
\label{sec:introduction}


Numerous modern-day science and engineering problems require the simulation of
complex systems with tens of millions of unknowns.  Despite improved
algorithms and the availability of massively parallel computing resources,
``high-fidelity'' models are, in practice, often too computationally expensive
for use in time-critical settings such as design, fast turnaround analysis and
control. The situation is particularly grave in applications involving
multi-query analyses (e.g., optimization, uncertainty quantification), which
require simulations to be repeated many times to explore the design space or to
properly characterize uncertainty. The necessary calculations can present an
intractable computational burden even with the projected growth in computing
power as we approach the exascale computing age.

Reduced-order modeling is a promising strategy for reducing the computational
cost of such simulations while preserving high levels of fidelity.
Reduced-order models (ROMs) are models constructed from
high-fidelity simulations that retain the essential physics and dynamics of
their corresponding full order models (FOMs), but have a much lower
computational cost.  Although some consider a data-fit or low-fidelity model a
ROM, herein the term ``reduced order model'' refers to a projection-based ROM.
In projection-based model reduction, the state variables are restricted to
reside in a low-dimensional subspace, typically computed offline through a
data-compression process performed on a set of snapshots collected from a
high-fidelity simulation or physical experiment, followed by truncation. There
are numerous approaches in the literature for computing a low-dimensional
subspace e.g., Proper Orthogonal Decomposition (POD) \cite{sirovich,holmes},
Dynamic Mode Decomposition (DMD) \cite{rowley1,Schmid_JFM_2010} balanced POD
(BPOD) \cite{rowley, willcox}, balanced truncation \cite{gugercin, moore}, and
the reduced basis method (RBM)  \cite{rozza, veroy}.  Herein, without loss of
generality, we restrict attention to the POD approach for calculating the
reduced bases due to its prevalence and simplicity.  Once a reduced basis is
computed, the ROM dynamical system is obtained by projecting the governing
equations, or some discretized form of these equations, onto the
low-dimensional subspace. For nonlinear problems, an additional approximation,
referred to as ``hyper-reduction'', is usually required to gain a
computational speed-up.   Since ROMs are, by construction, low-dimensional and
inexpensive to evaluate, these models can enable real-time analysis and
alleviate the computational burden posed by high-dimensional uncertainty
quantification (UQ) problems critical to many applications.

In order to serve as a viable predictive tool, a ROM should possess
certain fundamental mathematical properties, such as consistency
(with respect to its corresponding high-fidelity model), stability (in space
as well as in time),
and convergence (to the solution of its corresponding high-fidelity
model).  
Galerkin projection, the most popular approach, is
considered continuous optimal, as the resulting ROM minimizes the
time-continuous residual in the $l^2$ norm \cite{carlbergGalDiscOpt}, and preserves problem
structure \cite{kevin_lagrangian, Carlberg_SIAM_2015,lall} for
certain classes of problems (e.g., Lagrangian dynamical systems).
However, 
the method
can give rise to nonphysical instabilities \cite{Rempfer_JFM_1994,
baronekalash, baroneJCP, ec_ldrd_sand, romLDRDSAND, maciej,
maciejirina,kevin_gnat} and inaccurate long-time responses \cite{sirisup,
noack_papas_monkewitz_2005,carlbergGalDiscOpt}.  Additionally, it  lacks in general an \textit{a priori}
convergence guarantee \cite{rowley1}. The least-squares
Petrov--Galerkin (LSPG) projection method has been proposed
\cite{CarlbergGappy} to remedy some of these difficulties
via symmetrization of the discrete Jacobian. This
method performs projection at the level of the fully discrete
partial differential equations (PDEs), i.e., after the PDEs have
been discretized in space and time, and computes a solution that
minimizes the $\ell^2$-norm of the time-discrete residual arising in
each time step.  This procedure ensures that adding basis vectors
yields a monotonic decrease in the least-squares objective function
defining the underlying minimization problem.  The method can
maintain efficiency through the incorporation of hyper-reduction approaches
such as gappy POD
\cite{gappy}, the discrete empirical interpolation method
(DEIM) \cite{deim} or the Energy-Conserving mesh Sampling and Weighing (ECSW) method
\cite{FarhatChapman:2015, Chapman:2017, FarhatCortial:2014},
introduced to maintain
efficiency of the method for nonlinear problems.  When combined with
gappy POD, the LSPG approach is equivalent to the Gauss--Newton with
approximated tensors (GNAT) method \cite{kevin_gnat}.  While LSPG
projection does not necessarily guarantee \textit{a priori} accuracy
and stability,
it has been shown
numerically to possess better stability and accuracy properties than
the Galerkin method for a variety of problems  
\cite{CarlbergGappy,kevin_gnat, fy16_sand, fy17_sand, fy18_sand,
carlbergGalDiscOpt}. We note that the LSPG method has been extended in the
time domain \cite{kevin_space-time,parish2021windowed}, to allow for nonlinear trial manifolds rather than linear trial
subspaces \cite{kookjin}, to operate in a domain-decomposition setting \cite{hoang2021domain}, and to account for physics-based constraints \cite{schein2021preserving}.


The aim of the present paper is to develop a methodology for
improving the accuracy of ROMs constructed using the LSPG projection
method for a wide range of applications through the introduction of
preconditioning. As shown in \cite{carlbergGalDiscOpt}, LSPG errors are
subject to a stability constant that is dictated by the residual.
The use of preconditioning within the LSPG formulation improves the
conditioning of the system, which improves the stability constants.
Intuitively, this aligns the $\ell^2$-norm of the residual more
closely with the $\ell^2$-norm of the (time-local) state error.
To motivate the main contributions of this paper, we first provide 
a brief overview of related work involving projection-based ROMs 
and preconditioning.

\subsection{Overview of related work} \label{sec:overview_past_work}

The idea of preconditioning began to make an appearance within the
model-reduction literature approximately a decade ago. The bulk of
the literature on this subject falls into two categories, which we
overview succinctly below.

The first category of methods are aimed at developing ROM-based
preconditioners for high-fidelity models.
 In \cite{carlbergKrylov2015,pasetto}, 
the authors present a novel
class of ROM-based preconditioners for the iterative high-fidelity
solution of transient parabolic and self-adjoint PDEs.  The
preconditioners are obtained by nesting appropriate projections of
ROMs into the classical preconditioned Conjugate Gradient (CG)
iteration within the FOM.  Several authors explore the idea of
developing a ROM-based preconditioner to speed-up the numerical
solution of the pressure system arising in high-fidelity reservoir
modeling \cite{astrid, jiang}. In \cite{astrid}, Astrid \textit{et
al.} employ ROM concepts to develop a new preconditioner to speed up
the pressure solution during a two-stage preconditioning procedure
in a fully implicit time-integration scheme within a high-fidelity
reservoir simulator.  In \cite{jiang}, Jiang introduces a ROM-based
preconditioner within a stationary Richardson iteration for the
pressure system in a similar high-fidelity reservoir simulation.

The second category of methods explore the idea of introducing a
preconditioner within the ROM workflow to speed up the online
evaluation of the ROM.  This work is motivated by the observation
that the repeated solution of linear systems, characterized by dense matrices, can be the main
computational bottleneck in efficient scaling of ROMs \cite{singh1}.
Commonly, the linear systems arising within the ROM workflow are
solved using direct methods.  This approach can be reasonable
 for small ROMs, but for larger ROMs containing numerous
parameters, it was shown in \cite{elman1} that a direct method can
be more expensive to solve than the original FOM. Preconditioned
iterative methods have been developed for alleviating this
difficulty.  In \cite{elman2} and \cite{elman1}, preconditioners
that are precomputed in the offline stage of the model reduction are
introduced into a reduced basis ROM formulation with and without the
application of discrete empirical interpolation to handle the
nonlinear terms in the ROM.  These preconditioners are derived from
preconditioners used in solving the FOM, and are shown to be more
efficient than direct solvers for problems depending on a moderately
large number of parameters.  Singh \textit{et al.} develop
strategies for reusing preconditioners in projection-based model reduction
algorithms for parametric \cite{singh2} and non-parametric
\cite{singh1} ROMs to improve scalability of ROM solution over
direct methods.  These authors introduce a novel Sparse Approximate Inverse
(SPAI) preconditioner for ROMs constructed via moment matching using
the Bilinear Iterative Rational Krylov Algorithm (BIRKA) algorithm.
These preconditioners are aimed at systems where the linear system
may change substantially between nonlinear iterations, as happens
with highly nonlinear problems.  Several other authors consider
preconditioners constructed using IRKA methods, for which the most
computationally demanding part is the solution of linear shifted
systems.  The ROM matrices in this case are sparse.  Ahmad
\textit{et al.} \cite{ahmad} propose variants of the IRKA that use
preconditioned versions of the shifted BiCG algorithm, and develop a
polynomial preconditioner that can simultaneously be applied to all
shifted systems.  The idea of updating incomplete factorization
preconditioners within the IRKA approach is explored in \cite{anzt}.
The use of randomized matrices to construct preconditioners for ROMs
has been studied recently by several authors \cite{zahm, balabanov}.
In \cite{zahm}, the authors propose to use randomized linear algebra
for the efficient construction of a preconditioner by interpolation
of matrix inverse.  It is shown that, in addition to improving the
online performance of a given ROM, these preconditioners can be also
used to improve the quality of residual-based error estimates in the
context of the reduced basis method.  Building on \cite{zahm},
Balabanov \cite{balabanov} addresses the efficient construction
(using random sketching) of parameter-dependent preconditioners that
can be used to improve the quality of Galerkin projections or for
effective error certification for reduced basis and POD ROMs with
ill-conditioned operators.  

It is worthwhile to remark that the idea of
preconditioning an LSPG ROM is related to several other recently-introduced methodologies for improving
ROM accuracy, namely the idea of scaling the
ROM residual prior to its minimization and the idea of modifying 
the norm in which the residual is minimized.  In \cite{kyle}, 
Washabaugh proposed an approach for improving the accuracy of ROMs
for systems of PDEs in which 
different
variables have drastically different scales (e.g., applications
where dimensional equations are employed).   Specifically, it was shown 
in this reference that 
scaling the ROM residual to get all the equations to be
roughly the same order of magnitude (equivalent to the application
of a diagonal preconditioner) can result in more accurate ROMs.
In a similar vein, it was demonstrated in \cite{Blonigan:2022}
that normalizing snapshot data prior to computing POD modes in order 
to remove dimensional effects increases numerical robustness of ROMs for 
hypersonic aerodynamic simulations, especially in
cases with non-equilibrium thermo-chemical effects.
More recently, in \cite{Parish:2022}, Parish and Rizzi developed a set of dimensionally-consistent inner products
and demonstrated that these inner products have a positive 
impact on both Galerkin and LSPG ROMs with application to the compressible Euler equations.
In \cite{Collins:2020}, Collins \textit{et al.} examined an optimal test basis 
for LSPG ROMs that are obtained by solving an adjoint-like system, which
transforms the residual error minimization problem to be equivalent to a least-squares
state error minimization problem.  The resulting test functions can be seen 
as incorporating the action of a preconditioner matrix that approximates
the inverse of the Jacobian matrix from a ROM's corresponding FOM.
In \cite{Abgrall:2018}, Abgrall \textit{et al.}  proposed a residual minimization model reduction approach, in which
the norm in which the residual is minimized is the $L^1$ norm.  If the resulting $L^1$ residual minimization
is solved iteratively using reweighted least-squares, the procedure can be interpreted as a least-squares solve.

\subsection{Contributions and organization} \label{sec:contrib}

A distinguishing feature of the present work that differentiates it
from most of the references described above is the premise that the
introduction of  a preconditioner into the LSPG formulation can not
only lead to performance gains (by reducing the condition number of the 
reduced Jacobian), but can actually improve a ROM's accuracy.
Unlike prior work, e.g., \cite{elman1, elman2,singh1}, the idea herein is not to simply precondition 
the linear systems arising within the LSPG algorithm to improve linear solver performance, but rather to 
insert a preconditioner matrix directly into the minimization problem
underlying the LSPG formulation.  
Doing this has the
effect of changing the norm defining the residual minimization
problem underlying the LSPG method, which can improve the
residual-based stability constant bounding the method's error.
We additionally demonstrate that: 
\begin{itemize}
\item an ideal preconditioned LSPG ROM (a ROM in which the 
preconditioner is the inverse of the Jacobian of its corresponding FOM)
emulates projection of the FOM solution increment onto 
the reduced basis, which defines a lower bound on the error of a ROM for a given  
reduced basis; 
\item it is possible to obtain a ROM whose error approaches
this lower bound by improving the quality of the preconditioner employed 
within the LSPG formulation;
\item 
		the addition of certain preconditioners into the LSPG formulation can have the effect of scaling 
the 
components of the residual being minimized to make them roughly of the same magnitude, which
can minimize bias and reduce the number of nonlinear iterations required for 
convergence  when applying the LSPG method to problems with disparate scales 
(e.g., dimensional PDEs, multi-physics problems);
\item while our approach is not specifically based on preconditioning the linear 
systems arising within the LSPG iteration process, it has the effect of reducing 
the condition numbers of an LSPG ROM's reduced Jacobians.
\end{itemize}
We evaluate the performance and predictive 
accuracy of several preconditioned 
POD/LSPG ROMs on one mechanical and two thermo-mechanical problems 
implemented within the open-source HPC multi-physics finite element code known
as \albany{} \cite{albany}.  We point out that the vast majority of
the current ROM literature uses Galerkin projection for the 
reduction step of the ROM workflow when considering solids and structures
\cite{Lall:2003, Tiso:2013, FarhatChapman:2015, Chapman:2017, FarhatCortial:2014,
CarlbergAIAA:2012, Carlberg:2015}; 
herein, we demonstrate that the preconditioned LSPG 
method can also be effective for this class of problems.
For the mechanical problem, 
we demonstrate that
preconditioning is needed to obtain 
a convergent POD/LSPG ROM.  For the thermo-mechanical problems, 
	we
demonstrate that, 
while it is possible 
to obtain a convergent unpreconditioned POD/LSPG ROM for certain basis
dimensions, 
the addition of preconditioning reduces the error by between two and seven orders
of magnitude, all while reducing the overall wall time by as much as $12 \times$.  
The wall time improvements can be attributed to a preconditioner's ability
to reduce not only the condition number of the reduced Jacobians arising within
the LSPG algorithm, but also the total number of nonlinear iterations 
required for ROM convergence.

The remainder of this paper is organized as
follows.  In Section \ref{sec:formulation}, we provide an overview of the
POD/LSPG approach to model reduction applied to a
generic system of nonlinear algebraic
equations, which might arise from the discretization of a PDE.
We introduce the notion of preconditioning
within the LSPG residual minimization problem in Section
\ref{sec:precond}.  We demonstrate that the use of a preconditioner
has the effect of changing the norm in which the LSPG residual
minimization problem is solved, and show that the use of an ideal
preconditioner (the inverse of a Jacobian matrix in a particular
Newton iteration) will lead to an $l^2$-optimal projection of the
solution increment onto the reduced basis.  In Section
\ref{sec:impl}, we describe succinctly the \albany{} multi-physics
finite element code and our implementation of model reduction
capabilities within this code, including 
a partitioned/``blocking vector" approach for applying ROM 
Dirichlet boundary conditions  within a finite element 
code that does not remove the constrained (Dirichlet) degrees 
of freedom from the global finite element system prior to performing
its numerical solution.
We study the performance of the
proposed preconditioned LSPG ROMs in the context of several
mechanical and thermo-mechanical solid mechanics problems in Section
\ref{sec:numerical}, demonstrating the efficacy of the proposed 
preconditioning strategy.  Finally, conclusions are offered in Section
\ref{sec:conclusion}.

\section{Problem formulation}
\label{sec:formulation}

\subsection{Full-order model} \label{sec:fom}


Since the methodology described herein can be applied to a wide
range of FOMs arising in a variety of
applications, let us consider the following generic system of nonlinear
algebraic equations defining our FOM:
\begin{equation} \label{fom}
\res(\sol) = \boldsymbol{0}.
\end{equation}
Here, $\sol\in\RR{\ndof}$ is the state vector and
$\res:\RR{\ndof}\rightarrow\RR{\ndof}$ is the nonlinear residual
operator.  Systems of the form \eqref{fom} are obtained by fully
discretizing a set of governing PDEs using a
spatial discretization method and, in the case of a dynamics problem,
a time-integration scheme.  Assuming we solve \eqref{fom} using a
(globalized) Newton's method, the following sequence of solutions is
generated:
\begin{equation}\label{fom_newton}
  \begin{aligned}
     \jk\delta\sol\newton{k} &= -\rk,\quad k=1, \ldots,K,\\
     \sol\newton{k} &= \sol\newton{k-1} + \lengthk\delta\sol\newton{k}.
  \end{aligned}
\end{equation}
Here, 
$\jk\defeq \frac{\partial \res}{\partial \sol}\left(\sol\newton{k}\right)\in\RR{\ndof\times\ndof}$ is the FOM Jacobian,
$\rk\defeq \res\left(\sol\newton{k}\right)\in\RR{\ndof}$ is the FOM residual,
$\sol\newton{0}$ is an initial guess for the solution, 
$\lengthk\in\RR{}$ is the step length (obtained using a line-search method, or set to one, as commonly done),
and $K$ is the number of Newton iterations.




\subsection{Reduced order models} \label{sec:rom}

Our task is to build a ROM for \eqref{fom} using the
projection-based model reduction approach. This approach consists of
three steps:
\begin{enumerate}
\item calculation of a reduced basis,
\item projection of the governing equations (in our case, \eqref{fom})
onto the subspace spanned by the reduced basis, and
\item hyper-reduction to handle efficiently the projection of the
nonlinear terms.
\end{enumerate}
We describe each of these steps succinctly in the following
subsections.  
As mentioned earlier in the introduction, 
we restrict our attention herein to reduced bases calculated using the POD.
The numerical implementation of the model reduction
algorithms outlined in this section within the \albany{}
multi-physics finite element code is discussed later, in Section
\ref{sec:numerical}.  


\subsubsection{Reduced basis calculation via the POD} \label{sec:pod}

The first step in the projection-based approach to model reduction is the
calculation of a basis of reduced dimension $\nrb \ll \ndof$ (where $\ndof$
denotes the number of degrees of freedom (dofs) in the full order model
\eqref{fom}) using the POD.  The POD is a mathematical procedure
that, given an ensemble of data and an inner product, denoted
generically  by $(\cdot, \cdot)$, constructs a basis for the
ensemble that is optimal in the sense that it describes more energy
(on average) of the ensemble in the chosen inner product than any
other linear basis of the same dimension $\nrb$.  The ensemble
$\{\sol^p:p=1,\dots, P\}$ is typically a set of $P$ 
instantaneous snapshots of a numerical solution field, collected for $P$ values of
a parameter of interest, or at $P$ different times. 
 Following the so-called ``method of snapshots" 
\cite{sirovich}, given a snapshot matrix $~W:=[~w^1, ..., ~w^P] \in \mathbb{R}^{N \times P}$, 
a POD basis $\rb_M$
of dimension $\nrb \leq P$ is obtained by first computing the (thin) singular 
value decomposition (SVD)
\begin{equation} \label{eq:svd}
	~W = ~U ~\Sigma ~V^T,
\end{equation}
where the left singular vector matrix $~U := [~u^1, ..., ~u^P] \in \mathbb{R}^{N \times P}$
is orthogonal ($~U^T ~U = ~I$), the diagonal singular value matrix $~\Sigma = \text{diag}(\sigma_i)$ 
contains the ordered singular values $\sigma_1 \geq \sigma_2 \geq ... \geq \sigma_P \geq 0$
of $~W$, and the right singular vector matrix $~V \in \mathbb{R}^{P \times P}$
is also orthogonal, like $~U$ ($~V^T ~V = ~I$).  Given \eqref{eq:svd},
the sought-after POD basis is obtained by selecting the first $M \leq P$ 
left singular vectors of $~W$: 
\begin{equation}
	\rb_M := [u^1, ..., u^P] \in \mathbb{R}^{N \times M}.
\end{equation}
It follows that $\rb_M$ has orthonormal columns and satisfies 
$\rb_M^T \rb_M = ~I$.  


Once $\rb_M$ is calculated, we approximate the solution to
\eqref{fom} by
 \begin{equation} \label{eq:subspace}
\approximate{\sol} = \reference{\sol} + \rb_M\reduced{\sol} = \reference{\sol} +
	 \sum_{i=1}^\nrb{\bm \phi}_{i}\solredi{i}
  \end{equation}
where ${\bm \phi}_i := ~u^i$ for $i=1,..., M$, $\reduced{\sol}\defeq\entryvec{\reduced{\solsymb}}{\nrb}\in\RR{\nrb}$
denote the generalized coordinates, and $\reference{\sol}\in\RR{\ndof}$
denotes a reference solution, often taken to be the initial
condition in the case of an unsteady simulation.  Substituting the approximation \eqref{eq:subspace} into
\eqref{fom} yields
\begin{equation} \label{fomOver}
\res(\reference{\sol} + \rb_M\reduced{\sol}) = \boldsymbol{0}.
\end{equation}
\eqref{fomOver} is a system of $\ndof$ equations in $\nrb$ unknowns $\reduced{\sol}$.
As this is an over-determined system, it may not have a solution.

\subsubsection{Least-Squares Petrov--Galerkin (LSPG) projection} \label{sec:lspg}

As discussed earlier, the LSPG approach \cite{CarlbergGappy}
to model reduction has shown some promise in mitigating the stability problems 
and other disadvantages of other approaches, e.g., the Galerkin projection method.
In the LSPG approach, the ROM solution in \eqref{eq:subspace} is obtained by  solving the following
 least-squares optimization problem:
 \begin{equation} \label{lspg}
     \PetrovGalerkin{\reduced{\sol}} \in \arg\min_{\dummySol\in\RR{\nrb}}\|\res(\reference{\sol} +
     \rb_M\dummySol)\|_2^2.
 \end{equation}
Here, the approximate (Petrov--Galerkin) solution is 
$\PetrovGalerkin{\approximate{\sol}} \defeq \reference{\sol} + \rb_M\PetrovGalerkin{\reduced{\sol}}$.
The name ``LSPG" ROM comes from the observation that
solving \eqref{lspg} amounts to solving a nonlinear least-squares
problem. The two most popular approaches for this are the
Gauss--Newton approach and the Levenberg--Marquardt (trust-region)
method. Following the work of Carlberg \textit{et al.}
\cite{kevin_gnat}, we adopt the Gauss--Newton approach\footnote{The
LSPG approach is the basis for the Gauss--Newton with Approximated
Tensors (GNAT) method of Carlberg \textit{et al.}
\cite{kevin_gnat}.}. This approach implies solving a sequence of
linear least-squares problems of the form
\begin{gather} \label{lspgit}
    \delta\PetrovGalerkin{\reduced{\sol}}\newton{k} \in
    \arg\min_{\dummySol\in\RR{\nrb}}\|\jk\rb_M\dummySol +
    \rk\|_2^2,\quad k=1,\ldots,\PetrovGalerkin{K},\\
    \label{lspgup1}\PetrovGalerkin{\reduced{\sol}}\newton{k} = \PetrovGalerkin{\reduced{\sol}}\newton{k-1} + \lengthk\delta\PetrovGalerkin{\reduced{\sol}}\newton{k},\\
    \label{lspgup2}\PetrovGalerkin{\approximate{\sol}}\newton{k} = \reference{\sol} + \rb_M\PetrovGalerkin{\reduced{\sol}}\newton{k-1},
\end{gather}
where $\PetrovGalerkin{K}$ is the number of Gauss--Newton iterations. It
can be shown that the approximation upon convergence is
$\PetrovGalerkin{\approximate{\sol}} = \PetrovGalerkin{\approximate{\sol}}\newton{\PetrovGalerkin{K}}$ 
and 
$\PetrovGalerkin{\reduced{\sol}} = \PetrovGalerkin{\reduced{\sol}}\newton{\PetrovGalerkin{K}}$.\footnote{In
the event of an unsteady simulation, the initial guess for the
generalized coordinates is taken to be the generalized coordinates
at the previous time step.} Note that the normal equations for
\eqref{lspgit} take the form
\begin{gather} \label{lspgitne}
    \PetrovGalerkin{\jac}\newton{k}\delta\PetrovGalerkin{\reduced{\sol}}\newton{k} =
    -\PetrovGalerkin{\res}\newton{k},
    \quad k=1,\ldots,\PetrovGalerkin{K},
\end{gather}
where 
\begin{equation} \label{reducedJac}
	\PetrovGalerkin{\jac}\newton{k} := {\lb_M\newton{k}}^T \jk\rb_M, \hspace{0.5cm}
	\PetrovGalerkin{\res}\newton{k} := {\lb_M\newton{k}}^T \res\newton{k},
\end{equation}
and 
\begin{equation} \label{eq:psi}
	\lb_M\newton{k}:= \jk\rb_M.
\end{equation}
Equation \eqref{lspgit} 
can be solved numerically using a variety of methods,
including solving the normal equations \eqref{lspgitne}, or by
employing more numerically stable methods such as the QR
decomposition and the SVD; see
\cite{kevin_gnat} for more details. 

\vspace{0.5cm}

\noindent {\bf \textit{Remark 1.}} It is clear from \eqref{reducedJac} and \eqref{eq:psi} that LSPG projection can be interpreted as a 
Petrov--Galerkin process of the Newton
    iteration with trial basis  $\rb_M$ and test basis $\lb_M\newton{k}$
		defined in \eqref{eq:psi}.

\vspace{0.5cm}

\noindent {\bf \textit{Remark 2.}} When Newton's method
is employed to solve the reduced state equations for the LSPG method, 
the ROM solution procedure can be viewed as
minimizing at each iteration the error in the search direction:
\begin{equation} \label{opt}
\delta\reduced{\sol} \in \arg \min_{\dummySol\in\RR{\nrb}}\|\rb_M\dummySol
- \delta\sol\newton{k} \|_{{\bm \Theta}}^2 = \arg
\min_{\dummySol\in\RR{\nrb}}\|\rb_M\dummySol + \jkinv
\rk\|_{{\bm \Theta}}^2,
\end{equation}
where ${\bm \Theta}: = \jkT\jk$ 
and $|| \mathbf{x}||_{{\bm \Theta}}:= \sqrt{\mathbf{x}^T {\bm
\Theta}\mathbf{x}}$.  Since ${\bm \Theta}$ is always a symmetric
positive-definite matrix, the LSPG approach has the effect of symmetrizing 
the discrete ROM Jacobian.
As shown in \cite{CarlbergGappy}, when \eqref{opt} is satisfied, the error measure
$\|\rb_M\dummySol - \delta\sol\newton{k} \|_{{\bm \Theta}}$ decreases
monotonically as vectors are added to the POD basis.
\\






\subsubsection{Hyper-reduction} \label{sec:hyper}

The LSPG projection approach described in
Section \ref{sec:lspg} is inefficient for nonlinear problems.
This is because the solution of the ROM system requires algebraic
operations that scale with the dimension of the original full order
model $\ndof$. This problem can be circumvented through the use of
hyper-reduction. A number of hyper-reduction approaches have been
proposed, including DEIM \cite{deim}, ``best points" interpolation \cite{bp1, bp2},
collocation \cite{collocation}, gappy POD \cite{gappy} and 
the ECSW method \cite{FarhatChapman:2015, Chapman:2017, FarhatCortial:2014}.
The basic idea behind these approaches is to compute
the residual at some small number of points $q$ with $q \ll \ndof$,
encapsulated in a ``sampling matrix" $\boldsymbol Z$.  This set of
$q$ points is typically referred to as the ``sample mesh". The
``sample mesh" is computed offline using a greedy algorithm
following the notion of gappy data reconstruction \cite{gappy}, as
described in detail in \cite{kevin_gnat} and \cite{carlbergHawaii}.
The LSPG projection approach combined with gappy POD hyper-reduction
is equivalent to the GNAT method \cite{kevin_gnat}. In this method, the nonlinear least-squares
problem \eqref{lspg} is replaced with  
 \begin{equation} \label{gnat}
     \reduced{\sol}_{\text{GNAT}} \in \arg\min_{\dummySol\in\RR{\nrb}}\|(\boldsymbol{ZW})^+\boldsymbol{Z}\res(\reference{\sol} +
     \rb_M\dummySol)\|_2^2,
\end{equation}
where $\boldsymbol{W} \in \RR{\ndof\times\nrb}$ is a 
reduced basis for the residual $\res$ and the `+' symbol
denotes the pseudo-inverse.  It is straightforward to demonstrate that the Gauss--Newton iterations corresponding to \eqref{gnat}
take the form
\begin{gather} \label{gnatit}
	\delta\reduced{\sol}_{\text{GNAT}}\newton{k} \in
	\arg\min_{\dummySol\in\RR{\nrb}}\left| \left|(\boldsymbol{ZW}^+)\boldsymbol{Z} \left[ \jk\rb_M\dummySol +
	\rk\right]\right| \right|_2^2,\quad k=1,\ldots,K_{\text{GNAT}},\\
	\label{gnatup1}\reduced{\sol}_{\text{GNAT}}\newton{k} = \reduced{\sol}_{\text{GNAT}}\newton{k-1} + 
	\lengthk\delta\reduced{\sol}_{\text{GNAT}}\newton{k},\\
	\label{gnatup2}\approximate{\sol}_{\text{GNAT}}\newton{k} = \reference{\sol} + \rb_M\reduced{\sol}_{\text{GNAT}}\newton{k-1},
\end{gather}
where $K_{\text{GNAT}}$ denotes the number of Gauss--Newton iterations.  It is noted that the gappy POD approximation in \eqref{gnat} and \eqref{gnatit} 
aims to approximate the entire residual and Jacobian (via least-squares approximation) rather than simply sub-sample those quantities.



\section{Preconditioned reduced order models}
\label{sec:precond}

Having overviewed our basic ROM workflow, we now introduce the
concept of preconditioned ROMs to the general LSPG formulation
described in Section \ref{sec:lspg}.  Let $\M \in \RR{\ndof\times\ndof}$
denote a non-singular preconditioner matrix.  In the present work, 
a preconditioned LSPG ROM is obtained by inserting $\M$ into the least-squares
optimization problem \eqref{lspg} to yield: 
 \begin{equation} \label{lspg_precond}
     \PrecondPetrovGalerkin{\reduced{\sol}} \in \arg\min_{\dummySol\in\RR{\nrb}}\|\M\res(\reference{\sol} +
     \rb_M\dummySol)\|_2^2,
 \end{equation}
 where $\PrecondPetrovGalerkin{\reduced{\sol}}$ denotes the generalized
 coordinates of the preconditioned Petrov--Galerkin 
 ROM. It is noted that 
 $\M\res(\sol) = \boldsymbol{0}$
 is equivalent to \eqref{fom} provided
 $\M$ is non-singular.  

 Similar to the standard LSPG approach, if \eqref{lspg_precond} is solved 
 using a 
 Gauss--Newton algorithm, one obtains the following sequence of linear least-squares problems:
\begin{gather} \label{lspgit_precond}
    \delta\PrecondPetrovGalerkin{\reduced{\sol}}\newton{k} \in
	\arg\min_{\dummySol\in\RR{\nrb}}\left|\left|\Mk\left(\jk\rb_M\dummySol +
	\rk\right)\right|\right|_2^2,\quad k=1,\ldots,\PetrovGalerkin{K},\\
    \label{lspgup1}\PrecondPetrovGalerkin{\reduced{\sol}}\newton{k} = \PrecondPetrovGalerkin{\reduced{\sol}}\newton{k-1} + \lengthk\delta\PrecondPetrovGalerkin{\reduced{\sol}}\newton{k},\\
    \label{lspgup2}\PrecondPetrovGalerkin{\approximate{\sol}}\newton{k} = \reference{\sol} + \rb_M\PrecondPetrovGalerkin{\reduced{\sol}}\newton{k-1},
\end{gather}
where $\PrecondPetrovGalerkin{K}$ is the number of Gauss--Newton iterations. In \eqref{lspgit_precond}, 
the superscript $(k)$ has been added to $\Mk$ to indicated that, just like the Jacobian $\jk$, the preconditioner
may change in each Gauss--Newton iteration.  We will assume that the sequence of preconditioners
$\Mk$ is non-singular.
The
normal equations corresponding to \eqref{lspgit_precond} are
\begin{equation} \label{lspgitne_prec}
    \PrecondPetrovGalerkin{\jac}\newton{k} \delta\PrecondPetrovGalerkin{\reduced{\sol}}\newton{k} =
    -\PrecondPetrovGalerkin{\res}\newton{k},
    \quad k=1,\ldots,\PrecondPetrovGalerkin{K} \;
\end{equation}
where 
\begin{equation} \label{redjac_precond} \PrecondPetrovGalerkin{\jac}\newton{k}:= \rb_M^T (\jac\newton{k})^T (\M \newton{k})^T\Mk \jac\newton{k}\rb_M, \hspace{0.5cm}
	\PrecondPetrovGalerkin{\res}\newton{k}:= \rb_M^T (\jac \newton{k})^T (\M \newton{k})^T \Mk\res\newton{k}.
\end{equation}

In the remainder of this section, we bring to light and discuss the implications of several important properties
of preconditioned LSPG ROMs constructed via the approach described above.

\newtheorem*{thm1}{Theorem 1}

\begin{thm1}
An ``ideal" preconditioned ROM (a ROM in which $\Mk = \jkinv$ for each 
Gauss--Newton iteration $k$) emulates the projection of the FOM solution increment onto a given
POD basis:
\begin{equation} \label{eq:projsol3}
    \delta\approximate{\sol}\newton{k} = \rb_M \left( \rb_M^T \rb_M \right)^{-1} \rb_M^T \delta\sol\newton{k}.
\end{equation}
	Since \eqref{eq:projsol3} is the upper limit on the ROM accuracy given a basis $\rb_M$, it follows 
that the most accurate ROM possible is realizable through the introduction of preconditioning. \\
\end{thm1}

\begin{proof}
Assume $\jk$ is full rank and let $\Mk = \jkinv$.  
Substituting this $\Mk$ into \eqref{lspgitne_prec} reduces
these normal equations to
\begin{equation} \label{lspgitne_ideal_prec}
    \rb_M^T\rb_M\delta\PrecondPetrovGalerkin{\reduced{\sol}}\newton{k}  = -\rb_M^T\jkinv\rk, \quad k=1,\ldots,\PrecondPetrovGalerkin{K}
\;.
\end{equation}
Recognizing that the right-hand side of \eqref{lspgitne_ideal_prec}
contains the FOM solution increment $\delta\sol\newton{k} = -\jkinv\rk$,
\eqref{lspgitne_ideal_prec} is equivalent to
\begin{equation} \label{lspgitne_ideal_prec2}
    \delta\PrecondPetrovGalerkin{\reduced{\sol}}\newton{k}  = (\rb_M^T\rb_M)^{-1}\rb_M^T\delta\sol\newton{k}, \quad
    k=1,\ldots,\PrecondPetrovGalerkin{K} \;,
\end{equation}
so that, from \eqref{eq:subspace}, $\delta\approximate{\sol}\newton{k}$ is given by \eqref{eq:projsol3}. 

\end{proof}

\newtheorem*{thm2}{Theorem 2}

\begin{thm2}
Introducing preconditioning into the LSPG formulation
as in \eqref{lspgit_precond} has the effect of modifying the norm ${\bm \Theta}$ 
in which the solution increment
minimization is performed \eqref{opt}.  For a generic non-singular preconditioner matrix $\Mk$,
the LSPG ROM solution procedure can be can be viewed as minimizing at each iteration the error in the 
search direction 
\eqref{opt} in the norm given by ${\bm \Theta} = (\jk)^T (\Mk)^{T} \Mk \jk$. For an ideal 
	preconditioner ($\Mk = \jkinv$), ${\bm \Theta} = \I$, which leads to an $l^2$-optimal projection of the 
solution increment $\delta\sol\newton{k}$ onto the reduced basis $\rb_M$: 
\begin{equation} \label{opt_ideal_precond}
    \delta\PrecondPetrovGalerkin{\reduced{\sol}}\newton{k} \in \arg
\min_{\dummySol\in\RR{\nrb}}\|\rb_M\dummySol - \delta\sol\newton{k}\|_2^2
\end{equation}
\end{thm2}

\begin{proof}
 It is straightforward to see that
\begin{equation} \label{pf2}
\|\Mk\left(\jk\rb_M\dummySol +
        \rk\right)\|_2^2 = \|\Mk \jk \left(\rb_M\dummySol +
        \jkinv \rk\right)\|_2^2 = \| \rb_M\dummySol +
	\jkinv \rk\|_{{\bm \Theta}}^2,
\end{equation}
where ${\bm \Theta} = (\jk)^T (\Mk)^{T} \Mk \jk$.  The claim follows immediately.  
Substituting $\Mk = \jkinv$ and $\delta\sol\newton{k} = -\jkinv\rk$ into \eqref{pf2} 
gives \eqref{opt_ideal_precond}.  
\end{proof}

While selecting $\Mk = \jkinv$ is infeasible in practice, an important corollary of Theorem 1 
is that by selecting a preconditioner matrix that approximates $\jkinv$, it is possible to 
improve the accuracy of a given LSPG ROM, with the ROM solution approaching the most accurate 
ROM solution possible (for a given basis $\rb_M$) as $\Mk \to \jkinv$.  Moreover, if $\Mk$ is selected 
such that $\text{cond}\left(\Mk \jk\right) < \text{cond}\left(\jk\right)$ (where $\text{cond}(\boldsymbol{A})$ denotes
the condition number of a matrix $\boldsymbol{A}$), it is possible to 
improve the condition numbers of the reduced Jacobians \eqref{redjac_precond} arising in the 
Gauss--Newton iteration process.

We end this section by remarking that there is an additional benefit of
introducing a preconditioner into the LSPG formulation \eqref{lspg} 
as in \eqref{lspg_precond}.   
It has been observed \cite{kyle, Blonigan:2022, Parish:2022} that minimizing the raw
(unweighted) residual $\res$ can be problematic for systems of PDEs
where different variables have drastically different magnitudes,
which happens frequently in various applications where dimensional
equations are employed.  Residual components corresponding to
certain variables can be very large compared to the residual
components corresponding to other variables, thereby biasing the
minimization procedure.  As shown in \cite{kyle, Blonigan:2022, Parish:2022}, scaling the ROM
residual to get all the equations to be roughly the same order of
magnitude can remedy this problem.  From \eqref{lspg_precond}, it is can 
be seen that the introduction of
a preconditioner $\M$ can have the same effect.






\section{Implementation in the \albany{} multi-physics finite element
code} \label{sec:impl}

The preconditioned LSPG ROMs described herein have been
implemented with a version\footnote{Available on {\tt github} at the
following URL:
\url{https://github.com/sandialabs/Albany/releases/tag/MOR_support_end}.}
of the \albany{} code base \cite{albany_repo}, described succinctly
in this section. \albany{} is an open-source\footnote{\albany{}
is available on GitHub:
\url{https://github.com/sandialabs/Albany}.} C++
object-oriented, parallel, unstructured-grid, implicit finite
element code for solving general PDEs, 
developed using the ``Agile Components" code
development strategy \cite{AgileComponents} with
mature modular libraries from the
\trilinos{} \cite{trilinos} project\footnote{\trilinos{} is
available on GitHub: \url{https://github.com/trilinos/Trilinos}.}.
Over the
years, \albany{} has hosted a number of science and engineering
applications, including the \aeras{} global atmosphere code
\cite{Aeras:2015}, the Albany Land-Ice (\felix{}) \cite{Tezaur:2015}
ice sheet model solver, the Quantum Computer Aided Design (\qcad{})
\cite{QCAD:2013} simulator, the \ace{} thermo-mechanical terrestrial model of Arctic coastal erosion \cite{ACE_paper}, and the Laboratory for Computational
Mechanics (\lcm{}) \cite{LCM:2013, schwarz, Mota:2022} research code. This last project
comprises \albany{} \lcm{} and is specifically targeted at solid
mechanics applications, such as the ones considered in this paper. A
more detailed description of \albany{}, including a detailed
description of its underlying design and the physics implemented
therein, can be found in \cite{albany}.

In the subsections below, we describe the main ingredients enabling
the construction of preconditioned LSPG ROMs within the
\albany{} code base.  We present also a special partitioned/``blocking vector" approach for applying ROM 
Dirichlet boundary conditions within this code base, which does not remove the constrained dofs from 
the global finite element system prior to performing its 
numerical solution.  This implementation is extensible to other finite element codes 
with similar boundary condition treatment.  


\subsection{The model reduction workflow in \albany{}} \label{sec:workflow}

As described in Section \ref{sec:rom}, the basic workflow in
building a ROM consists of three steps: calculation of
a reduced basis, projection of the governing equations onto the
reduced basis, and calculation of a sample mesh using
hyper-reduction.  The implementation of this workflow within
\albany{} for a generic nonlinear problem is demonstrated in
Algorithm \ref{alg:albany_mor_workflow}.

\begin{algorithm}[h!]
    \begin{small}
\SetAlgoLined 
 Given: {\tt Albany}, {\tt AlbanyRBGen}, {\tt AlbanyMeshSample} executables,
 and corresponding input files: fomInput.xml, reducedBasisInput.xml, meshReduceInput.xml, romInput.xml  \\
 \vspace{0.2cm}
   \# Run the full order model to get the snapshots \\
   ./{\tt Albany} fomInput.xml \\
   \vspace{0.2cm}
   \# Compute reduced basis (RBGen required) \\
./{\tt AlbanyRBGen} reducedBasisInput.xml \\
 \vspace{0.2cm}
\# Compute sample nodes and reduced mesh\\
./{\tt AlbanyMeshSample} meshReduceInput.xml\\
 \vspace{0.2cm}
\# Run the ROM with hyper-reduction (sample mesh) \\
 ./{\tt Albany} romInput.xml
    \caption{\albany{} model reduction workflow.  Example input files can be
    found in the GitHub repository corresponding to this version of \albany{} \cite{albany_repo}, in the directory
    Albany-MOR\_support\_ends/tests/small/MOR/MOR\_MechanicalCube.}
    \label{alg:albany_mor_workflow}
    \end{small}
\end{algorithm}

First, the main {\tt Albany} executable is run in FOM mode to
collect the relevant snapshots for a given problem, which are
written to an output \exodus{} file.  Next, the reduced basis
and sample mesh are calculated in \albany{} as a pre-processing step
by utilities known as {\tt AlbanyRBGen} and {\tt AlbanyMeshReduce},
respectively, using data collected from a high-fidelity simulation.
The {\tt AlbanyRBGen} utility employs the \rbgen{} module within the
\anasazi{} package of \trilinos{} \cite{anasazi} to calculate a
reduced basis $\rb_M$ using the POD algorithm outlined in Section
\ref{sec:pod}.  While it is possible to use {\tt AlbanyRBGen} to generate 
a scalar POD basis (i.e., a separate POD basis built for each of the unknowns
to a vector-valued problem), the utility was designed to construct 
vector POD bases, and vector POD bases were employed in
the numerical experiments in Section \ref{sec:numerical}.
We note that {\tt AlbanyRBGen} does not orthonormalize the POD modes with 
respect to the mass matrix for dynamic solid mechanics problems,
meaning the resulting bases are $l^2$ (as opposed to $L^2(\Omega)$) orthogonal.

Hyper-reduction in \albany{} is also handled as a pre-processing
step using another \albany{} utility, {\tt AlbanyMeshSample}, which
creates the ``sample mesh'' by using either the collocation or the 
gappy POD approach \cite{CarlbergGappy, kevin_gnat} applied to the snapshot data (see
Section \ref{sec:hyper}).  Following the calculation of $\rb_M$ and
sample mesh, the main {\tt Albany} executable is run again, this
time in ROM mode, with $\rb_M$ and the sample mesh as inputs. During this
step, \albany{} solves the optimization problem \eqref{lspg} (or \eqref{lspg_precond} 
if preconditioning is used) using the Gauss--Newton method,
and evaluates the ROM solution.  
The
details of the implementation of the projection step of the model
reduction procedure in \albany{} are described below in Section
\ref{sec:wrapper}.
The linear least-squares problem
arising within each Gauss--Newton step is solved 
using the normal equations 
approach (\eqref{lspgitne} or  \eqref{lspgitne_prec}).  
Solutions to these linear systems are obtained 
using Krylov iterative methods
(Conjugate Gradient or GMRES), implemented within the \aztec{} library \cite{aztec} in 
\trilinos{}.



\subsection{Projection via the {\tt Albany::ModelEvaluator} wrapper}
\label{sec:wrapper}

The implementation of the projection step of the model reduction
procedure within \albany{} relies heavily on the {\tt
EpetraExt::ModelEvaluator} nonlinear model abstraction available
within this code \cite{epetraext_me}, whose concrete implementation
in Albany is referred to hereafter as the {\tt
Albany::ModelEvaluator}. The {\tt Albany::ModelEvaluator} provides
an extensible interface to the underlying application code
(\albany{}) for a variety of solver and analysis algorithms within
\trilinos{}.  Suppose we are solving a nonlinear algebraic system
that can be written as
\begin{equation} \label{eq:resid}
	\res(\sol, {\bm \mu}) = \boldsymbol{0},
\end{equation}
where $\res$ is the residual, $\sol$ is the solution and ${\bm \mu}$ is a vector of parameters,
using the basic
Newton method.  The Newton solution procedure requires querying the
application for various quantities, including the residual
\eqref{eq:resid}, and the Jacobian, defined as $d\res/d\sol$, calculated
in \albany{} using automatic differentiation within the \sacado{}
library. The objective of the model evaluator interface is to
provide these quantities to \trilinos{} solvers through the {\tt
Albany::ModelEvaluator} interface, which is agnostic to the physics,
model and data structure that are needed to support the relevant
matrix and vector abstractions within the solver. The {\tt
Albany::ModelEvaluator} has utility beyond this simple nonlinear
solver example: it works in a similar fashion to perform time
integration, continuation, sensitivity analysis, stability analysis,
optimization and uncertainty quantification.  For these more
sophisticated analysis cases, it returns additional quantities
besides $\res$ and $d\res/d\sol$, e.g., the time-derivative of $\sol$
(denoted by $\dot{\sol}$), sensitivities of $\res$ with respect to
parameters ${\bm \mu}$ (denoted by $d\res/d{\bm \mu}$), etc.

The model reduction capabilities in \albany{} are implemented by creating a
wrapper around the {\tt Albany::ModelEvaluator} class within
this code, known as the {\tt MOR::ReducedOrderModelEvaluator}.  This
class is only activated when reduced order model options are enabled
in an \albany{} input file, and takes the objects returned by the
{\tt Albany::ModelEvaluator}, e.g., the residual $\res$
\eqref{eq:resid} and Jacobian $d\res/d\sol$, together with a reduced basis
$\rb_M$, and returns the reduced residual and Jacobian defined in
\eqref{lspgit} for the LSPG method.
The ROM formulation within \albany{} is generic
enough to create ROMs for various solvers and physics sets enabled
in \albany{}; herein we focus on quasi-static mechanical and
thermo-mechanical problems from the \lcm{} suite within \albany{}
solved using the \loca{} continuation package within
\trilinos{} \cite{loca}. 
There are currently two projection options implemented in Albany for
the type of reduced order model, Galerkin and LSPG projection.
Herein, attention is restricted to the latter approach.

\subsection{Implementation of boundary conditions} \label{sec:bcs}


In model reduction, it is customary to work with equations governing the \textit{unconstrained}
dofs, i.e., the dofs modulo any Dirichlet boundary conditions imposed.
While in the standard finite element method formulation, constrained dofs are
typically removed from the nonlinear system being solved, it is not uncommon for finite element
codes to retain these dofs.  This approach is taken in the boundary condition
implementation within \albany{}.  Specifically, in place of the system \eqref{fom},
the governing equations in Albany take the form
\begin{equation} \label{eq:fom_bcs}
  \bar{\res}(\bar{\sol}) = \boldsymbol{0}.
\end{equation}
In \eqref{eq:fom_bcs}, $\bar{\sol} \in \RR{\bar{\ndof}}$ and $\bar{\res}: \RR{\bar{\ndof}} \to \RR{\bar{\ndof}}$
denote the ``extended" state and residual, respectively, with $\bar{\ndof} = \ndof + \ndof_{\text{DBC}}$, where
$\ndof_{\text{DBC}}$ and $N$ are the number of constrained (Dirichlet)  and unconstrained dofs in the system, respectively. The extended variables
can be rearranged to take the form:
\begin{equation}\label{eq:bc_defs}
  \bar{\sol} := \begin{varray} 
                  \sol \\
                  \sol_{\text{DBC}}
                \end{varray}
  \hspace{0.5cm}
  \bar{\res} := \begin{varray} 
                  \res \\
                  \res_{\text{DBC}}
                \end{varray},
\end{equation}
where $\sol_{\text{DBC}} \in \RR{\ndof_{\text{DBC}}}$ denote the constrained dofs that extend the solution $\sol \in \RR{\ndof}$.
Assuming \eqref{eq:fom_bcs} is solved using a globalized Newton method as in \eqref{fom_newton},
  \albany{} solves the following sequence of problems
\begin{equation}\label{fom_newton_bcs}
  \begin{aligned}
    \bar{\jac}\newton{k}\delta\bar{\sol}\newton{k} &= -\bar{\res}\newton{k}, \quad k=1,\ldots,K,\\
    \bar{\sol}\newton{k} &= \bar{\sol}\newton{k-1} + \lengthk\delta\bar{\sol}\newton{k},
  \end{aligned}
\end{equation}
where
\begin{equation} \label{dbc_jac}
  \bar{\jac}\newton{k} =
  \begin{marray} 
    \jk            & \boldsymbol{0} \\
    \boldsymbol{0} & \bar{\jac}\newton{k}_{\text{DBC}}
  \end{marray},
  \hspace{0.5cm}
  \bar{\res}\newton{k} := \begin{varray} 
                            \res\newton{k} \\
                            \boldsymbol{0}
                          \end{varray}.
\end{equation}
In \eqref{dbc_jac}, the entries of $\bar{\jac}\newton{k}_{\text{DBC}} \in \RR{\ndof_{\text{DBC}} \times
\ndof_{\text{DBC}}}$ are defined
as follows: $[\bar{\jac}\newton{k}_{\text{DBC}}]_{ij}\defeq \left[\frac{\partial \bar{\res}}{\partial \sol_{\text{DBC}}}\left(\sol_{\text{DBC}}\newton{k}\right)\right]_{ij} \delta_{ij}$ for $i = 1, ..., \ndof_{\text{DBC}}$ and $j = 1, ..., \ndof_{\text{DBC}}$, 
where $\delta_{ij}$ denotes the Kronecker delta.  In other words, $\bar{\jac}\newton{k}_{\text{DBC}}$
is the diagonal part of the Dirichlet boundary condition (DBC) restriction of the FOM Jacobian $\jk \in \RR{\ndof\times\ndof}$. 
The remaining variables appearing in \eqref{fom_newton_bcs} were defined earlier in \cref{sec:fom}.
Remark that the implementation \eqref{fom_newton_bcs} maintains symmetry of the original system:
  since $\delta \sol_{\text{DBC}}\newton{k} = \boldsymbol{0}$ for all rows corresponding to Dirichlet boundary conditions,
  these rows are fully decoupled from all the other rows,
  and the corresponding columns can be trivially
zeroed.

A consequence of the above discussion is that care must be taken to ensure boundary conditions are properly
applied in ROMs constructed within the \albany{} code.  When computing the reduced basis $\rb_M$
using \albany{}'s {\tt AlbanyRBGen} utility (step 2 of Algorithm \ref{alg:albany_mor_workflow}),
all constrained dofs (dofs for which Dirichlet BCs are imposed)
are excluded from the snapshot matrix on which POD is performed.
The resulting basis modes are then augmented with a set of modes, denoted 
by $\rb_{\text{DBC}} \in \RR{\ndof_{\text{DBC}}\times\nrb_{\text{DBC}}}$, defined on a so-called \textit{blocking region}
of the domain when performing the projection-based reduction.  
These modes, $\rb_{\text{DBC}}$, are used as modes for the Dirichlet dofs,
  with common dofs grouped together into \textit{blocks}.
For instance, individual components ($x$, $y$, $z$) of a certain nodeset might be grouped into a block.
Specifically,
\begin{equation} \label{eq:dfns_gal_bcs}
  \rb_{\bar{M}} := \begin{marray} 
                   \rb_M & \boldsymbol{0} \\ 
                   \boldsymbol{0} & \rb_{\text{DBC}}
                 \end{marray},
\end{equation}
where the columns of $\rb_{\text{DBC}} \in \RR{\ndof_{\text{DBC}}\times\nrb_{\text{DBC}}}$ are constructed from the (normalized) summation of the canonical unit vectors of each \textit{block} of common Dirichlet dofs, 
i.e. $\rb_{\text{DBC}} := [\uvec_1, \uvec_2, \ldots, \uvec_{\nrb_{\text{DBC}}}]$,
  where $\uvec_i := \bvec_i/\norm{\bvec_i}$, $\bvec_i = \sum_{j \in C_i} \evec_j$,
  and $C_i$ is the set of dofs for each block of dofs comprising the blocking region as described above.

Performing an LSPG projection amounts to
solving the following system
\begin{equation} \label{eq:lspg_bcs}
  \lb_{\bar{M}}\newton{k}(\reference{\bar{\sol}}
                        +
                        \rb_{\bar{M}} \PetrovGalerkin{\bar{\reduced{\sol}}})^T
  \bar{\res}(\reference{\bar{\sol}}
             +
             \rb_{\bar{M}}\PetrovGalerkin{\bar{\reduced{\sol}}})
  = 0,
\end{equation}
where
\begin{equation} \label{eq:dfns_lspg_bcs}
  \lb_{\bar{M}}\newton{k}(\sol) := \bar{\jac}\newton{k} \lb_{\bar{M}}
                               = \begin{marray} 
                                   \jk\rb_M & \boldsymbol{0} \\
                                   \boldsymbol{0} & \bar{\jac}\newton{k}_{\text{DBC}}\rb_{\text{DBC}}
                                 \end{marray},
  \hspace{0.5cm}
  \PetrovGalerkin{\bar{\reduced{\sol}}} := \begin{varray} 
                                             \PetrovGalerkin{\reduced{\sol}} \\
                                             \reduced{\sol}_{\text{DBC}}
                                           \end{varray}.
\end{equation}
and $\PetrovGalerkin{\reduced{\sol}} \in \RR{\nrb}$ and $\reduced{\sol}_{\text{DBC}} \in \RR{\nrb_{\text{DBC}}}$ are the ROM solution 
(equations \eqref{eq:subspace} and \eqref{lspg}) and \textit{blocking solution}, respectively.


\subsection{Preconditioned ROMs in \albany{}} \label{sec:precond_eval}

To precondition the LSPG ROMs constructed by the \albany{} code base, 
we rely on the \ifpack{} suite of preconditioners \cite{ifpack-guide} within \trilinos{}. 
The specific preconditioners tested herein are summarized in Table \ref{tab:preconds}.  
In all \ifpack{} preconditioned cases, 
the effect of preconditioning is to pre-multiply the terms $\jk\rb$ and $\rk$ 
in \eqref{lspgit}  by an approximation
of the FOM Jacobian inverse, $\Mk \approx \jkinv$.
We limit our attention to three \ifpack{} preconditioners: a Jacobi 
preconditioner, a Gauss-Seidel preconditioner and an incomplete LU (ILU) preconditioner with a level-of-fill 
of 1.
The Jacobi 
preconditioner approximates $\jkinv$ by an inverse of the diagonal 
of the Jacobian $\jk$, whereas the Gauss-Seidel preconditioner 
is constructed by inverting the elements of the upper triangular part of $\jk$.  
Our ILU preconditioner is another simple preconditioner,
obtained by performing an incomplete LU 
factorization of $\jk$ and allowing only one non-zero diagonal above
and one non-zero diagonal below the main diagonal.  
For a detailed description of these preconditioners and their implementations
within \ifpack{}, 
the interested reader is referred to \cite{ifpack-guide}.
By considering several preconditioners, 
we are able to study numerically how the ROM error changes
as the preconditioner is $\Mk$ is improved.
In addition to the preconditioners provided by \ifpack{}, 
we consider also, when computationally feasible to calculate, the ideal preconditioned case, 
in which $\Mk$ is calculated by computing the action of the inverse of the Jacobian matrix $\jkinv$.
As shown earlier, in Theorem 1, calculating $\Mk = \jkinv$ is equivalent to directly computing
the FOM solution increment $\delta \sol^k $ by solving $\jk\delta \sol^k = -\rk$,
then projecting that solution increment onto the ROM space to obtain the
ROM solution increment $\delta\PetrovGalerkin{\reduced{\sol}}\newton{k}$,
as in \eqref{eq:projsol3}. Since this approach solves the FOM system at every
iteration instead of a modified ROM system as with the other approaches,
the dominant source of error with this approach will be in the projection,
i.e., how much the basis $\rb_M$ has been truncated, and we expect the errors with
this approach to be very low.  


\begin{table}[H]
    \centering
\caption{Summary of preconditioners evaluated}
    \begin{tabular}{c|c}
Preconditioner Name & Description of Preconditioner \\ 
\hline 
	    {\tt Jacobi} & Jacobi\\
	    {\tt Gauss-Seidel} & Gauss-Seidel\\
	    {\tt ILU} & Incomplete LU factorization with 1 level-of-fill \\
	    {\tt Ideal} & $\jkinv$ (equivalent to projected solution increment \eqref{eq:projsol3})\\
\end{tabular}
\label{tab:preconds}
\end{table} 

\section{Numerical examples}
\label{sec:numerical}





\subsection{Governing equations} \label{sec:eqns}

The predictive ability 
of the \albany{} model reduction capabilities described in Section
\ref{sec:impl} are evaluated on three problems based on the quasi-static mechanics
and thermo-mechanics PDEs implemented within \albany{}'s \lcm{}
suite.  The last of these cases is a problem of practical scale 
and involves a realistic geometry.
Before presenting our test cases, we briefly summarize the PDEs being solved.

\subsubsection{Mechanical formulation} \label{sec:mechanical}

Consider a body defined by an open set $\Omega \subset \RR{3}$
undergoing a motion described by the mapping $~x = ~\varphi(~X):
\Omega \times  \RR{3}$,
where $~X \in \Omega$.  Let $~F :=
\Grad ~\varphi$ be the deformation gradient.  Let also $\rho_0 ~B:
\Omega \to \RR{3}$ be the body force, with $\rho_0$ the mass density
in the reference configuration.  It is straightforward to show
\cite{schwarz} that the Euler-Lagrange equations for a canonical
quasi-static mechanical problem take the form:
\begin{equation} \label{eq:euler-lagrange}
  \Div ~P + \rho_0 ~B = ~0, \hspace{0.2cm} \text{in } \Omega.
  \end{equation}
Here, $~P = \partial W / \partial ~F$ denotes the first
Piola-Kirchhoff stress and, assuming an elastic material model\footnote{The \albany{}
\lcm{} code contains a wide range of constitutive models for solid mechanics, and includes models
for both elastic and plastic materials.  Since the results presented herein assume 
an elastic material, we restrict our presentation to the elastic case.  For plastic 
materials, the Helmholtz free-energy density is also a function of the internal variables, denoted by $~Z$: $W=W(~F, ~Z)$.}, $W = W(~F)$ is the Helmholtz
free-energy density.
Assume that the boundary of the body is $\Gamma =
\partial_{~\varphi} \Omega \cup \partial_{~T} \Omega$,
where $\partial_{~\varphi} \Omega$ is a prescribed position
boundary, $\partial_{~T} \Omega$ is a prescribed traction boundary,
and $\partial_{~\varphi} \Omega \cap \partial_{~T} \Omega =
\emptyset$. Denoting the prescribed boundary positions or Dirichlet
boundary conditions by $~\chi :
\partial_{~\varphi} \Omega \to \RR{3}$, and the prescribed boundary
tractions or Neumann boundary conditions are $~T: \partial_{~T}
\Omega  \to \RR{3}$, \eqref{eq:euler-lagrange} is supplemented
by the following constraints:
\begin{equation} \label{eq:initial-conditions}
    \begin{array}{c}
         ~\varphi(~X)
       =
      ~\chi
       \text{  on  }
      \partial_{~\varphi} \Omega,
        \hspace{0.2cm} ~P ~N
       =
      ~T
       \text{  on  }
      \partial_{~T} \Omega.
    \end{array}
\end{equation}
Embedded within the Helmholtz free-energy density $W(~F)$ is a
constitutive or material model.  
For the numerical examples given
here, we employ a Neohookean constitutive model extended to the compressible regime, which is consistent
with a nonlinear elastic material \cite{Holzapfel}.  Following common practice in 
solid mechanics, we decompose $W(~F)$ into its volumetric and deviatoric components: 
\begin{equation} \label{eq:Wdecomp}
	W(~F) = W^{\text{vol}}(~F) + W^{\text{dev}}(~F).
\end{equation}
For a Neohookean material, the components comprising $W(~F)$ in \eqref{eq:Wdecomp} are given by: 
\begin{equation} \label{eq:neohookean} 
	W^{\text{vol}}(~F):=\frac{1}{2}\kappa \left[\frac{1}{2}(J^2 - 1) - \log J \right], \hspace{0.5cm} 
	W^{\text{dev}}(~F):=\frac{1}{2}\mu \left[J^{-2/3} \text{trace} (\boldsymbol{b}) - 3\right].
\end{equation}
In \eqref{eq:neohookean}, $\kappa$ denotes the bulk modulus, given by 
\begin{equation} \label{eq:bulk}
	\kappa = \frac{E}{3(1-2\nu) },
\end{equation}
where $E$ is the elastic modulus and $\nu$ is Poisson's ration; $\mu$ denotes the 
shear modulus, given by 
\begin{equation} \label{eq:shear}
	\mu = \frac{E}{2(1+\nu)}.
\end{equation}
$J:=\text{det}(~F)$ is the determinant of the deformation gradient, and 
\begin{equation} \label{eq:b}
	\boldsymbol{b}:=~F_M ~F_M^T,
\end{equation}
is the left Cauchy-Green deformation, tensor, 
where $~F_M$ denotes the mechanical part of the deformation gradient.  For a 
pure mechanics problem, $~F_M = ~F$. The mechanical equations \eqref{eq:euler-lagrange} are 
solved for the displacements, the primary unknowns in these equations.

\subsubsection{Thermo-mechanical formulation}
\label{sec:thermo-mechanical}

In the case of a quasi-static thermo-mechanical problem, the
mechanics equations \eqref{eq:euler-lagrange} are augmented with the
following generic steady-state heat conduction equation:
\begin{equation} \label{temp}
\nabla \cdot ~q(T) = ~f, \hspace{0.2cm} \text{in } \Omega,
\end{equation}
where $~q$ denotes the heat flux, $T$ denotes the temperature and $~f$
is a source term. Equation \eqref{temp} can include boundary conditions
of both the Dirichlet and Neumann type, e.g.,
\begin{equation}
 T(~X) = T_p  \text{  on  }
      \partial_{~\varphi} \Omega, \hspace{0.2cm}~q(T) = ~q_p
       \text{  on  }
      \partial_{~T} \Omega,
\end{equation}
where $T_p$ and $~q_p$ are prescribed temperature and heat flux
values, respectively.  In this case, a temperature dependence is introduced 
into the Helmholtz free-energy density \eqref{eq:neohookean}: $W = W(~F, T)$.  For the Neohookean 
formulation considered herein, temperature effects are introduced in the form 
of thermal expansion.  More specifically, the mechanical part of the deformation 
gradient appearing in \eqref{eq:b} now takes the form:
\begin{equation}
	~F_M := ~F ~F_T^{-1},
\end{equation}
where 
\begin{equation} \label{eq:FT}
	~F_T:= \exp\left[\alpha(T - T_{\text{ref}})\right] ~I 
\end{equation}
is the thermal part of the deformation gradient, with $~I$ denoting the identity tensor. In \eqref{eq:FT}, 
$\alpha$ is the thermal expansion coefficient and $T_{\text{ref}}$ is the 
reference temperature.
For the
problems considered here, the heat flux in \eqref{temp} takes the
form:
\begin{equation} \label{heat_flux}
	~q(T) = \boldsymbol{K} \cdot \nabla T,
\end{equation}
where $\boldsymbol{K}$ is the thermal diffusivity tensor.

\subsection{Evaluation criteria} \label{sec:eval}

Recall that we denote the FOM solution as $\sol$ and the approximate ROM solution as $\approximate{\sol} = \reference{\sol} + \rb_M\reduced{\sol}$ 
\eqref{eq:subspace}.
Let $\sol_i$ and $\approximate{\sol}_i$ denote the FOM and ROM solutions, respectively, corresponding to time $t_i$, for $i = 1, ..., P$, where $P$ 
denotes the number of snapshots collected.  
Building on these definitions, we define the \textit{global relative error} in the approximate ROM solution as:
\begin{equation}
    \label{pterr}
	\mathcal{\epsilon} :=    \frac{\displaystyle\sum_{i=0}^{P}\left\|\sol_i-\approximate{\sol}\right\|_2}{\displaystyle\sum_{i=0}^{P}\left\|\sol_i\right\|_2}.
\end{equation} 
For mechanical problems, $\sol$ and $\approximate{\sol}$ represent the displacement field solutions, 
whereas for thermo-mechanical problems, these vectors include also the temperature field.
In the subsequent sections, the proposed preconditioned LSPG ROMs are evaluated using the metric \eqref{pterr}.  We emphasize that the evaluations
performed herein are in the \textit{predictive} regime.  That is, we calculate errors by running
the proposed ROMs with \textit{different} 
parameter values than the parameter values used to generate the training data from which the ROMs were constructed.  More details on the training/testing
for each problem considered are provided in the subsections below.
Following the usual convention in solid mechanics, the FOMs and ROMs for all three test cases considered
are run dimensionally.

In addition to reporting the global relative error for each of the ROMs evaluated, we also report the total CPU-time of 
the online stage of the model reduction procedure for each ROM.  The specifications of the run-time environment 
for each benchmark problem are described in the following subsections.
As the purpose of this work is to present a preliminary
study examining the viability of preconditioning within the LSPG framework,
hyper-reduction, which would potentially introduce another source of error, was
not employed in our study.  Evaluating the performance of LSPG ROMs with 
preconditioning and hyper-reduction will be the subject of future work.  It is noted that 
special hyper-reduction techniques, e.g., the ECSW method
\cite{FarhatChapman:2015, Chapman:2017, FarhatCortial:2014} of Chapman, Farhat \textit{et al.} and the structure-preserving/``matrix gappy POD" 
method of Carlberg \textit{et al.} \cite{CarlbergAIAA:2012, Carlberg:2015}, are needed to preserve the Lagrangian structure of ROMs 
for solid mechanics problems such as the ones discussed herein.  These are not available within the 
\albany{} code base at the present time.

Lastly, for each problem considered, we report also either the average condition number of all reduced Jacobians
obtained during the Gauss--Newton iteration process or the total number of nonlinear iterations required 
for ROM convergence.  The latter metric is adopted for the larger thermo-mechanical pressure 
vessel problem (Section \ref{sec:vessel}), as estimating reduced Jacobian condition numbers is not feasible
for a problem of that size.



\subsection{Mechanical and thermo-mechanical beam} \label{sec:beam}

Our first example involves a simple three-dimensional (3D) beam geometry $\Omega$ of size 
$0.16\times0.016\times0.032\len$ centered at $(0,0,0)$, discretized using 320 hexahedral elements, 
which gives rise to 525 nodes (Figure \ref{fig:brickGeometry}).  
The geometry $\Omega$ is 
divided into five material blocks, as shown in in Figure \ref{fig:brickGeometry}(a).
The purpose of having 
several material blocks is to enable the specification of different material models and/or material parameters
in different parts of the domain.  As mentioned earlier, we consider a nonlinear elastic material model known as 
the Neohookean model, which is prescribed in all the blocks comprising $\Omega$; however, as discussed below, different 
material parameters are specified in different blocks.
\begin{figure}[h!tb]
 \centering
   \begin{tabular}{cc}
	   \includegraphics[width=0.4\textwidth]{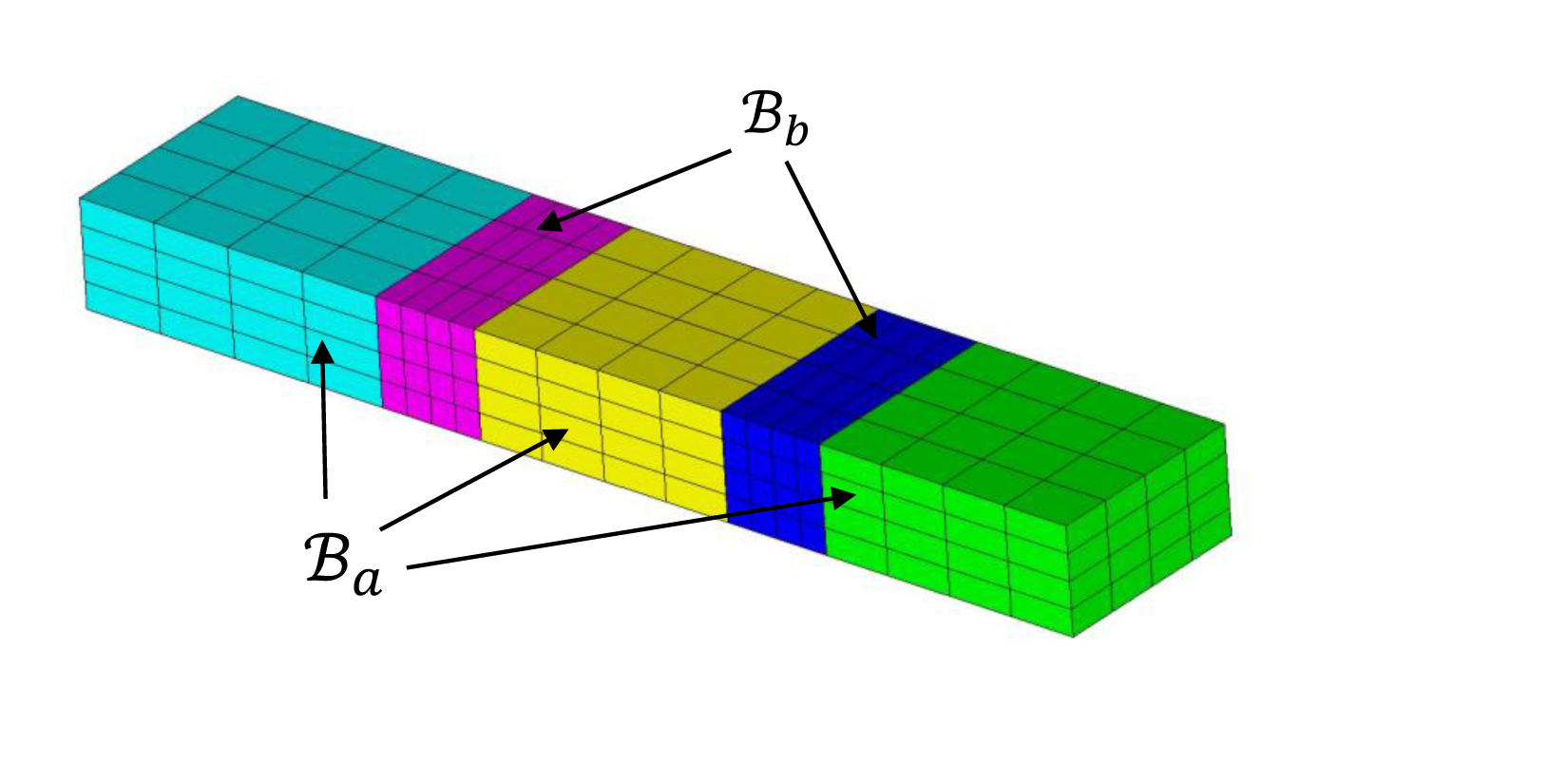}&
   \includegraphics[width=0.4\textwidth]{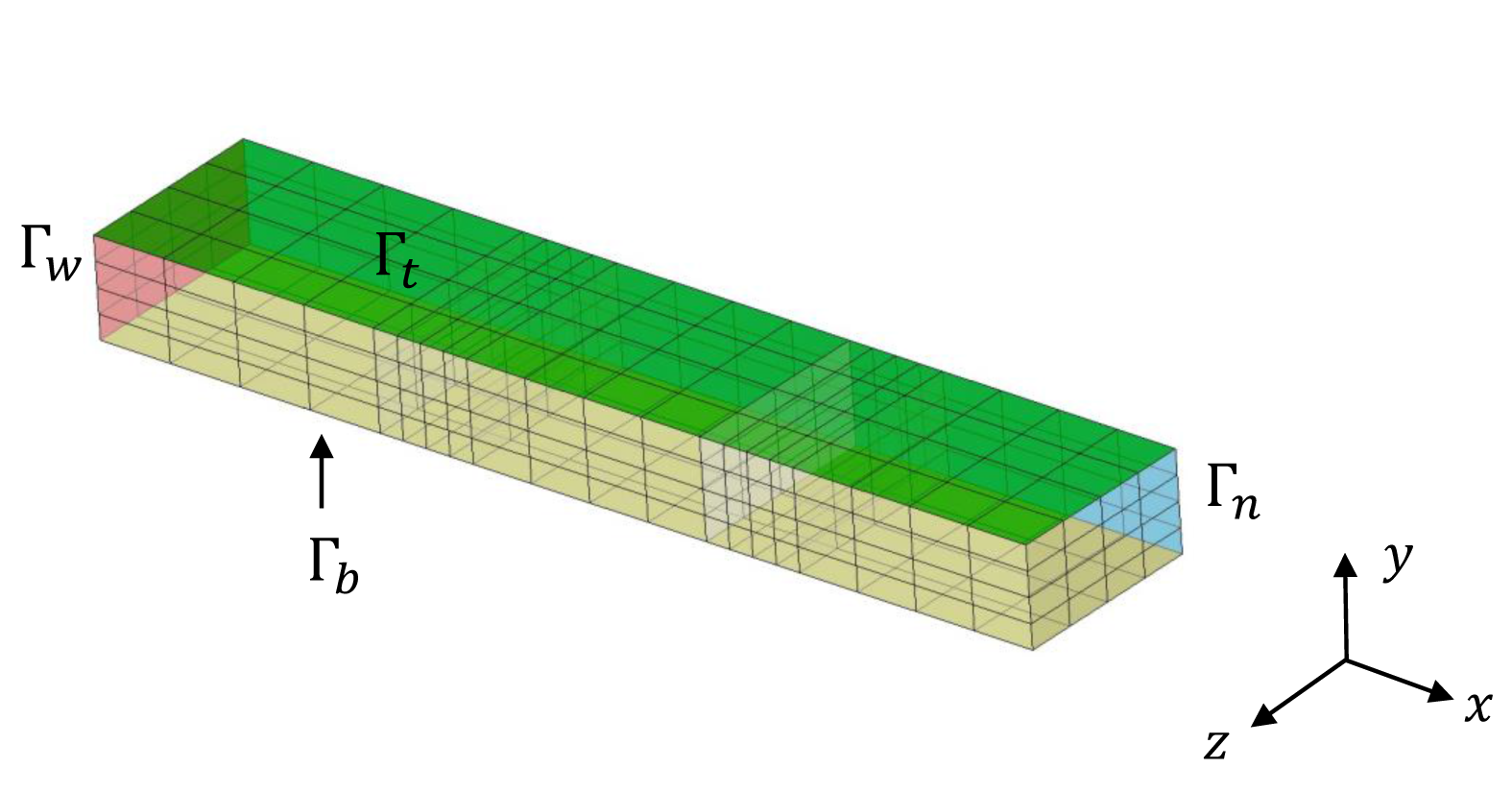} \\
	   (a) Discretization of $\Omega$ & (b) $\Omega$'s boundaries  \\
  \end{tabular}
	\caption{Computational domain for the beam problem.  In subfigure (a), different colors indicate different material blocks.
	Subfigure (b) distinguishes in color the different boundaries on which boundary conditions are prescribed.}
  \label{fig:brickGeometry}
\end{figure}

In order to complete the problem definition, it is necessary to specify well-posed boundary conditions 
on our beam geometry $\Omega$.  We prescribe boundary conditions on four of 
$\Omega$'s boundaries, shown in Figure \ref{fig:brickGeometry}(b) and denoted by $\Gamma_w$, $\Gamma_b$, $\Gamma_t$ and $\Gamma_n$.
$\Gamma_w$ corresponds to $x=-0.08$ and is shown in red; $\Gamma_b$ corresponds to $y=-0.008$ and is shown in yellow; 
$\Gamma_n$ corresponds to $z = -0.016$ and is shown in blue; $\Gamma_t$ corresponds to $y = 0.008$ and is shown in green.
In both versions of our problem formulation (mechanical and thermo-mechanical; Sections \ref{sec:beam_mech} and \ref{sec:beam_thermomech}, respectively), $\Gamma_w$ is 
fixed in the $x$-direction, $\Gamma_n$ is fixed in the $z$-direction and $\Gamma_t$ is fixed in the $y$-direction.  
Additionally, the following linearly-varying time-dependent pressure Neumann boundary condition is applied on $\Gamma_b$
\begin{equation} \label{eq:beam_pbc} 
	P(t) = (7.2599 \times 10^{4}) t\press, \hspace{0.5cm} t \geq 0,
\end{equation}
where $t$ is a pseudo-time variable, described in more detail below.
For the thermo-mechanical version of the beam problem (Section \ref{sec:beam_thermomech}), an additional temperature boundary condition 
is prescribed on $\Gamma_w$; we defer discussion of this boundary condition until Section \ref{sec:beam_thermomech}.  For both variants of this problem, the mechanical and the thermo-mechanical beam, the system, initially at 
rest, is solved quasi-statically by performing a homotopy continuation with respect to the pseudo-time variable
$t$ appearing in \eqref{eq:beam_pbc}.  This amounts to incrementing the 
applied boundary conditions (e.g., the pressure load \eqref{eq:beam_pbc}) with respect to $t$.


\subsubsection{Mechanical beam} \label{sec:beam_mech}
\label{subsubsec:MechBeamResults}
We will first investigate the purely mechanical response of the beam described above; hence, the governing PDEs are those 
described in Section \ref{sec:mechanical}.  For this problem, the displacement boundary conditions represent 235 constrained dofs 
out of the 1575 total,
    leaving this problem with 1340 free dofs.  


For the study performed herein, different material parameters are specified within different sets of blocks.  Let block $\mathcal{B}_a$ denote the union of 
blocks 1, 2 and 4 (shown in green, yellow and cyan, respectively, in Figure
\ref{fig:brickGeometry}(a)) and let block $\mathcal{B}_b$ denote the union of
blocks 3 and 5 
(shown in magenta and blue, respectively, 
in Figure  \ref{fig:brickGeometry}(a)).
In block $\mathcal{B}_a$, we set the Young's modulus, Poisson's ratio and density to the following values:
$E_a=1.103\e{11}\press$, $\nu_a=0.32$ and $\rho_a=7920\dens$, respectively.  The values of these same parameters in block $\mathcal{B}_b$, 
denoted by $E_b$, $\nu_b$ and $\rho_b$, were varied for the purpose of ROM training and testing, are provided in 
Table \ref{tab:MechBeamPredicCases}.  
These values were generated by performing Latin Hypercube (LHC) sampling
within the following parameter ranges: $E_b \in [1.27725 \times 10^{11}\press, 2.12875 \times 10^{11}\press]$, 
$\nu_b \in [0.24, 0.4]$, $\rho_b \in [5940\dens, 9900\dens]$.  
After generating the training data by quasi-static advancement of the problem to time $T = 7200\time$
with a time-step of 10$\time$,
and building a set of LSPG POD ROMs (with and without preconditioning)
for a wide range of ROM sizes $\nrb$ (ranging from 1 to 721) from the resulting set of 3605 snapshots, the ROMs were tested in four additional 
regimes, corresponding to four different values of the material parameters in block $\mathcal{B}_b$.  
The values of $E_b$, $\nu_b$ and $\rho_b$
for each of these training cases are provided in Table 
\ref{tab:MechBeamPredicCases}.
It is emphasized that the parameter variations in Table 
\ref{tab:MechBeamPredicCases} led to nontrivial differences in the displacement 
field, which varied by as much as 20\% between the different training and testing cases.

\begin{table}[H]
    \centering
	\caption{
		Mechanical beam problem: summary of parameters specified in the material model in block $\mathcal{B}_b$ 
	for the training and testing stages of the ROM process.}
    \begin{tabular}{r c | c c c}
	    Regime & Case & $E_b(\e{11})\units{\press}$ & $\nu_b$ & $\rho_b\units{\dens}$ \\ \hline
        \multirow{5}{*}{training} & 1 & 1.38002 & 0.28028 & 9194.74 \\
        \multirow{5}{*}{} & 2 & 2.11826 & 0.332646 & 7683.22 \\ 
        \multirow{5}{*}{} & 3 & 1.82559 & 0.395908 & 6150.4 \\ 
        \multirow{5}{*}{} & 4 & 1.56036 & 0.350415 & 9067.35 \\ 
        \multirow{5}{*}{} & 5 & 1.68463 & 0.256473 & 7466.27 \\ \hline
        \multirow{4}{*}{testing} & 1 & 1.50293 & 0.244704 & 6466.96 \\
        \multirow{4}{*}{} & 2 & 1.54545 & 0.304329 & 6774.12 \\
        \multirow{4}{*}{} & 3 & 1.47145 & 0.367092 & 8362.44 \\
        \multirow{4}{*}{} & 4 & 1.703 & 0.32 & 7920 \\
    \end{tabular}
    \label{tab:MechBeamPredicCases}
\end{table}




The main results for the mechanical beam problem are summarized in Figures 
\ref{fig:MechBeamErr}--\ref{fig:MechBeamCond}.  
The reader can observe that the classical (unpreconditioned) LSPG solution 
does not appear in these plots.  The LSPG solution is not included in our 
results summary because the unpreconditioned LSPG ROMs were not able to converge
to a nonlinear solver tolerance smaller than $\mathcal{O}(1)$ for any of the basis 
sizes considered.
Figure \ref{fig:MechBeamErr}
plots the ROM global relative error $\mathcal{\epsilon}$ (calculated using \eqref{pterr}) as a function of the basis dimension for various preconditioned LSPG ROMs, including the projected solution increment ROM, indicated with a black
dashed line.  It can be seen that, in general, as the preconditioner is improved (from {\tt Jacobi}, 
to {\tt Gauss-Seidel}, 
to {\tt ILU}), the accuracy of the resulting ROM solution is improved, approaching the ideal projected 
solution increment solution, as expected from the discussion in Section \ref{sec:precond}.  
Next, in Figure \ref{fig:MechBeamTimes}, we report wall times for the four ROMs evaluated for 
each of our testing cases.
The wall times reported are averaged over
four processors of a Linux workstation having 20 Intel Xeon CPU E5-2670 v2 CPUs.
The reader can observe that the more sophisticated {\tt ILU} preconditioner results
in slightly larger wall times in general compared to the {\tt Jacobi} and {\tt Gauss-Seidel} 
preconditioners.  This is due to a larger preconditioner construction time associated 
with the {\tt ILU} preconditioner.  As expected, the ideal preconditioner solution, 
which requires calculating the action of $\jkinv$
in every iteration as discussed in Section
\ref{sec:precond_eval}, is in general the costliest to obtain.
The results in Figures \ref{fig:MechBeamErr} and \ref{fig:MechBeamTimes} 
are combined in Figure \ref{fig:beam_mech_pareto}, which shows a Pareto plot for the mechanical 
beam problem, in which the global relative error $\epsilon$ is plotted as a function of the 
wall time.  From this plot, it is possible to identify the optimal preconditioner to use, 
based on one's error and CPU-time requirements.  The blue line in this plot traces points
having Pareto optimality, which define the so-called Pareto front.
\begin{figure}[h!tb]
 \centering
   \begin{tabular}{cc}
	   \includegraphics[width=0.5\textwidth]{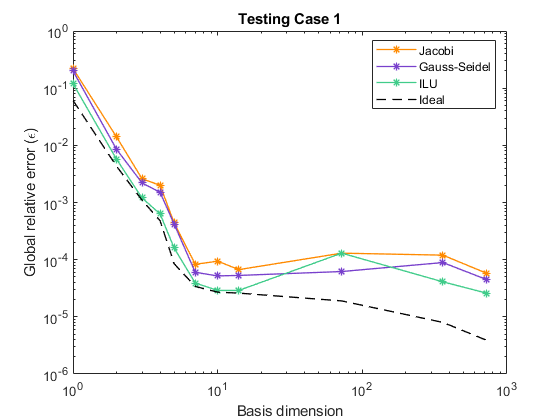}&
   \includegraphics[width=0.5\textwidth]{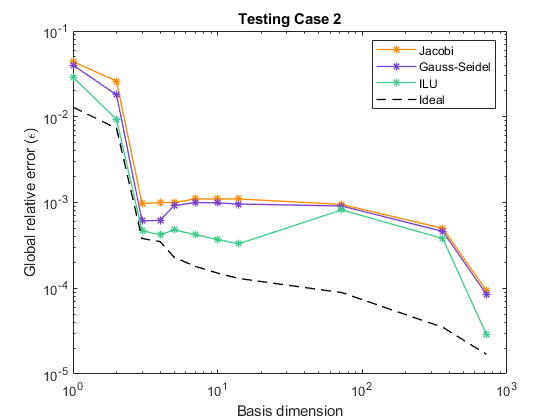} \\
	   (a) Testing case 1 & (b) Testing case 2  \\
	   \includegraphics[width=0.5\textwidth]{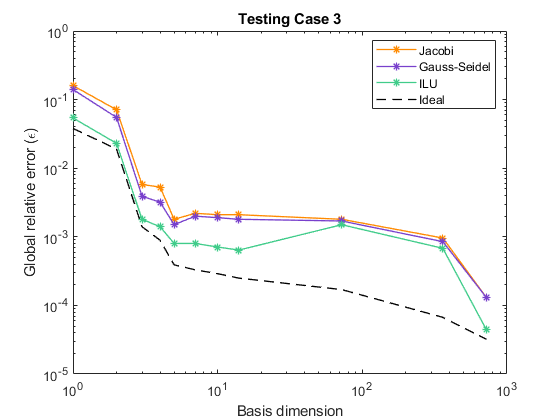}&
   \includegraphics[width=0.5\textwidth]{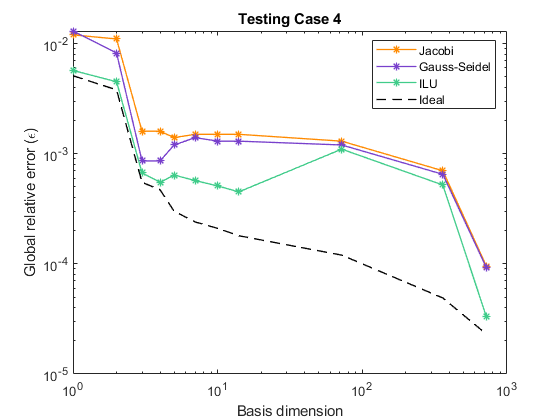} \\
	   (c) Testing case 3 & (d) Testing case 4  
   \end{tabular}  
	\caption{
		Mechanical beam problem: global relative errors $\mathcal{\epsilon}$ \eqref{pterr} for various ROMs as a function
	of the basis dimension for the four testing cases in Table \ref{tab:MechBeamPredicCases}.}
  \label{fig:MechBeamErr}
\end{figure}
Finally, in Figure \ref{fig:MechBeamCond}, we report the average condition number of the reduced Jacobian, 
$\PrecondPetrovGalerkin{\jac}\newton{k}$ (see \eqref{lspgitne_prec}) over all Gauss--Newton iterations
for each of the testing cases summarized in 
Table \ref{tab:MechBeamPredicCases}.  Although the LSPG ROMs were not convergent for this problem, 
we include condition numbers for the LSPG ROMs in this plot. 
The reader can observe that, as expected, the addition of a preconditioner $\Mk$ 
improves the condition number of the corresponding ROM relative to the baseline LSPG ROM.  
For all ROMs evaluated, condition numbers grow steadily with the basis dimension $\nrb$; however, preconditioning using the relatively  simple 
preconditioners considered herein (the first three rows of Table \ref{tab:preconds}) 
is able to reduce the condition number of the reduced Jacobian of the LSPG ROM by up to an order of magnitude.
By design, 
the ideal preconditioner (equivalent to the projected 
solution increment) gives rise to a ROM system having a perfect condition 
number of one for all basis dimensions.  

\begin{figure}[h!tb]
 \centering
   \begin{tabular}{cc}
	   \includegraphics[width=0.5\textwidth]{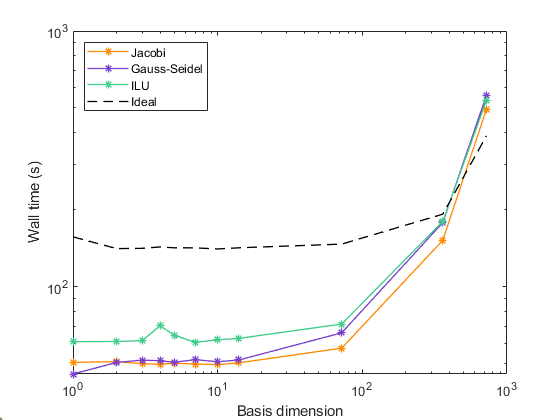}&
   \includegraphics[width=0.5\textwidth]{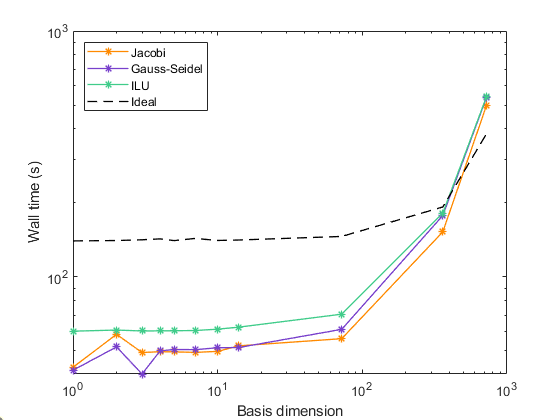} \\
	   (a) Testing case 1 & (b) Testing case 2  \\
	   \includegraphics[width=0.5\textwidth]{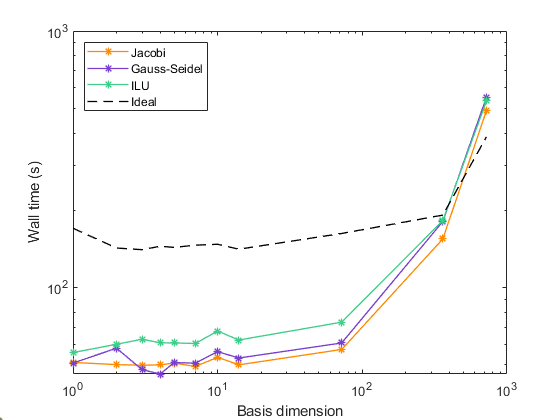}&
   \includegraphics[width=0.5\textwidth]{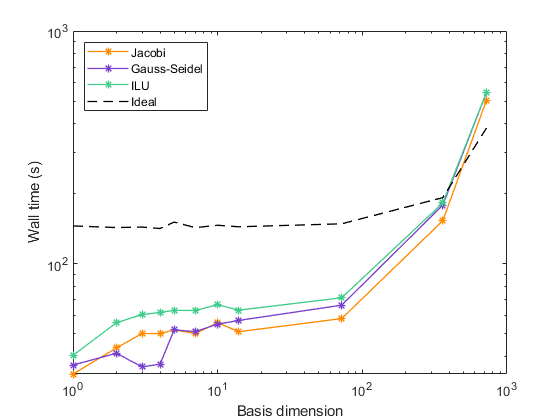} \\
	   (c) Testing case 3 & (d) Testing case 4  \\
   \end{tabular}  
	\caption{Mechanical beam problem: wall times (in s, averaged over 4 processors) for various ROMs as a function of the basis dimension for the four testing cases in Table \ref{tab:MechBeamPredicCases}.}
  \label{fig:MechBeamTimes}
\end{figure}

\begin{figure}[h!tb]
 \centering
\includegraphics[width=1.0\textwidth]{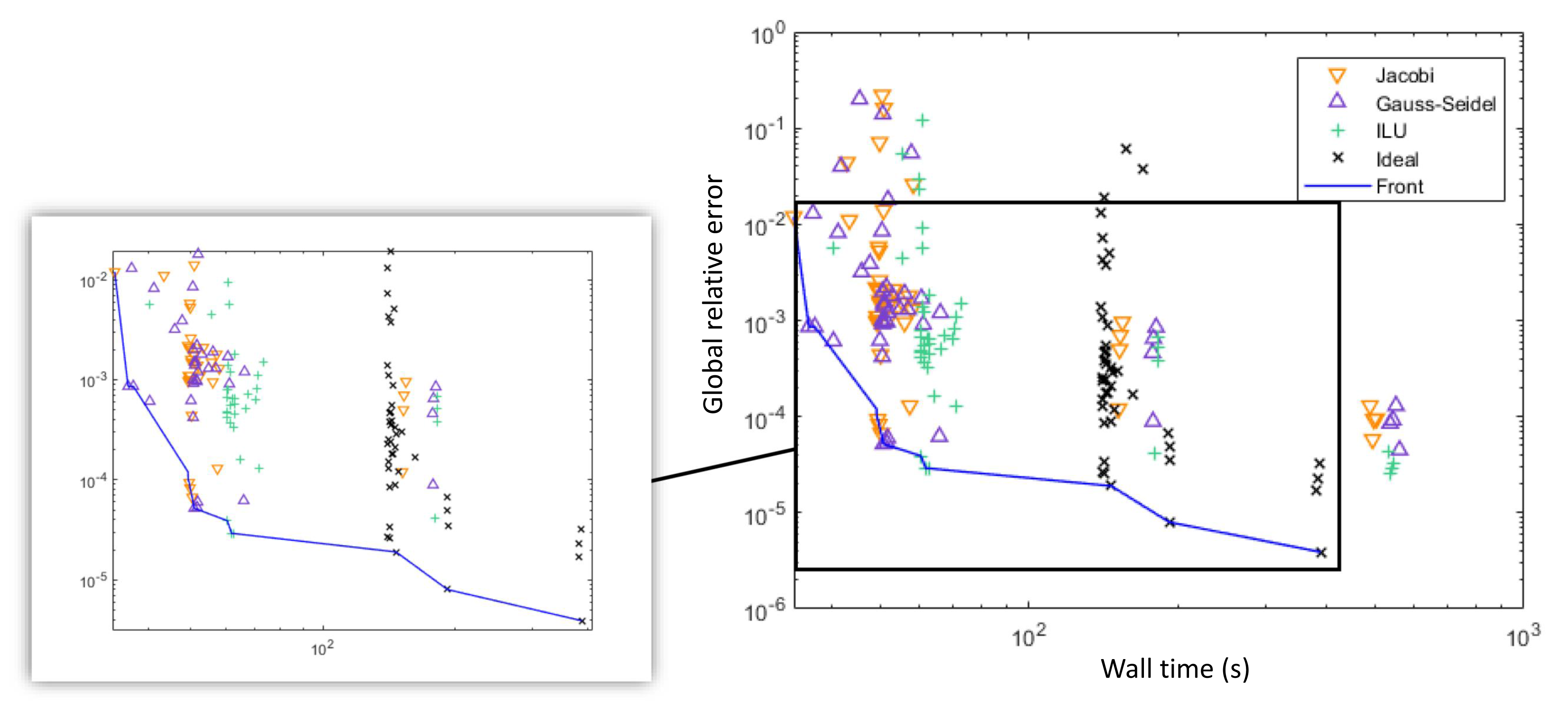}
	\caption{Mechanical beam problem: Pareto plot showing the total  wall time (in s, averaged over 4 processors) 
	versus the global relative error $\epsilon$.  Results for all four testing cases in Table
	\ref{tab:MechBeamPredicCases} are used to generate this figure.} 
  \label{fig:beam_mech_pareto}
\end{figure}

\begin{figure}[h!tb]
 \centering
   \begin{tabular}{cc}
	   \includegraphics[width=0.5\textwidth]{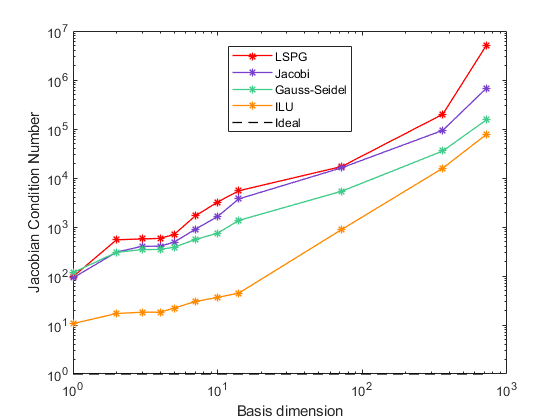}&
   \includegraphics[width=0.5\textwidth]{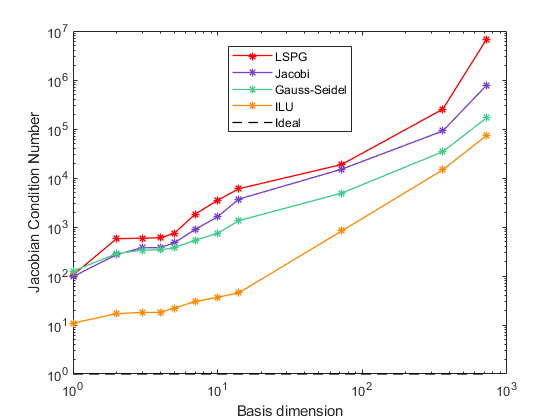} \\
	   (a) Testing case 1 & (b) Testing case 2  \\
	   \includegraphics[width=0.5\textwidth]{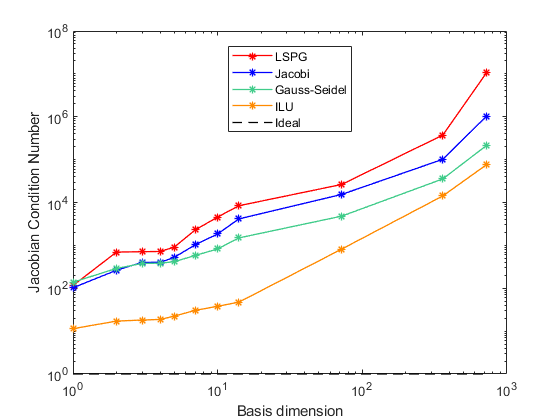}&
   \includegraphics[width=0.5\textwidth]{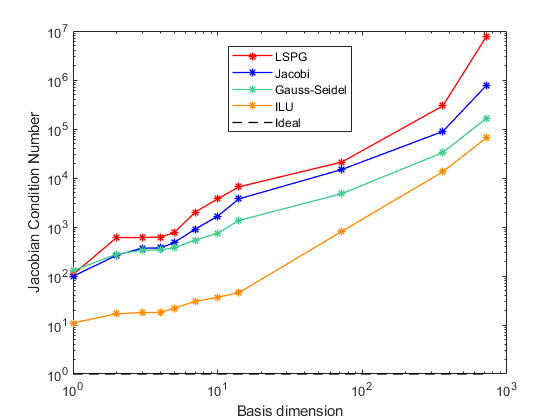} \\
	   (c) Testing case 3 & (d) Testing case 4  \\
  \end{tabular}
	\caption{Mechanical beam problem: average reduced Jacobian condition number for various ROMs as a function
	of the basis dimension for the four testing cases in Table \ref{tab:MechBeamPredicCases}.}
  \label{fig:MechBeamCond}
\end{figure}

\subsubsection{Thermo-mechanical beam} \label{sec:beam_thermomech}


We now turn our attention to the thermo-mechanical version of the beam problem, in which the mechanical equations of Section \ref{sec:mechanical} are 
augmented with the thermal equation in Section \ref{sec:thermo-mechanical}.  In the current \albany{} implementation, 
the thermal and mechanical equations are coupled monolithically within the code, with 
the coupling occurring at the level of the material model, as described earlier in Section \ref{sec:thermo-mechanical}.   


In addition to the displacement and pressure \eqref{eq:beam_pbc}
boundary conditions described previously in Section \ref{sec:beam_mech}, we impose a
temperature Dirichlet boundary condition on $\Gamma_w$ in Figure \ref{fig:brickGeometry}(a).  The value of the prescribed 
temperature varies linearly from 293 to 393$\temp$ over the course of the simulation: 
\begin{equation} \label{eq:beam_Tbc}
	T(t) = (100t + 293)\temp, \hspace{0.5cm} t \geq 0. 
\end{equation}
As before, the variable $t$ in \eqref{eq:beam_Tbc} is a pseudo-time variable which is incremented 
quasi-statically via homotopy continuation to emulate time-dependent behavior in a quasi-static framework.
Note that the temperature boundary condition \eqref{eq:beam_Tbc} is not the only one 
advanced in time quasi-statically; the pressure Neumann boundary condition \eqref{eq:beam_pbc} 
is also incremented, as discussed earlier in Section \ref{sec:beam_mech}.   
As before, we use the Neohookean material model for the mechanical problem, and specify different material properties in block $\mathcal{B}_a$ and $\mathcal{B}_b$.  In block $\mathcal{B}_b$, 
the mechanical properties  are $E_a=1.103\e{9}\press$ (Young's modulus), 
$\nu_a=0.32$ (Poisson's ratio) and $\rho_a=7.92\e{-5}\dens$ (density).  As discussed in Section \ref{sec:thermo-mechanical}, the Helmholtz free-energy 
density $A$ appearing in  \eqref{eq:euler-lagrange}  is a function of the temperature $T$ in the thermo-mechanical case. 
This thermo-mechanical formulation requires one to specify a reference temperature.  In block $\mathcal{B}_a$, we 
use a reference temperature of $T_{a, \text{ref}} = 293 \temp$. 
In block $\mathcal{B}_b$, the mechanical properties of the Neohookean  material as well as the reference temperature, denoted by 
$T_{b, \text{ref}}$ are varied for both the training and the testing cases, as summarized in Table \ref{tab:ThermoMechBeamPredicCases}.  
Additionally, we assume the material is isotropic, i.e., $\boldsymbol{K} := k \I_3$ in \eqref{heat_flux}, 
where $\I_3$ is the $3 \times 3$ identity matrix, with thermal diffusivity $k_a = k_b = 1.1870\times 10^3$ m$^2$/s in 
both blocks.  We also assume the thermal expansion coefficient \eqref{eq:FT} is the same in both blocks, taking the value 
$\alpha_a = \alpha_b = 1\e{-5}\temp^{-1}$.
There are no applied forces and the system is initially at rest, with an initial temperature of 293$\temp$.  
For the quasi-static advancement of the system, we use a 
step size of 1$\time$ and run the problem for a total of 720 steps.  
The displacement and temperature boundary conditions represent 260 constrained dofs out of the 2100 total dofs, leaving this problem with 1840 free dofs.
As for the mechanical version of this problem, the parameter ranges in Table \ref{tab:ThermoMechBeamPredicCases} 
are selected such that nontrivial variations in the solution are observed; specifically,
the displacement varies by as much as 60\% as the material properties are perturbed according to Table \ref{tab:ThermoMechBeamPredicCases}.

The ROM was trained with six sets of parameters in block $\mathcal{B}_b$, chosen using LHC sampling and summarized in Table \ref{tab:ThermoMechBeamPredicCases}.  The following parameter ranges for the LHC sampling were employed:
$E_b \in [1.27725 \times 10^9\press, 2.12875 \times 10^9\press]$, $\nu_b \in [0.24, 0.40]$, $\rho_b \in [5.94\times 10^{-5}\dens, 
9.9\times 10^{-5}\dens]$, 
$T_{b, \text{ref}} \in [219.75\temp, 366.25\temp]$.  
This training generated a total of 4318 snapshots, which were used to build various LSPG/POD ROMs (with and without preconditioning) ranging in size from 1 to 721 modes, like in the mechanical variant of this problem.
These ROMs were tested in the predictive regime over four additional sets of parameters, also summarized in Table \ref{tab:ThermoMechBeamPredicCases}. 
\begin{table}[H]
    \centering
	\caption{Thermo-mechanical beam problem: summary of parameters specified in the material model in block $\mathcal{B}_b$ 
	for the training and testing stages of the ROM process.}
    \begin{tabular}{r c | c c c c}
	    Regime & Case & $E_b(\e{9})\units{\press}$ & $\nu_b$ & $\rho_b(\e{-5})\units{\dens}$ & $T_{b, \text{ref}}\units{\temp}$ \\ \hline
        \multirow{6}{*}{training} & 1 & 2.01313 & 0.285907 & 7.94827 & 273.657 \\
        \multirow{6}{*}{} & 2 & 1.71637 & 0.332083 & 6.93965 & 318.406 \\ 
        \multirow{6}{*}{} & 3 & 1.96881 & 0.3478 & 9.37181 & 301.406 \\ 
        \multirow{6}{*}{} & 4 & 1.28954 & 0.29427 & 9.14636 & 365.378 \\ 
        \multirow{6}{*}{} & 5 & 1.61326 & 0.262464 & 6.32164 & 223.434 \\ 
        \multirow{6}{*}{} & 6 & 1.54724 & 0.374118 & 7.31561 & 245.778 \\ \hline
        \multirow{4}{*}{testing} & 1 & 1.52473 & 0.27925 & 8.80694 & 266.674 \\
        \multirow{4}{*}{} & 2 & 1.31153 & 0.345538 & 7.58234 & 333.462 \\
        \multirow{4}{*}{} & 3 & 1.37015 & 0.246513 & 7.73303 & 345.942 \\
        \multirow{4}{*}{} & 4 & 1.703 & 0.32 & 7.92 & 293 \\
    \end{tabular}
    \label{tab:ThermoMechBeamPredicCases}
\end{table}

The main results for the thermo-mechanical beam problem are summarized in 
Figures \ref{fig:ThermoMechBeamErr1}--\ref{fig:ThermoMechBeamCond}.
It is interesting to observe that, despite the multi-physics and multi-scale
nature of this problem, the classical (unpreconditioned) LSPG ROMs, 
plotted in red in these figures, deliver convergent
solutions for the smaller basis dimensions considered.

Figures \ref{fig:ThermoMechBeamErr1} and \ref{fig:ThermoMechBeamErr2} 
report the global relative error $\mathcal{\epsilon}$ in each of the ROM solutions, 
computed using \eqref{pterr}.  The reader 
can observe by examining this figure that, while convergence with basis refinement is observed up to a point 
(until $\nrb \approx 7$ modes) for 
the unpreconditioned LSPG ROM, 
for larger basis dimensions, the error grows steadily before reaching a point where
a lack of convergence is observed.  
The results are markedly different for the preconditioned ROMs.  By introducing 
preconditioning, it is possible to reduce $\mathcal{\epsilon}$ by between two and six orders of 
magnitude, depending on the basis dimension, and there are no convergence issues, like for the unpreconditioned
cases.  Moreover, one can see 
that all the preconditioned LSPG ROMs achieve errors which are close to 
(less than one order of magnitude greater than)  the error obtained
by the projected solution increment ROM, which represents an ideal preconditioned LSPG ROM. 



\begin{figure}[h!tb]
 \centering
   \begin{tabular}{c}
	   \includegraphics[width=1.0\textwidth]{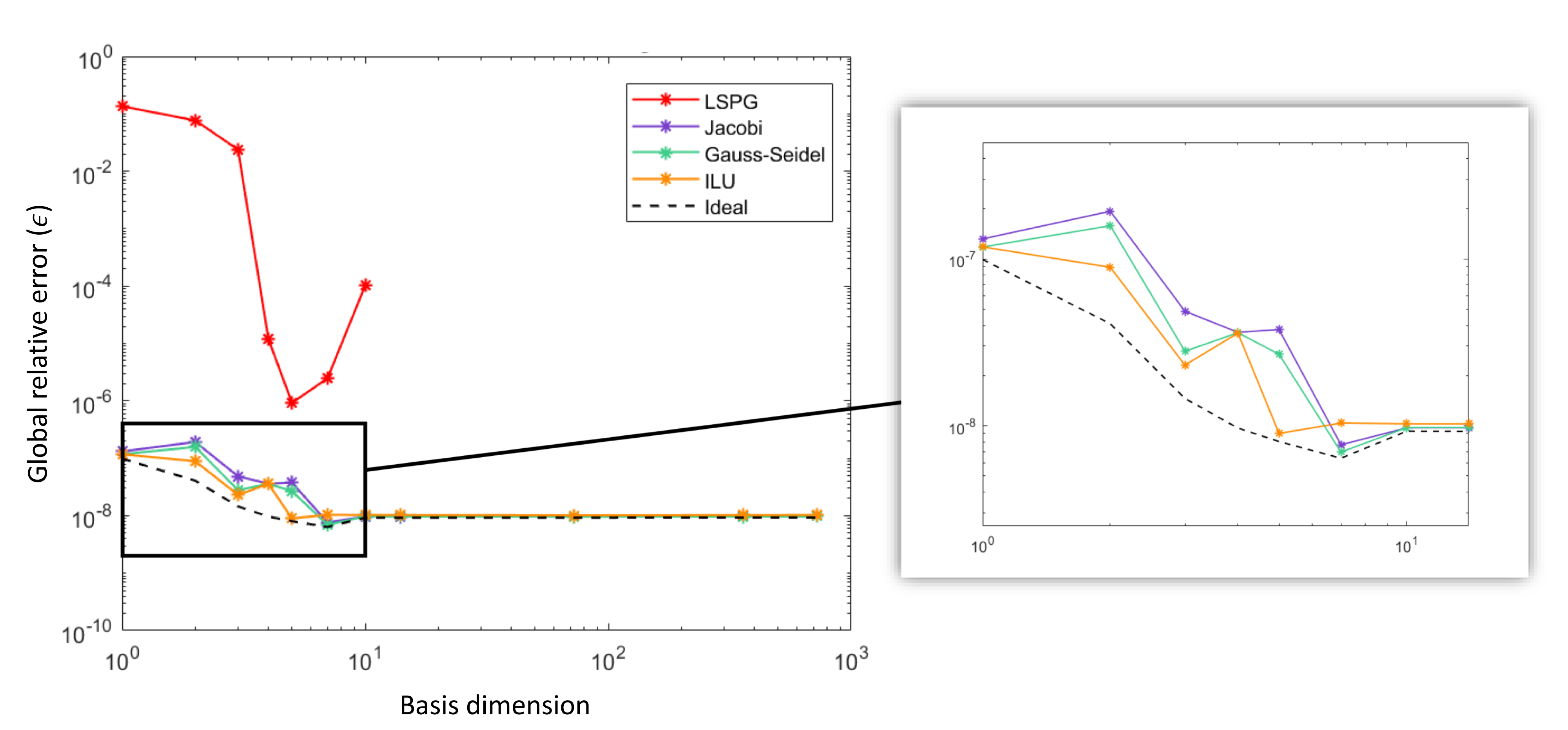} \\
	   (a) Testing case 1  \\
	   \includegraphics[width=1.0\textwidth]{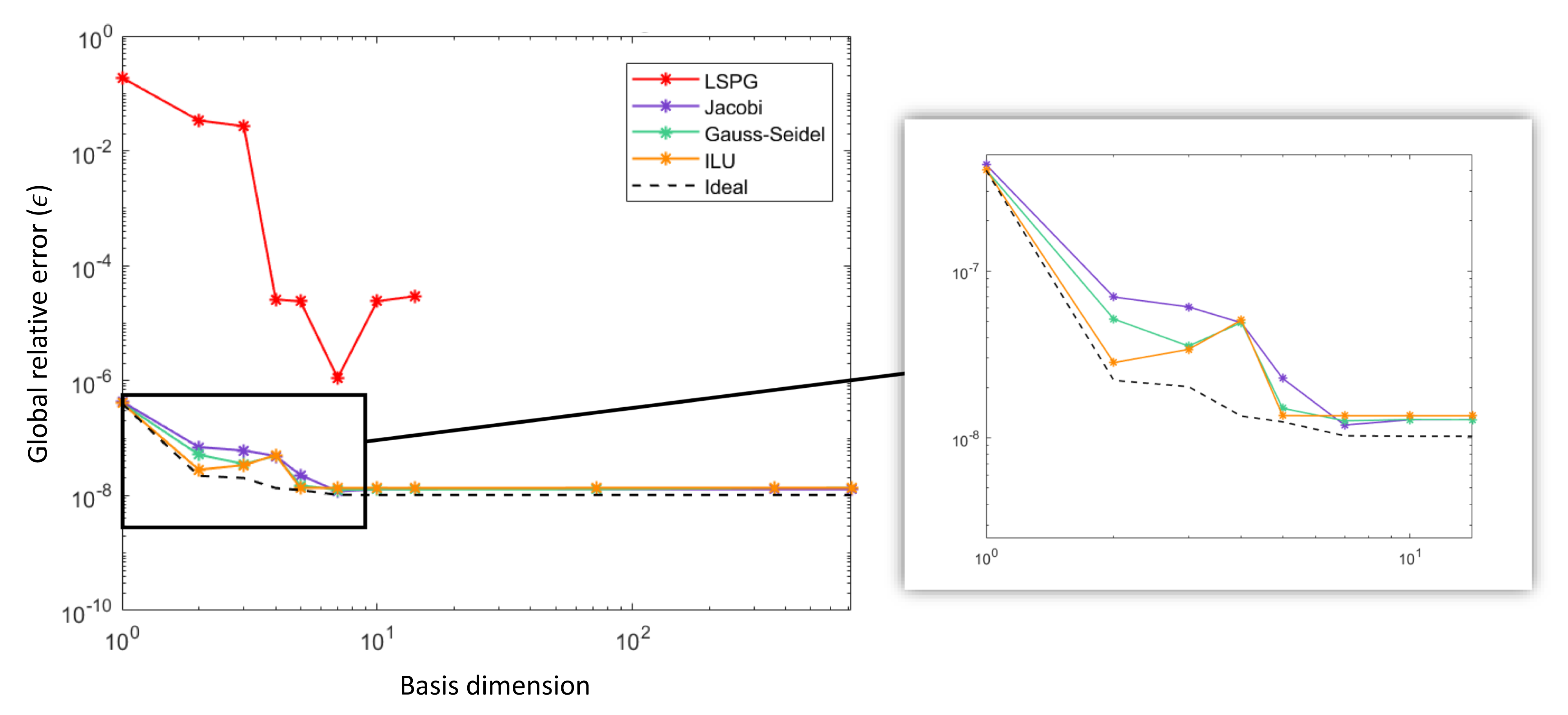} \\
	   (b) Testing case 2  
   \end{tabular} 
	\caption{Thermo-mechanical beam problem: global relative errors $\mathcal{\epsilon}$ \eqref{pterr} for various ROMs as a function
	of the basis dimension for testing cases 1--2 in Table \ref{tab:ThermoMechBeamPredicCases}.}
  \label{fig:ThermoMechBeamErr1}
\end{figure}

\begin{figure}[h!tb]
 \centering
   \begin{tabular}{c}
	   \includegraphics[width=1.0\textwidth]{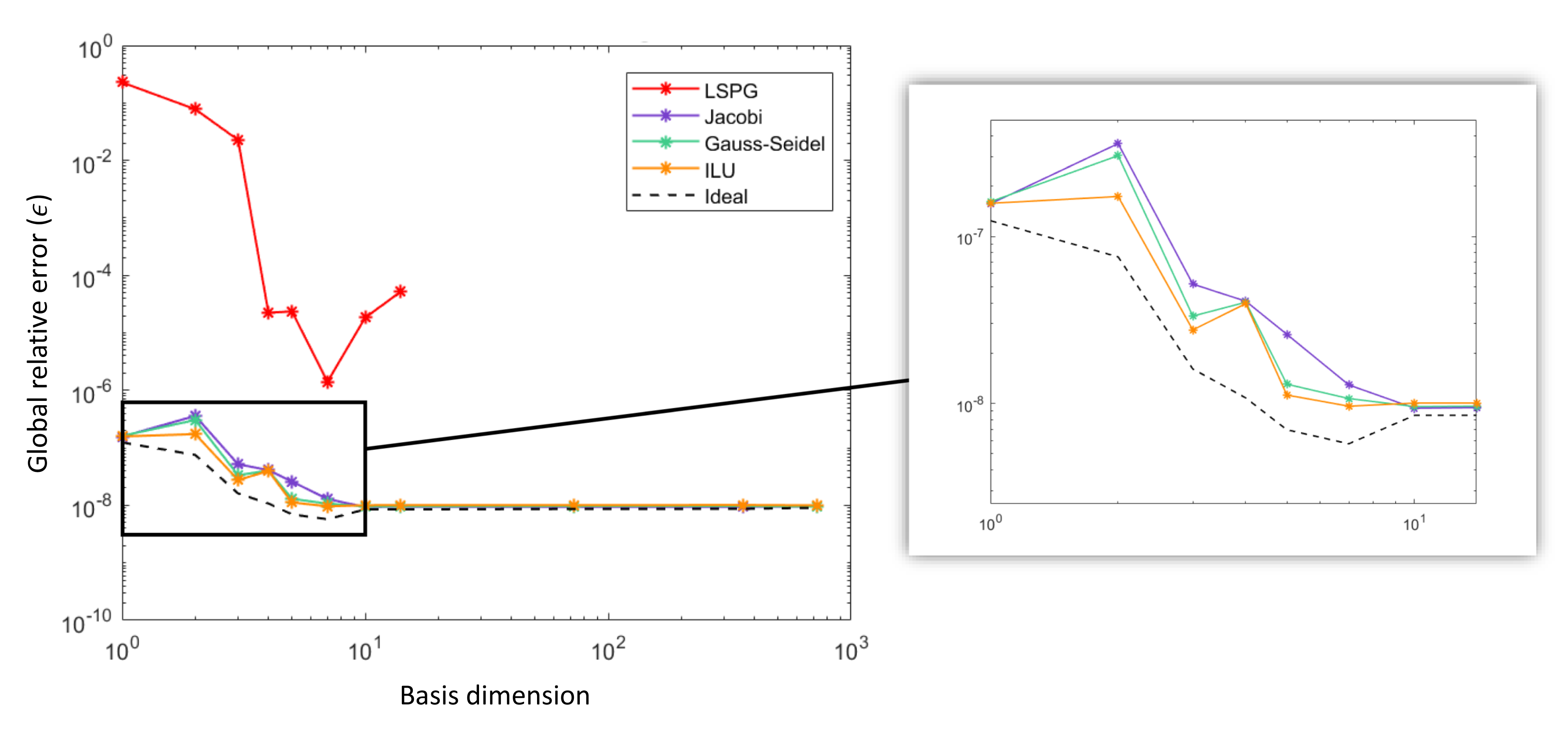} \\
	   (a) Testing case 3 \\
	   \includegraphics[width=1.0\textwidth]{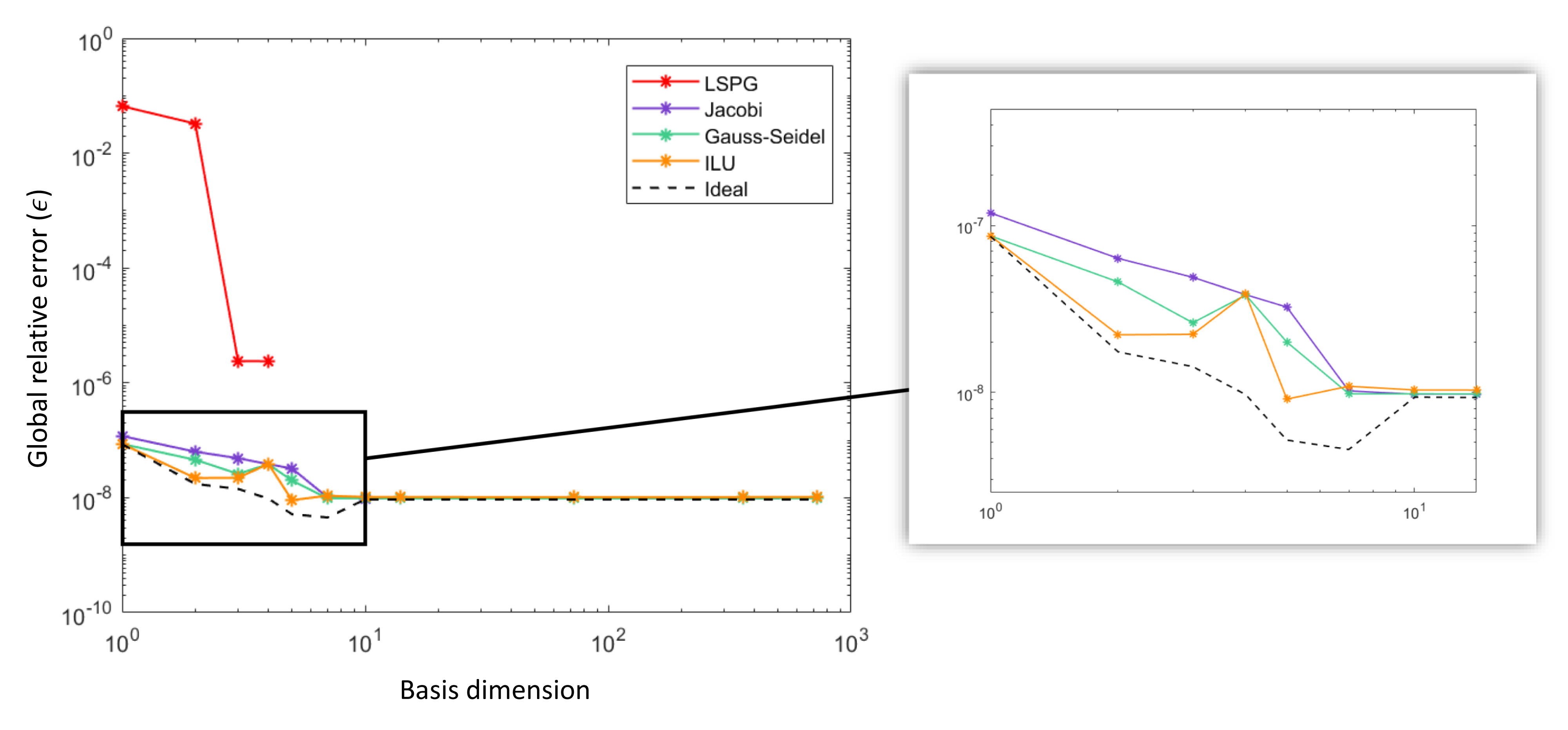} \\
	   (b) Testing case 4 
   \end{tabular} 
	\caption{Thermo-mechanical beam problem: global relative errors $\mathcal{\epsilon}$ \eqref{pterr} for various ROMs as a function
	of the basis dimension for testing cases 3--4 in Table \ref{tab:ThermoMechBeamPredicCases}.}
  \label{fig:ThermoMechBeamErr2}
\end{figure}

We next investigate the effect of preconditioning on the overall 
run time of the ROM problem (Figure \ref{fig:ThermoMechBeamTimes}).
As for the mechanical version of this problem, we report wall times averaged over 
eight processors of a Sandia National Laboratories' Linux cluster known as Uno, which has 201 dual-socket 
eight-core 2.7 GHx Intel Sandy Bridge CPUs.
The reader can observe from Figure \ref{fig:ThermoMechBeamTimes} that for the majority of testing cases and basis dimensions, 
the preconditioned
LSPG ROMs considered achieve smaller wall time than the baseline (unpreconditioned) LSPG ROM.  
Combining these 
results with our earlier error results (Figures \ref{fig:ThermoMechBeamErr1}--\ref{fig:ThermoMechBeamErr2}), we can confidently conclude that preconditioning
is extremely advantageous for the thermo-mechanical beam problem: in general, for the same computational cost, one is able to achieve
a much smaller global relative error with preconditioning.  This result is confirmed in the Pareto plot
shown in Figure \ref{fig:thermo-mech_beam_pareto}.  
As before, the black line in Figure \ref{fig:ThermoMechBeamTimes}
represents the ideal preconditioned ROM.  
This ROM has a much higher computational 
cost than the other ROMs evaluated, as expected.




\begin{figure}[h!tb]
 \centering
   \begin{tabular}{cc}
	   \includegraphics[width=0.5\textwidth]{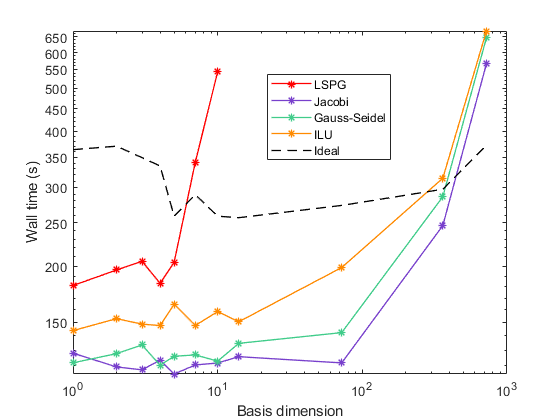}  & 
	   \includegraphics[width=0.5\textwidth]{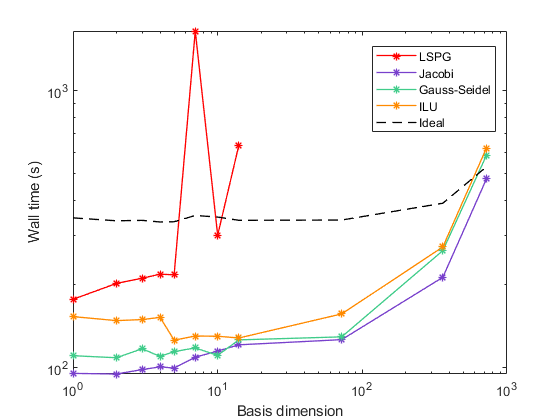}  \\ 
	   (a) Testing case 1 & (b) Testing case 2  \\
	   \includegraphics[width=0.5\textwidth]{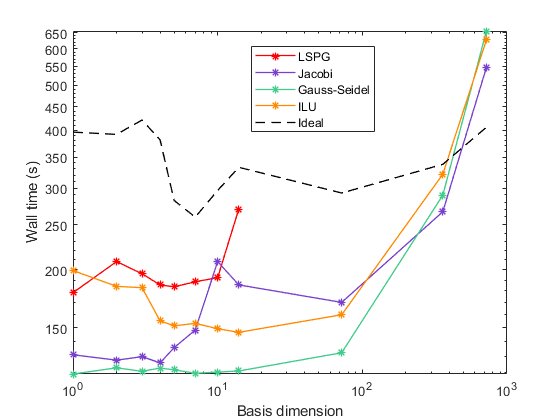}  & 
	   \includegraphics[width=0.5\textwidth]{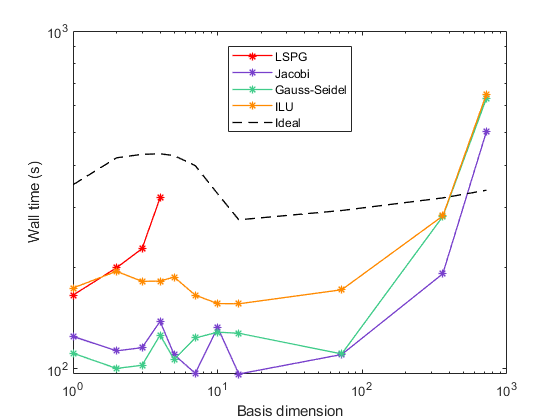}  \\
	   (c) Testing case 3 & (d) Testing case 4  \\
   \end{tabular} 
	\caption{Thermo-mechanical beam problem: wall times (in$\time$, averaged over 8 processors) for various ROMs as a function of the basis dimension for the four testing cases in Table \ref{tab:MechBeamPredicCases}.}
  \label{fig:ThermoMechBeamTimes}
\end{figure}

\begin{figure}[h!tb]
 \centering
\includegraphics[width=1.0\textwidth]{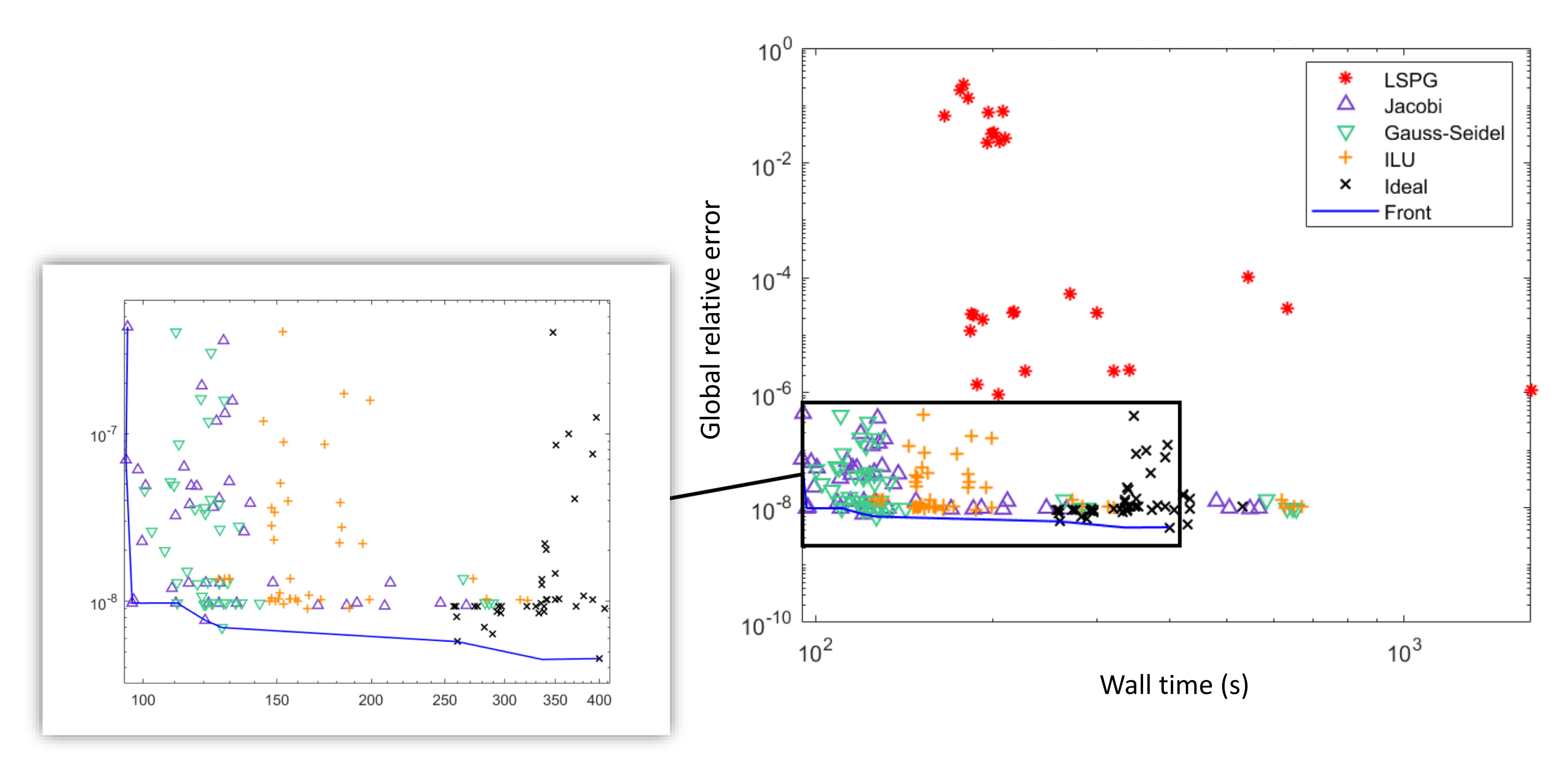}
	\caption{Thermo-mechanical beam problem: Pareto plot showing the total wall time (in s, averaged over 4 processors) 
	versus the global relative error $\epsilon$.  Results for all four testing cases in Table
	\ref{tab:ThermoMechBeamPredicCases} are used to generate this figure.} 
  \label{fig:thermo-mech_beam_pareto}
\end{figure}


Finally, in Figure \ref{fig:ThermoMechBeamCond}, we report average condition numbers of the 
reduced Jacobians $\PrecondPetrovGalerkin{\jac}\newton{k}$ (see \cref{lspgitne_prec}) encountered
during the Gauss--Newton iteration process
for each of the ROMs considered.  The most striking observation that can be made from this figure is that the regular LSPG Jacobians are very ill-conditioned, 
with condition numbers ranging between $\mathcal{O}(10^{14})$ and $\mathcal{O}(10^{20})$ depending on $\nrb$, the basis dimension.  These condition numbers 
are between seven to ten times greater than the condition numbers for the mechanical variant of this problem (see Figure \ref{fig:MechBeamCond}).  The extreme 
ill-conditioning exhibited by the thermo-mechanical beam problem can be attributed to the extreme differences in scales between the thermal 
and the mechanical problems: whereas the temperature solution is $\mathcal{O}(100)$, 
the displacement solution varies between $\mathcal{O}(10^{-7})$ and $\mathcal{O}(10^{-5})$ 
during the duration of the simulation, a difference of between seven and nine orders of magnitude. 
Our results demonstrate that, 
by introducing a simple preconditioning 
strategy into the LSPG formulation, it is possible to bring down the reduced Jacobian condition numbers by as many as ten orders of magnitude, thereby 
alleviating to a large extent the scaling issue introduced by the large difference in scales between the thermal and mechanical problems 
(a similar result was discovered by 
Washabaugh in \cite{kyle}, but for a different application, namely computational fluid dynamics).
Experience 
involving other  solid mechanics applications within \albany{} \cite{schwarz} suggests that further improvements can be made by introducing physics-specific 
block preconditioners into the framework; however, a study of this sort goes beyond the scope of this work.  As before, the condition number 
of the reduced Jacobian for the projected solution increment ROM is identically one regardless of the basis dimension.  


\begin{figure}[h!tb]
 \centering
   \begin{tabular}{cc}
	   \includegraphics[width=0.5\textwidth]{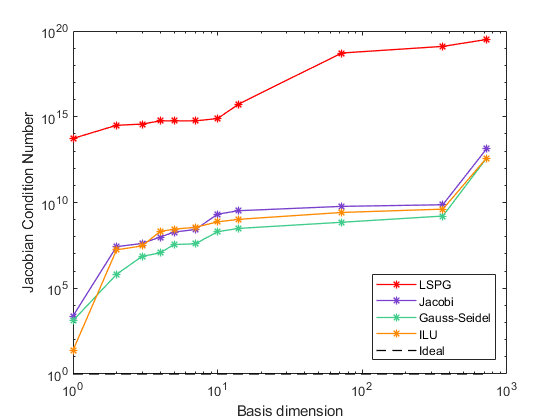}&
	   \includegraphics[width=0.5\textwidth]{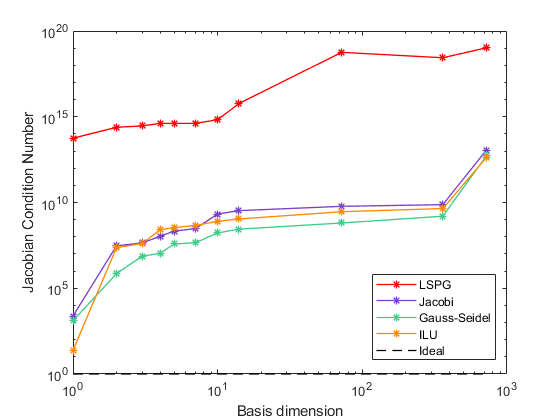} \\
	   (a) Testing case 1 & (b) Testing case 2  \\
	   \includegraphics[width=0.5\textwidth]{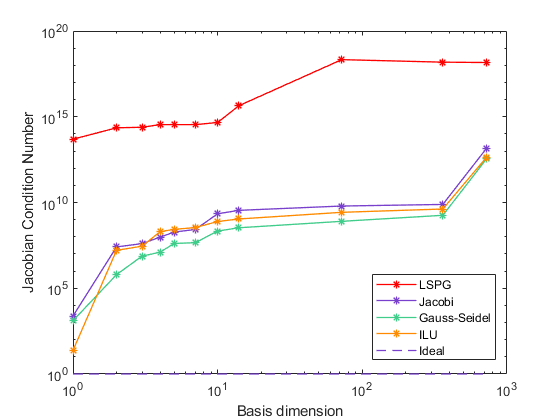}&
   \includegraphics[width=0.5\textwidth]{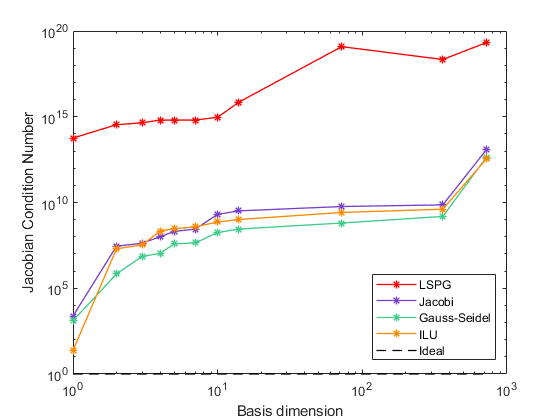} \\
	   (c) Testing case 3 & (d) Testing case 4  \\
   \end{tabular}
	\caption{Thermo-mechanical beam problem: average reduced Jacobian condition number for various ROMs as a function
	of the basis dimension for the four testing cases in Table \ref{tab:ThermoMechBeamPredicCases}.}
  \label{fig:ThermoMechBeamCond}
\end{figure}


\subsection{Thermo-mechanical pressure vessel}  \label{sec:vessel}

The last test case considered aims at studying the thermo-mechanical response of a 
3D pressured cylindrical vessel shown in Figure \ref{fig:PCAPDBCs}(b).
This problem is substantially larger and more 
realistic than the beam problem considered earlier in Section \ref{sec:beam}.  
Our approach herein is to model a quarter of the pressure vessel (denoted by $\Omega$), 
as shown in Figure \ref{fig:PCAPDBCs}(a), and apply appropriate symmetry boundary 
conditions on the relevant sides to emulate a simulation on the pressure vessel in its entirety.  
We discretize our domain $\Omega$ using a mesh comprised of 77,768 hexahedral elements, which give rise 
to 92,767 nodes.  Similar to the beam problem, the domain $\Omega$ is comprised of five element blocks, 
depicted in green (block 1), yellow (block 2), blue (block 3), magenta (block 4) and cyan (block 5), as shown in Figure \ref{fig:PCAPBlocks},
so as to enable to specification of different materials and/or material properties in different parts of the domain.  
\begin{figure}[h!tb]
	\begin{minipage}[c]{0.69\textwidth}
		\subfig{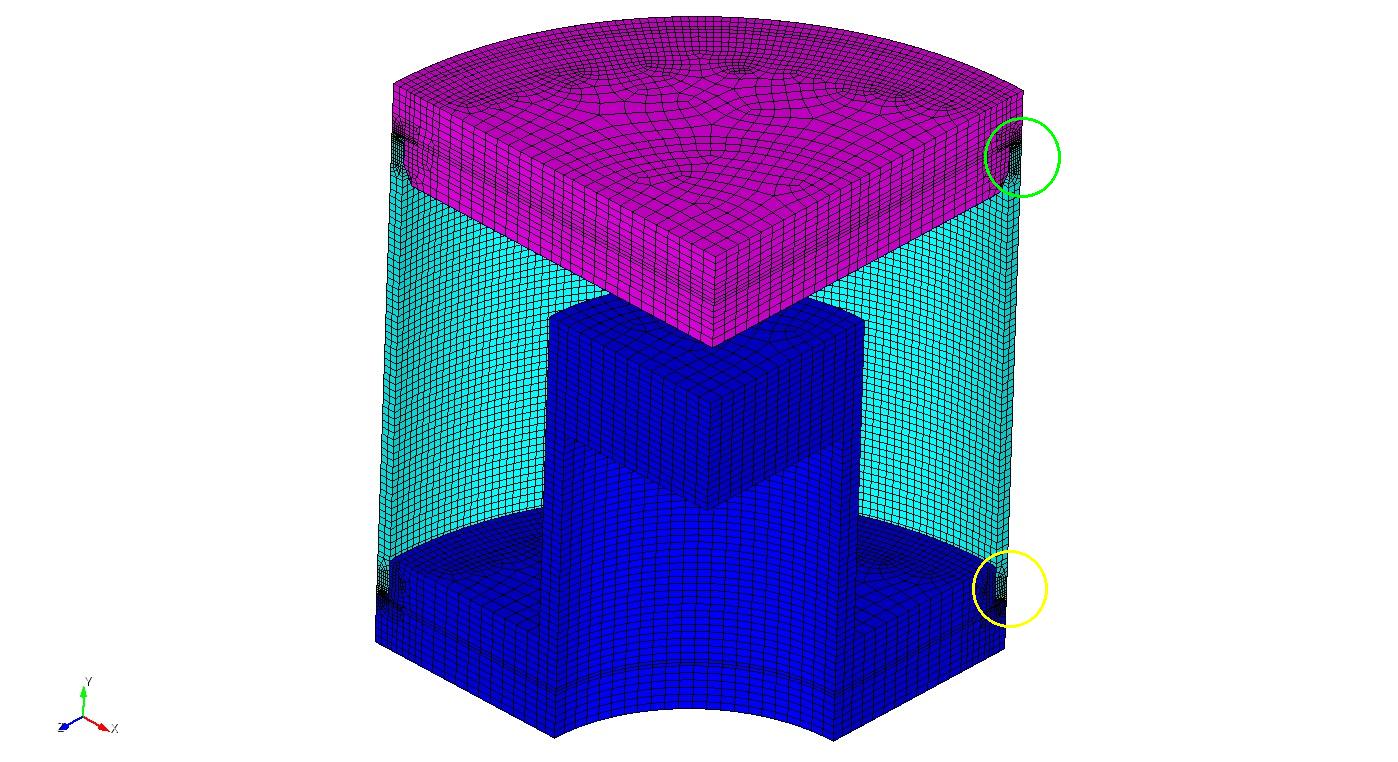}{All material blocks}{fig:PCAPBlocks_main}
	\end{minipage}
	\begin{minipage}[c]{0.30\textwidth}
		\subfig{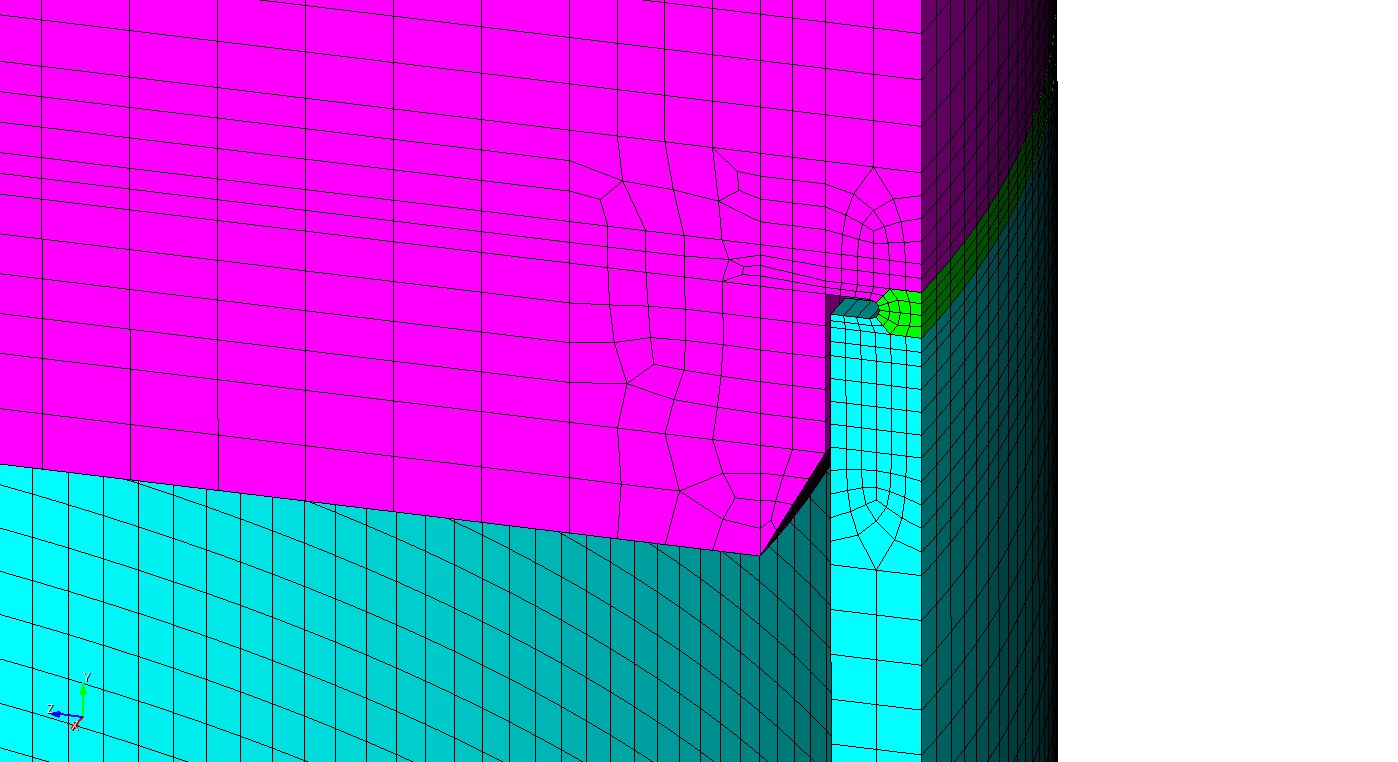}{Focus on block 1}{fig:PCAPBlocks_zoom1}

		\subfig{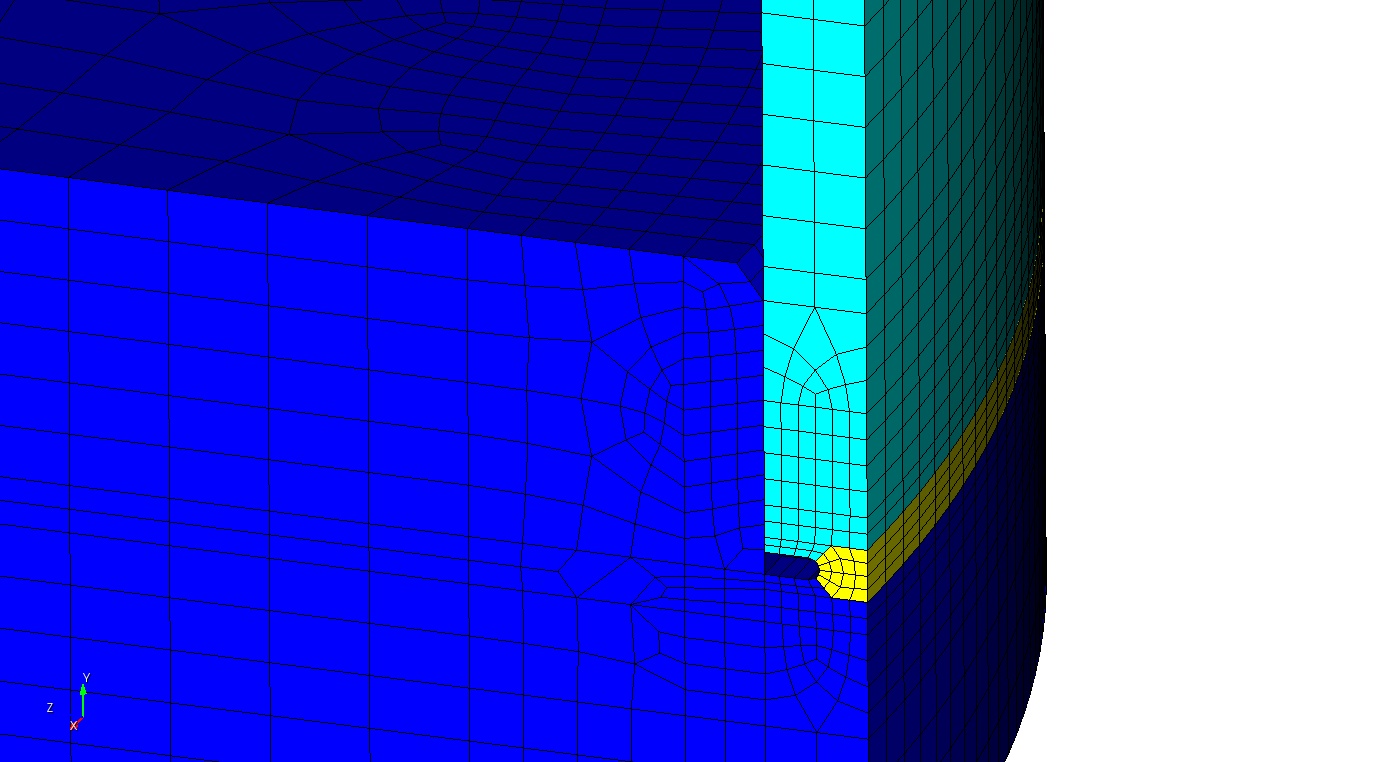}{Focus on block 2}{fig:PCAPBlocks_zoom2}
	\end{minipage}
	\caption{Thermo-mechanical pressure vessel problem: computational domain.  Different colors indicate different material
	blocks; blocks 1, 2, 3,  4 and 5 are represented in green, yellow, magenta, cyan and blue, 
	respectively. Subfigures (b) and (c) 
	zoom in on the regions circled in green and yellow respectively in subfigure (a), and depict the location of blocks 1 and 2, respectively.  }
	\label{fig:PCAPBlocks}
\end{figure}

Figure \ref{fig:PCAPDBCs} depicts the boundaries on which boundary conditions are applied 
for the thermo-mechanical pressure vessel problem.  Let $\Gamma$ denote the boundary of our computational domain $\Omega$. 
As 
shown in Figure \ref{fig:PCAPDBCs}(a), we decompose $\Gamma$ as $\Gamma:=\Gamma_1 \cup \Gamma_2 \cup \Gamma_3 \cup \Gamma_4 \cup \Gamma_5$ where 
$\Gamma_i \cap \Gamma_j = \emptyset$ for $i, j = 1, ..., 5$ and $i\neq j$.  In Figure \ref{fig:PCAPDBCs}(a), $\Gamma_1$ is shown in green, $\Gamma_2$ is 
shown in blue, $\Gamma_3$ is shown in yellow, $\Gamma_4$ is shown in gray and $\Gamma_5$ is shown in red.  It is noted that $\Gamma_5$ consists of a single 
point located in the center of the top of the pressure vessel.  Since we are modeling a quarter of the pressure vessel, symmetry boundary conditions 
are applied on the displacements on $\Gamma_1$, $\Gamma_2$ and $\Gamma_5$.  These BCs emulate performing 
a simulation on the full cylindrical pressure vessel geometry (Figure \ref{fig:PCAPDBCs}(b)).  
\begin{figure}[h!tb]
 \centering
   \begin{tabular}{cc}
	   \includegraphics[width=0.3\textwidth]{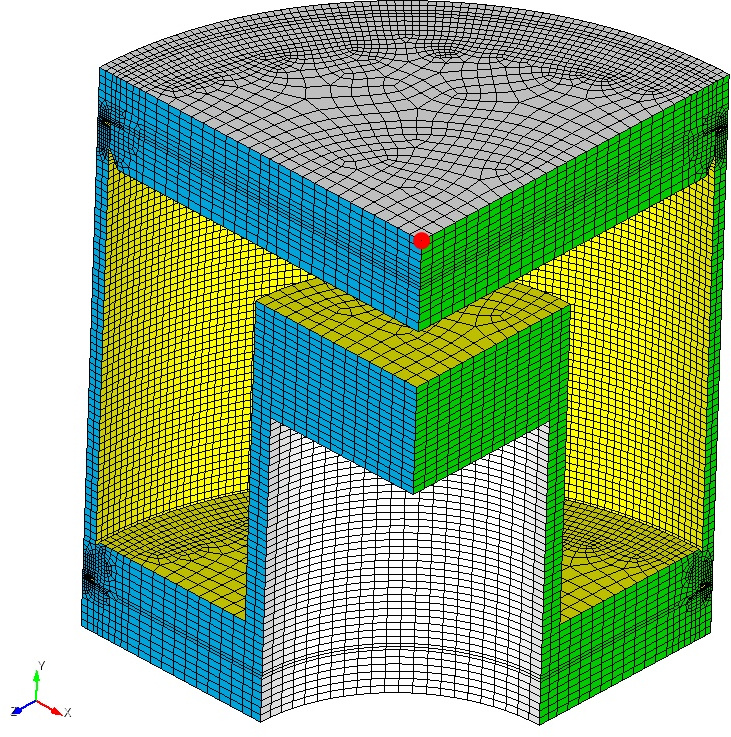} &  \hspace{2cm}
	   \includegraphics[width=0.3\textwidth]{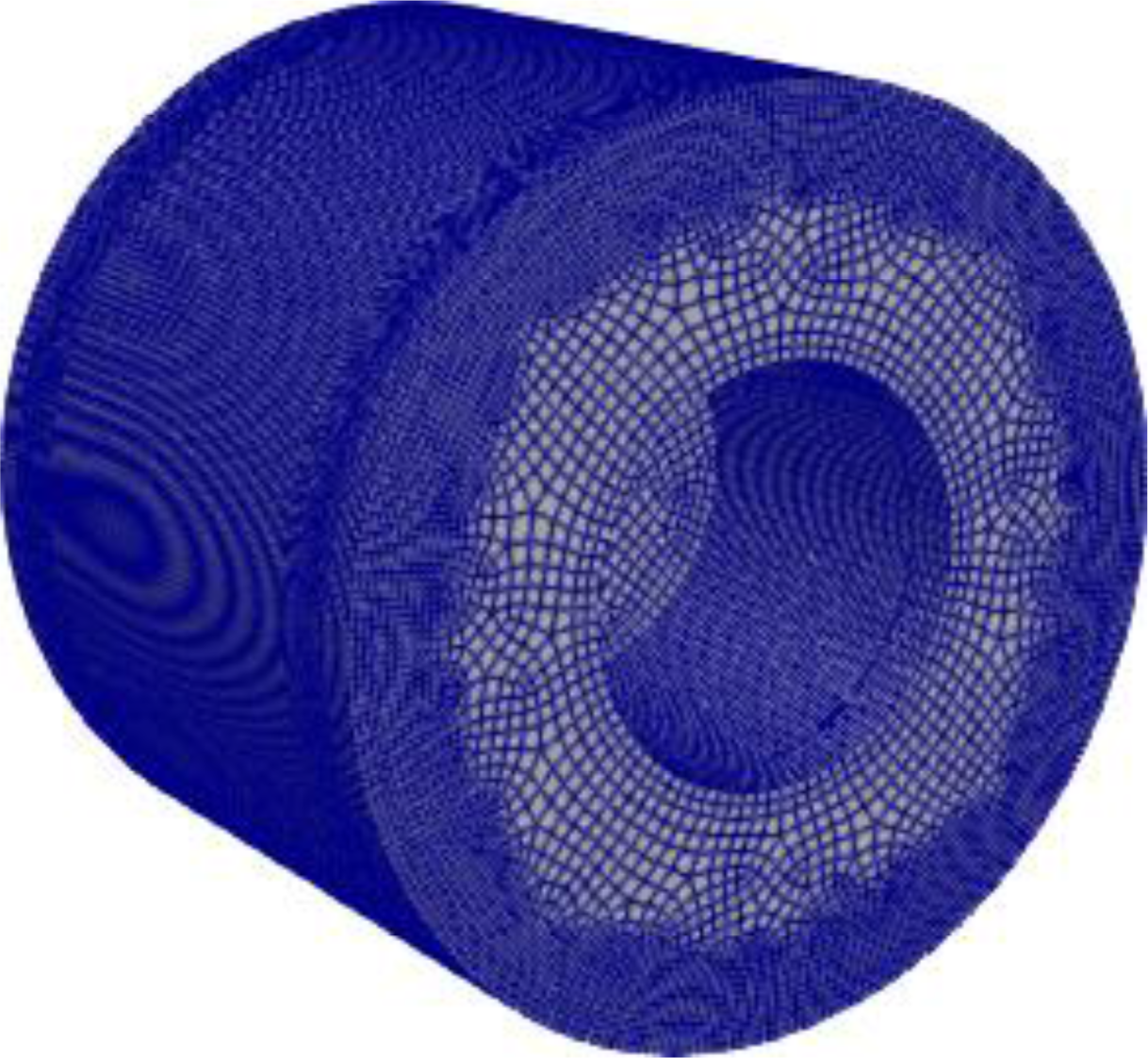}  \\ 
	   \\
	   (a) Boundaries & \hspace{2cm} (b) Full geometry with mesh\\
   \end{tabular}
	\caption{Thermo-mechanical pressure vessel problem.  (a) Illustration of boundaries on which boundary conditions are applied, 
	distinguished by different colors.  $\Gamma_1$ is shown in green, $\Gamma_2$ is
	shown in blue, $\Gamma_3$ is shown in yellow, $\Gamma_4$ is shown in gray and $\Gamma_5$ is shown in red. (b) Illustration
	of full pressure vessel geometry following application of symmetry boundary conditions on $\Gamma_1$, $\Gamma_2$ and $\Gamma_5$.}
  \label{fig:PCAPDBCs}
\end{figure}
In addition to the displacement boundary conditions mentioned previously, we impose the time-dependent 
pressure Neumann boundary condition \eqref{eq:beam_pbc} and 
a time-dependent temperature Dirichlet boundary condition \eqref{eq:beam_Tbc} on $\Gamma_3$.  This amounts to 
heating and pressurizing our vessel from the inside.
The boundary conditions account for a total of 5020 constrained dofs out of the 371,068 total dofs, leaving this problem with 366,048 free dofs.

Similar to the beam problem (Section \ref{sec:beam}), we  group together different element blocks for the purpose of 
specifying the material model parameters for the mechanical problem's material model, which we take to be Neohookean.  
Let $\mathcal{B}_a$ denote the union of blocks 3 and 5 (magenta and blue, respectively, in Figure \ref{fig:PCAPDBCs}) and let $\mathcal{B}_b$ denote the 
the union of blocks 1, 2 and 4 (green, yellow and cyan, respectively, in Figure \ref{fig:PCAPDBCs}).
We fix the mechanical parameters in block $\mathcal{B}_a$ to have the following values for the Young's modulus, Poisson's ratio, 
and density, respectively: $E_a=1.103\e{9}\press$, $\nu_a=0.32$, and density $\rho_a=7.92\e{-3}\dens$.  The thermal equation requires us 
to specify a reference temperature, which is set to $T_{a, \text{ref}} = 293 \temp$ within the thermal model.  
Similar to the thermo-mechanical version of the beam problem (Section \ref{sec:beam_thermomech}), these analogous parameters in block $\mathcal{B}_b$ 
are varied using LHC for the purpose of training and testing our ROMs; the values of these parameters
are given in in Table \ref{tab:ThermoMechPCAPPredicCases} and were sampled from the following ranges:
$E_b \in [1.27725 \times 10^9\press, 2.12875 \times 10^9\press]$, $\nu_b \in [0.24, 0.40]$, 
$\rho_b \in [5.94\times 10^{-3}\dens, 9.9\times 10^{-3}\dens]$,
$T_{b, \text{ref}} \in [219.75\temp, 366.25\temp]$.  
The following additional parameters are also 
specified within the thermal model in both blocks $\mathcal{B}_a$ and $\mathcal{B}_b$. 
As for the thermo-mechanical beam problem (Section \ref{sec:beam_thermomech}), 
we assume the materials are isotropic in both blocks ($\boldsymbol{K}:=k~I_3$), with thermal diffusivity $k_a = k_b = 1.1870\times 10^3$ m$^2$/s.
Additionally, we assume the thermal expansion coefficient \eqref{eq:FT} is the same in both blocks, and prescribe it the value $\alpha_a = \alpha_b = 1\e{-5}\temp^{-1}$.  
There are no applied forces, and the system is initially at rest with an initial temperature of $293\temp$.
The system is advanced using a quasi-static approach, in which the
a homotopy
continuation is performed with respect to the pseudo-time variable $t$ in \eqref{eq:beam_Tbc}.  The simulation 
proceeds for a total of 720 steps with a step size of $1\time$.  

\begin{figure}[h!tb]
 \centering
	\begin{tabular}{cc}
   \begin{tabular}{cc}
	   \includegraphics[width=0.15\textwidth]{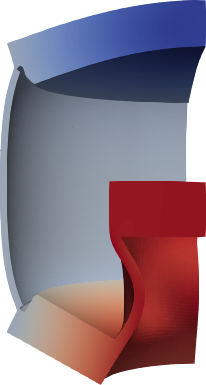} &  \hspace{2cm}
	   \includegraphics[width=0.15\textwidth]{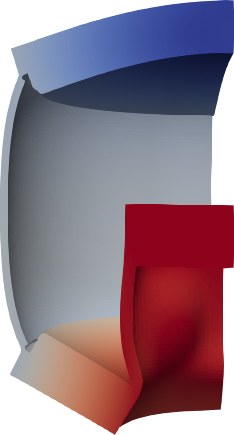}  \\ 
	   \\
	   (a) Testing case 1&  \hspace{2cm}(b) Testing case 2\\
   \end{tabular} \hspace{1cm} & 
		\includegraphics[width=0.09\textwidth]{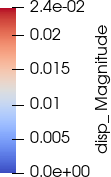}   
	\end{tabular}
	\caption{Thermo-mechanical pressure vessel problem: final displacement and configuration for the two testing cases in Table \ref{tab:ThermoMechPCAPPredicCases},
	obtained by running \albany{} in FOM mode.}
  \label{fig:PCAP_test_cases}
\end{figure}

In this example, we trained the ROM over four sets of parameters, summarized in the top part of 
Table \ref{tab:ThermoMechPCAPPredicCases}.  Snapshots were collected every continuation step, and 
POD bases of various numbers of modes, from 2 to 790, were constructed using our ensemble of 2160 snapshots. 
Once our ROMs were constructed, they were evaluated in the predictive regime on two additional cases, 
summarized in the bottom part of Table \ref{tab:ThermoMechPCAPPredicCases}. The final displacement FOM solutions
and configurations for the two testing cases are plotted in Figure \ref{fig:PCAP_test_cases}.  The reader
can observe that the solutions are noticeably different, indicating that our parameter variations 
led to nontrivial changes in the problem solution.




%
\begin{table}[H]
    \centering
	\caption{Thermo-mechanical pressure vessel problem: summary of parameters specified in the material model in block $\mathcal{B}_b$ 
	for the training and testing stages of the ROM process.}
    \label{tab:ThermoMechPCAPPredicCases}
    \begin{tabular}{r c | c c c c}
	    Regime & Case & $E_b(\e{9})\units{\press}$ & $\nu_b$ & $\rho_b(\e{-3})\units{\dens}$ & $T_{b, \text{ref}}\units{\temp}$ \\ \hline
        \multirow{4}{*}{training} & 1 & 1.64424 & 0.39524 & 8.33058 & 311.094 \\
        \multirow{4}{*}{} & 2 & 1.77118 & 0.300065 & 9.67843 & 267.396 \\ 
        \multirow{4}{*}{} & 3 & 1.9893 & 0.32161 & 7.17625 & 223.746 \\ 
        \multirow{4}{*}{} & 4 & 1.45551 & 0.266385 & 6.67746 & 331.116 \\ \hline
        \multirow{2}{*}{testing} & 1 & 2.06416 & 0.391368 & 7.79804 & 252.102 \\
        \multirow{2}{*}{} & 2 & 1.703 & 0.32 & 7.92 & 293 \\
    \end{tabular}
\end{table}

The main results for thermo-mechanical pressure vessel problem are summarized in 
Figures \ref{fig:PCAP_ROM_errors}--\ref{fig:PCAP_ROM_nonlin_iters}.  Similar to the thermo-mechanical
beam problem considered in Section \ref{sec:beam_thermomech}, the unpreconditioned LSPG ROMs fail to converge 
for the two larger basis dimensions evaluated.  The ideal preconditioned ROM solution is not plotted, 
as calculating this solution is too computationally expensive for this problem.
Figure \ref{fig:PCAP_ROM_errors} plots the global relative 
errors $\epsilon$ \eqref{pterr} as a function of the basis dimension for the two testing cases of interest.
The reader can observe that the unpreconditioned LSPG approach is unable to deliver a solution 
with a relative error smaller than $\mathcal{O}(10\%)$. In contrast, all three preconditioned LSPG 
ROMs produce solutions that are between four and six orders of magnitude more accurate.  
This result is confirmed by Figure \ref{fig:PCAP_ROM_sols}, which plots the final 
displacement and configuration obtained using the FOM and two LSPG ROMs having 32 modes, 
an unpreconditioned LSPG ROM and an LSPG ROM preconditioned using the {\tt Gauss-Seidel} 
preconditioner, for testing case 1.  It can be seen that the preconditioned LSPG ROM (Figure \ref{fig:PCAP_ROM_sols}(c))
is indistinguishable from the FOM (Figure \ref{fig:PCAP_ROM_sols}(a)); in contrast, the unpreconditioned
ROM solution (Figure \ref{fig:PCAP_ROM_sols}(b)) is visibly incorrect.
Turning our attention back to Figure \ref{fig:PCAP_ROM_errors}, it can be seen that,
by improving the preconditioner, it is possible to improve solution accuracy.
\begin{figure}[h!tb]
 \centering
	\begin{tabular}{cc}
	   \includegraphics[width=0.5\textwidth]{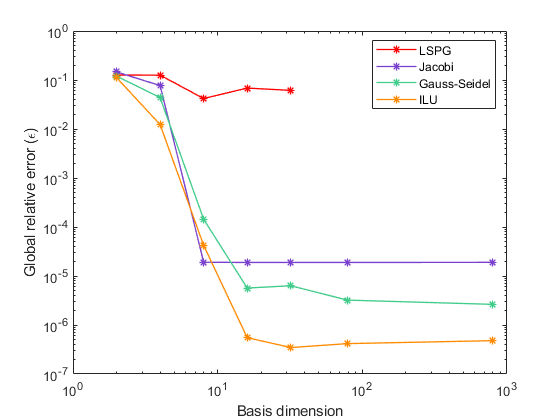} &  
	   \includegraphics[width=0.5\textwidth]{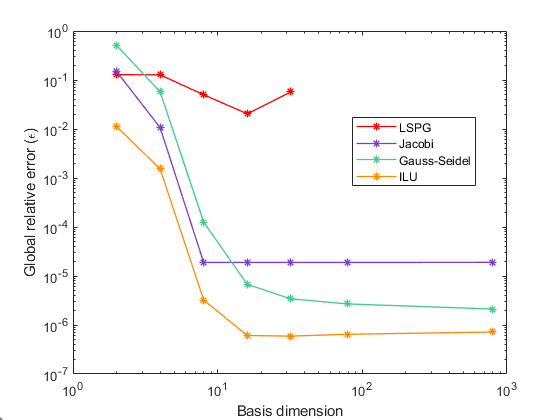} \\
	   (a) Testing case 1 & (b) Testing case 2  \\
	\end{tabular}
	\caption{Thermo-mechanical pressure vessel problem: global relative errors $\mathcal{\epsilon}$ \eqref{pterr} for various ROMs as a function
	of the basis dimension for the two testing cases in Table \ref{tab:ThermoMechPCAPPredicCases}.}
  \label{fig:PCAP_ROM_errors}
\end{figure}
\begin{figure}[h!tb]
 \centering
	\begin{tabular}{cc}
   \begin{tabular}{cc}
	   \includegraphics[width=0.15\textwidth]{fom_test2.png} &  
	   \includegraphics[width=0.20\textwidth]{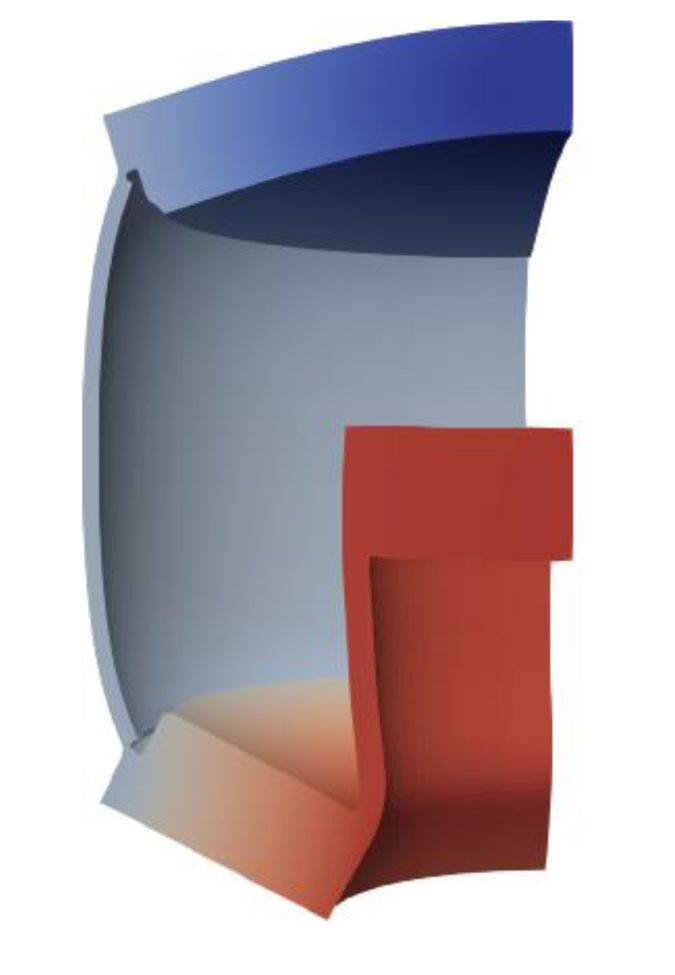} \\
	   (a) FOM & (b) LSPG ROM  \\
	   \\
	    \multicolumn{2}{c}{\includegraphics[width=0.20\textwidth]{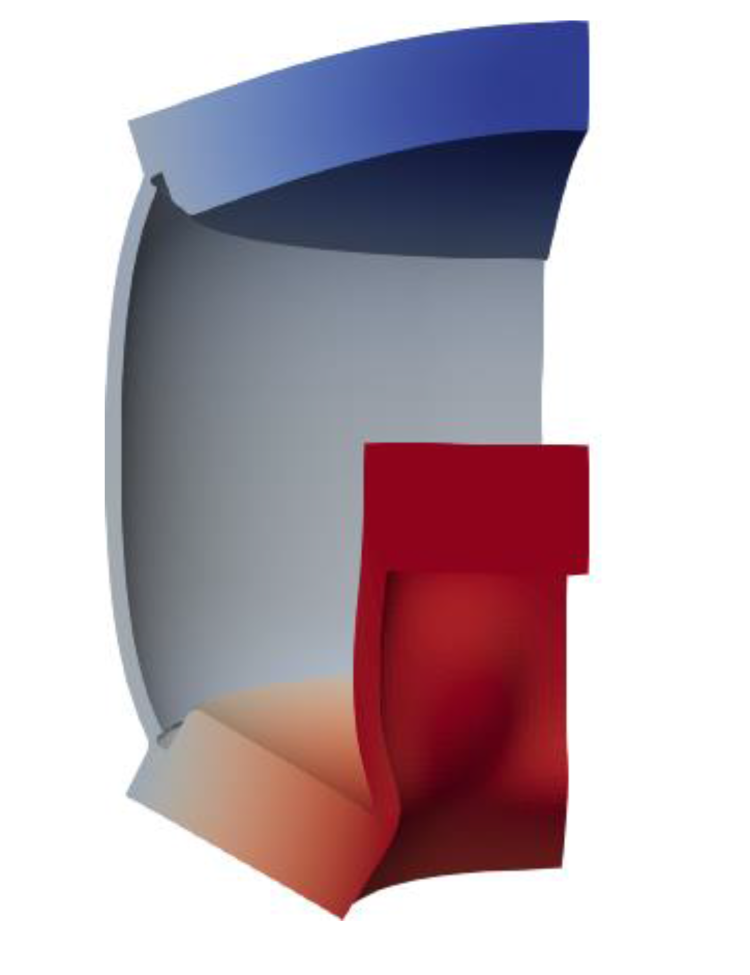}}  
	   \\
	   \multicolumn{2}{c}{(c) {\tt Gauss-Seidel} Preconditioned LSPG ROM} \\
   \end{tabular} & 
		\includegraphics[width=0.09\textwidth]{legend.png}   
	\end{tabular}
	\caption{Thermo-mechanical pressure vessel problem: comparison of FOM and $\nrb=32$ ROM solutions (final displacement and configuration) for testing case 1.}
  \label{fig:PCAP_ROM_sols}
\end{figure}
Figure \ref{fig:PCAP_ROM_times} depicts the wall times required to run the various LSPG ROMs evaluated
for each of the testing cases.  Our simulations for the thermo-mechanical pressure vessel problem
were performed on 64 cores (4 nodes) of the Skybridge high-capacity cluster located at Sandia National Laboratories, which 
contains 1848 nodes, each having 16 2.6 GHz Intel Sandy Bridge processors.  It can be seen from this figure 
that a speedup of as large as $12\times$ is attainable through the introduction of preconditioning
into the LSPG ROM formulation.  
Figures \ref{fig:PCAP_ROM_errors} and \ref{fig:PCAP_ROM_times} together demonstrate that LSPG ROMs are substantially
faster and more accurate than unpreconditioned LSPG ROMs.  This result is reinforced by Figure  
\ref{fig:pcap_pareto}, which depicts a Pareto plot and front for the thermo-mechanical pressure vessel problem.
As before, it is possible to use this figure to 
determine the ideal preconditioner to use based on one's error and CPU-time requirements.

\begin{figure}[h!tb]
 \centering
	\begin{tabular}{cc}
	   \includegraphics[width=0.5\textwidth]{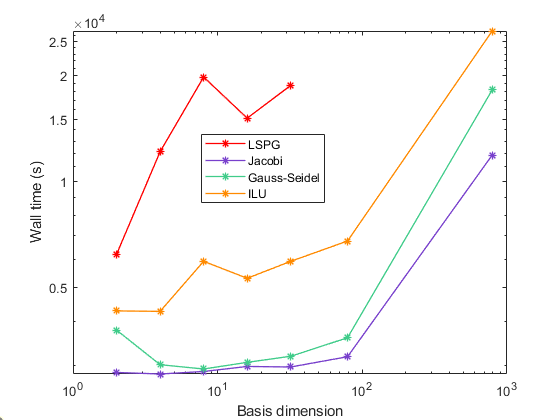} &  
	   \includegraphics[width=0.5\textwidth]{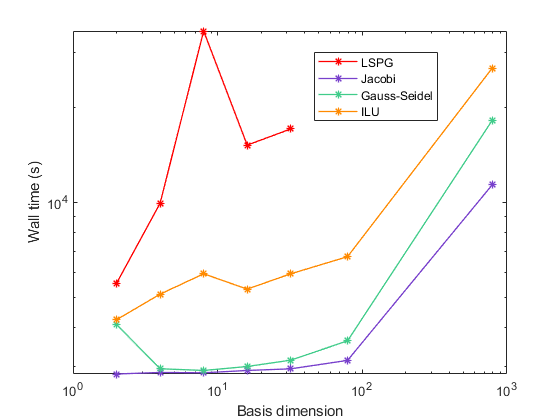} \\
	   (a) Testing case 1 & (b) Testing case 2  \\
	\end{tabular}
	\caption{Thermo-mechanical pressure vessel problem: wall times (in$\time$, averaged over 64 processors) for various ROMs as a function of the basis dimension for the two testing cases in Table \ref{tab:ThermoMechPCAPPredicCases}.}
  \label{fig:PCAP_ROM_times}
\end{figure}

\begin{figure}[h!tb]
 \centering
\includegraphics[width=1.0\textwidth]{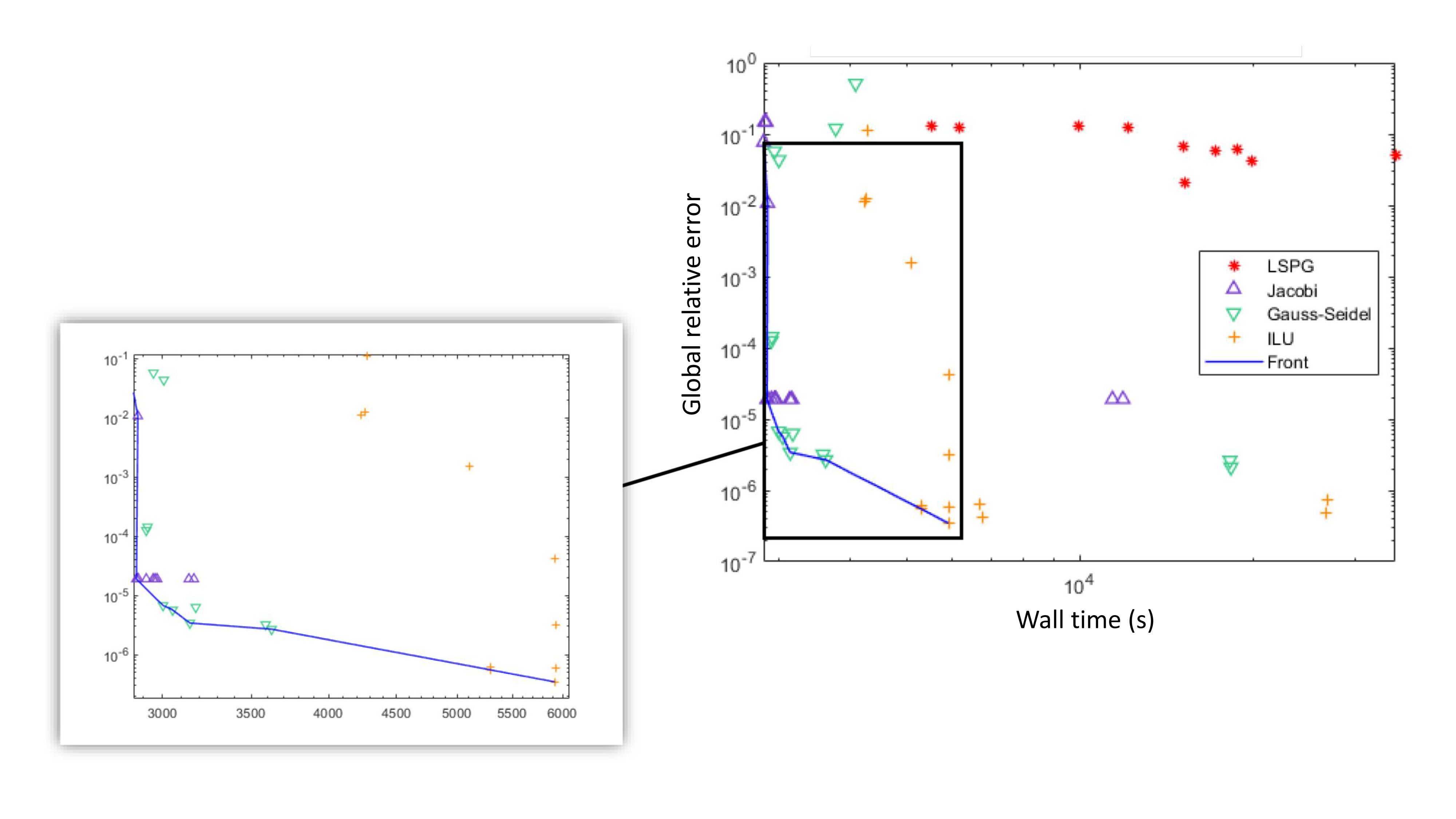}
	\caption{Thermo-mechanical pressure vessel problem: Pareto plot showing the total wall time (in s, averaged over 4 processors) 
	versus the global relative error $\epsilon$.  Results for all four testing cases in Table
	\ref{tab:ThermoMechBeamPredicCases} are used to generate this figure.} 
  \label{fig:pcap_pareto}
\end{figure}

Our final figure, Figure \ref{fig:PCAP_ROM_nonlin_iters}, plots the total number of nonlinear iterations required 
to attain a convergent solution for the various ROMs evaluated in this section.  We present this result in lieu of 
the average condition number of the reduced Jacobian $\PrecondPetrovGalerkin{\jac}\newton{k}$, as the thermo-mechanical 
pressure vessel problem is too large to calculate the latter quantity.  Comparing Figure \ref{fig:PCAP_ROM_nonlin_iters} 
with Figure \ref{fig:PCAP_ROM_times}, one can see that the wall time improvements achieved through the introduction 
of preconditioning are likely due to the reduction in the total number of nonlinear iterations
required to achieve a convergent ROM solution: by preconditioning the LSPG formulation, one can reduce the number 
of nonlinear iterations by up to a factor of $12\times$.  It is interesting to observe that, unlike 
for the classical LSPG formulation, the number of nonlinear 
iterations remains approximately constant regardless of the basis dimension for the preconditioned LSPG ROMs.
The reduction in the number of nonlinear iterations when preconditioning is employed 
can be attributed to the scaling effect of the ROM preconditioner for this multi-scale and multi-physics 
problem. 

\begin{figure}[h!tb]
 \centering
	\begin{tabular}{cc}
	   \includegraphics[width=0.5\textwidth]{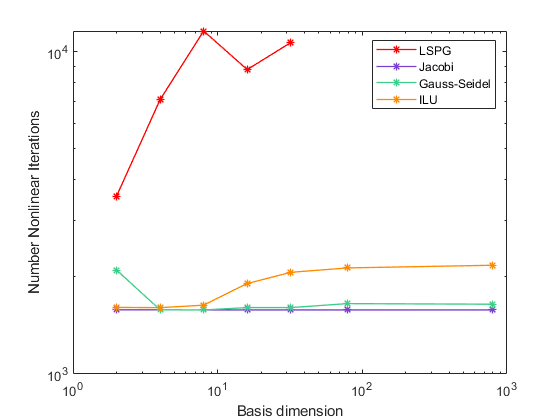} &  
	   \includegraphics[width=0.5\textwidth]{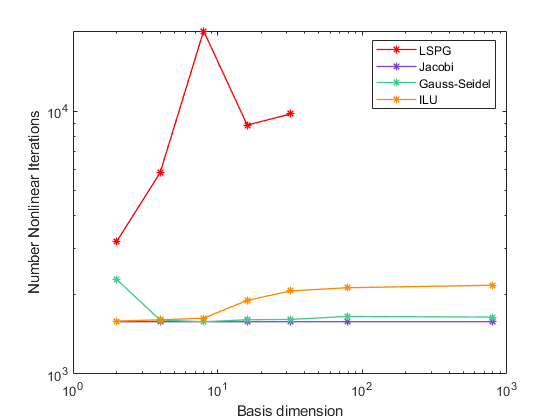} \\
	   (a) Testing case 1 & (b) Testing case 2  \\
	\end{tabular}
	\caption{Thermo-mechanical pressure vessel problem: total number of nonlinear iterations for various ROMs as a function of the basis dimension for the two testing cases in Table \ref{tab:ThermoMechPCAPPredicCases}.}
  \label{fig:PCAP_ROM_nonlin_iters}
\end{figure}

\section{Conclusions}
\label{sec:conclusion}

Projection-based model reduction can enable real-time and multi-query analyses
in a variety of applications.  
The LSPG method to building projection-based ROMs
has shown particular promise, as it can generate stable and accurate solutions
for applications where standard Galerkin techniques have failed \cite{CarlbergGappy,kevin_gnat, fy16_sand, fy17_sand, fy18_sand,
carlbergGalDiscOpt}.  Despite its superior performance 
over Galerkin projection, we have observed that, when run in the predictive regime,
the LSPG method sometimes exhibits a lack of convergence
or delivers a solution whose error surpasses the engineering tolerance required by the application.
The approach is most likely to suffer from these deficiencies when applied to problems with disparate
scales, such as dimensional PDEs or multi-physics problems \cite{kyle, Blonigan:2022, Parish:2022}.  This is because, when there are differences 
in scale between components of a PDE system, residual components corresponding to certain variables 
can be very large compared to the residual components corresponding to other variables, which can bias the 
residual minimization procedure underlying the LSPG formulation.  

In this paper, we have demonstrated that the accuracy, robustness, efficiency and convergence properties of 
the LSPG method for projection-based model reduction can be improved substantially through the introduction
of preconditioning.  While this approach can reduce the condition number 
of the matrix problems arising in the LSPG iteration procedure, it is unlike other existing ROM preconditioning
approaches (e.g., \cite{elman1, elman2,singh1}), as the preconditioning is introduced directly into the LSPG residual minimization problem.
Doing so ensures that all residual components are on approximately the same scale, even 
when the method is applied to multi-scale and multi-physics problems having variables of drastically different 
magnitudes.
We demonstrated that preconditioning the LSPG method in this way is equivalent to modifying the norm in which 
residual minimization is performed, which can improve the
residual-based stability constant bounding the method's error.  
We additionally showed that, by designing 
a preconditioner that approximates the inverse of the Jacobian matrix, it is possible to 
create an LSPG ROM whose accuracy approaches the upper limit on ROM accuracy for a given reduced basis:
an ideal preconditioned ROM, which is equivalent to 
the projection of the FOM solution increment onto the reduced basis.  

We studied the efficacy of the proposed approach using one mechanical and two thermo-mechanical 
examples simulated using the \albany{} HPC multi-physics finite element code \cite{albany}, into which
we have implemented the POD/LSPG method to model reduction.  
Our test cases include a realistic simulation involving a thermo-mechanical vessel
pressurized from the inside, having approximately 370,000 dofs.
We developed a partitioned/``blocking vector'' approach for applying ROM 
Dirichlet boundary conditions strongly within this code, 
which does not remove the constrained dofs from
the global finite element system prior to performing its numerical
solution.  We considered three simple 
preconditioners, a Jacobi preconditioner, a Gauss-Seidel preconditioner, and 
an ILU preconditioner, implemented within the \ifpack{} \trilinos{} library,
on which \albany{} has a dependence.
Our ROMs were evaluated in the predictive regime,
with prediction performed across the material parameter space.  
For the first test case, known as the mechanical beam problem, 
the classical LSPG method failed to deliver a convergent solution regardless of the basis dimension.
In contrast, all three preconditioned ROMs converged robustly and efficiently, 
with global relative errors as small as $\mathcal{O}(0.01\%)$.  
While the classical LSPG method achieved convergence for the smaller basis dimensions
considered in our two thermo-mechanical test cases, we demonstrated that the introduction of 
preconditioning can reduce the global relative error by up to six orders of magnitude and the wall time by as much as $12\times$.  
The wall time gains are attributed to a reduction in the reduced Jacobian condition number and
a reduction in the total number 
of nonlinear iterations required for convergence, both achieved via the proposed preconditioning
strategy.  For thermo-mechanical problems, which exhibit extreme differences in scale 
between the displacement and temperature solutions, the introduction of even a simple preconditioner
such as Jacobi reduced the reduced Jacobian condition number by as many as ten orders of magnitude.

The results discussed herein pave the way for several subsequent studies, which will be the subject
of future work.  First, we plan to extend our preconditioned LSPG formulation to include hyper-reduction using 
Lagrangian structure-preserving methods such as \cite{FarhatChapman:2015, Chapman:2017, FarhatCortial:2014, 
CarlbergAIAA:2012, Carlberg:2015}, and 
intend to repeat the present numerical experiments with hyper-reduction enabled.  It would 
additionally be 
interesting to study: (1) the efficacy of more sophisticated preconditioners within our framework, e.g.,  
 multi-grid preconditioners that take advantage 
of rigid body mode information and block preconditioners for multi-physics applications, 
and (2) the performance of our framework with more sophisticated reduced bases, e.g., 
nonlinear manifold reduced bases constructed using an auto-encoder or some other type of 
neural network \cite{fresca, kookjin, youngsoo}. 
Lastly, we plan to investigate the viability of our approach for a larger set of applications, 
including problems in the field of Computational Fluid Dynamics (CFD) and problems in which 
the ROM is used to perform prediction in time.  Our preliminary results \cite{fy18_sand, pers_corr}  suggest that preconditioning improves accuracy and efficiency in 
problems involving compressible cavity flows, hypersonic aerodynamics, thermal ablation and reacting hypersonic flows.

\section{Acknowledgements}
\label{sec:acknowledgements}

Support for this work was received through the U.S. Department of Energy's (DOE)
Advanced Simulation and Computing (ASC) Program at Sandia National
Laboratories. The writing of this manuscript was funded in part by the third author’s Presidential Early Career Award for Scientists and Engineers (PECASE).
Sandia National Laboratories is a multi-mission laboratory
managed and operated by National Technology and Engineering Solutions of Sandia, LLC.,
a wholly owned subsidiary of Honeywell International, Inc., for the U.S. Department of Energy's
National Nuclear Security Administration under contract DE-NA0003525.

The numerical results presented herein were generated by running the \albany{} open-source
HPC code.  The version of the code that was used can be downloaded from the following URL:
\url{https://github.com/sandialabs/Albany/releases/tag/MOR_support_end}.  Additionally,
to ensure transparency and reproducibility, we have made available the \albany{} input files
needed to  reproduce our results.  These input files 
can be downloaded from the following {\tt github} repository: \url{https://github.com/ikalash/PrecondLSPGROMsSM}.   

The authors wish to thank Alejandro Mota for assisting with the formulation
of the mechanical and thermo-mechanical problems described herein.  The authors would 
also like to thank Eric Parish for providing useful feedback on the first full 
draft of this work, which resulted in an improved manuscript.


\bibliography{ms}{}

\begin{thebibliography}{10}

\bibitem{Lall:2003}
Structure-preserving model reduction for mechanical systems.
\newblock {\em Physica D: Nonlinear Phenomena}, 184(1):304--318, 2003.

\bibitem{anasazi}
Anasazi: a block eigensolvers package.
\newblock
  \url{https://docs.trilinos.org/dev/packages/anasazi/doc/html/index.html},
  2020.

\bibitem{albany_repo}
{MOR} support end tag.
\newblock
  \url{https://github.com/SNLComputation/Albany/releases/tag/MOR_support_end},
  2020.

\bibitem{loca}
{NOX and LOCA: object-oriented nonlinear solver and continuation packages}.
\newblock \url{https://docs.trilinos.org/dev/packages/nox/doc/html}, 2020.

\bibitem{Abgrall:2018}
R.~Abgrall and R.~Crisovan.
\newblock Model reduction using l1-norm minimization as an application to
  nonlinear hyperbolic problems.
\newblock {\em International Journal for Numerical Methods in Fluids},
  87(12):628--651, 2018.

\bibitem{ahmad}
M.~Ahmad, D.~Szyld, and M.~van Guzen.
\newblock Preconditioned multishift {BiCG} for {H2}-optimal model reduction.
\newblock {\em SIAM J. Matrix Anal. Appl.}, 38:401--424, 2017.

\bibitem{anzt}
H.~Anzt, E.~Chow, J.~Saak, and J.~Dongarra.
\newblock Updating incomplete factorization preconditioners for model order
  reduction.
\newblock {\em Numer. Algor.}, 73:611--630, 2016.

\bibitem{astrid}
P.~Astrid, S.~Weiland, K.~Willcox, and T.~Backx.
\newblock Missing point estimation in models described by proper orthogonal
  decomposition.
\newblock {\em {IEEE} Transactions on Automatic Control}, 53(10):2237--2250,
  2008.

\bibitem{balabanov}
O.~Balabanov.
\newblock Randomized linear algebra for model order reduction.
\newblock Ph.D. Thesis, Universitat Politecnica de Catalunya, 2019.

\bibitem{maciejirina}
M.~Balajewicz, I.~Tezaur, and E.~Dowell.
\newblock Minimal subspace rotation on the stiefel manifold for stabilization
  and enhancement of projection-based reduced order models for the compressible
  navier-stokes equations.
\newblock {\em J. Comput. Phys.}, 321:224--241, 2016.

\bibitem{maciej}
M.J. Balajewicz, E.H. Dowell, and B.R. Noack.
\newblock Low-dimensional modelling of high-{R}eynolds-number shear flows
  incorporating constraints from the {N}avier-{S}tokes equation.
\newblock {\em J. Fluid Mech.}, 729:285--308, 2013.

\bibitem{romLDRDSAND}
M.F. Barone, I.~Kalashnikova, M.R. Brake, and D.J. Segalman.
\newblock Reduced order modeling of fluid/structure interaction.
\newblock Sandia National Laboratories Report, SAND No. 2009-7189, Sandia
  National Laboratories, Albuquerque, NM, 2009.

\bibitem{baroneJCP}
M.F. Barone, I.~Kalashnikova, D.J. Segalman, and H.~Thornquist.
\newblock Stable {G}alerkin reduced order models for linearized compressible
  flow.
\newblock {\em J. Comput. Phys.}, 288:1932--1946, 2009.

\bibitem{pers_corr}
P.~Blonigan.
\newblock {Personal correspondence}, 2021.

\bibitem{carlbergGalDiscOpt}
K.~Carlberg, M.~Barone, and H.~Antil.
\newblock Galerkin v.\ least-squares {P}etrov--{G}alerkin projection in
  nonlinear model reduction.
\newblock {\em J. Comput. Phys.}, 330:693--734, 2017.

\bibitem{carlbergHawaii}
K.~Carlberg, J.~Cortial, D.~Amsallem, M.~Zahr, and C.~Farhat.
\newblock The {GNAT} nonlinear model reduction method and its application to
  fluid dynamics problems.
\newblock {\em AIAA Paper 2011-3112, 6th AIAA Theoretical Fluid Mechanics
  Conference, Honolulu, HI}, June 27--30, 2011.

\bibitem{CarlbergGappy}
K.~Carlberg, C.~Farhat, and C.~Bou-Mosleh.
\newblock {Efficient non-linear model reduction via a least-squares
  Petrov--Galerkin projection and compressive tensor approximations}.
\newblock {\em International Journal for Numerical Methods in Engineering},
  86(2):155--181, April 2011.

\bibitem{kevin_gnat}
K.~Carlberg, C.~Farhat, J.~Cortial, and D.~Amsallem.
\newblock The {GNAT} method for nonlinear model reduction: Effective
  implementation and application to computational fluid dynamics and turbulent
  flows.
\newblock {\em J. Comput. Phys.}, 24(2):623--647, 2013.

\bibitem{carlbergKrylov2015}
K.~Carlberg, V.~Forstall, and R.~Tuminaro.
\newblock Krylov-subspace recycling via the {POD}-augmented conjugate-gradient
  method.
\newblock {\em SIAM Journal on Matrix Analysis and Applications},
  37(3):1304--1336, 2016.

\bibitem{CarlbergAIAA:2012}
K.~Carlberg, R.~Tuminaro, and P.~Boggs.
\newblock Efficient structure-preserving model reduction for nonlinear
  mechanical systems with application to structural dynamics.
\newblock {\em AIAA Paper 2012-1969, 53rd AIAA/ASME/ASCE/AHS/ASC Structures,
  Structural Dynamics and Materials Conference, Honolulu, HI}, April 23--26,
  2012.

\bibitem{kevin_lagrangian}
K.~Carlberg, R.~Tuminaro, and P.~Boggs.
\newblock Preserving {L}agrangian structure in nonlinear model reduction with
  application to structural dynamics.
\newblock {\em SIAM J. Sci. Comput.}, 2:B153--B184, 2015.

\bibitem{Carlberg_SIAM_2015}
K.~Carlberg, R.~Tuminaro, and P.~Boggs.
\newblock Preserving lagrangian structure in nonlinear model reduction with
  application to structural dynamics.
\newblock {\em SIAM J. Sci. Comput}, 47(2):B153--B184, 2015.

\bibitem{Carlberg:2015}
K.~Carlberg, R.~Tuminaro, and P.~Boggs.
\newblock Preserving {L}agrangian structure in nonlinear model reduction with
  application to structural dynamics.
\newblock {\em SIAM Journal on Scientific Computing}, 37(2):B153--B184, 2015.

\bibitem{Chapman:2017}
Todd Chapman, Philip Avery, Pat Collins, and Charbel Farhat.
\newblock Accelerated mesh sampling for the hyper reduction of nonlinear
  computational models.
\newblock {\em International Journal for Numerical Methods in Engineering},
  109(12):1623--1654, 2017.

\bibitem{deim}
S.~Chaturantabut and D.C. Sorensen.
\newblock Nonlinear model reduction via discrete empirical interpolation.
\newblock {\em SIAM J. Sci. Comput.}, 32:2737--2764, 2010.

\bibitem{Blonigan:2022}
D.~Ching, P.~Blonigan, F.~Rizzi, and J.~Fike.
\newblock Reduced order modeling of hypersonic aerodynamics with grid
  tailoring.
\newblock AIAA SciTech Forum, January 3-7, 2022, San Diego, CA and virtual,
  2022.

\bibitem{kevin_space-time}
Y.~Choi and K.~Carlberg.
\newblock Space-time least-squares petrov-galerkin projection for nonlinear
  model reduction.
\newblock {\em arXiv e-Print}, (1703.04560), 2017.

\bibitem{Collins:2020}
Gary Collins, Krzysztof Fidkowski, and Carlos~E. Cesnik.
\newblock {\em Petrov-Galerkin Projection-Based Model Reduction with an
  Optimized Test Space}.

\bibitem{elman2}
H.~Elman and Virginia Forstall.
\newblock Preconditioning techniques for reduced basis methods for
  parameterized elliptic partial differential equations.
\newblock {\em SIAM J. Sci. Comput.}, 37, 2015.

\bibitem{elman1}
Howard~C. Elman and Virginia Forstall.
\newblock Numerical solution of the parameterized steady-state navier-stokes
  equations using empirical interpolation methods.
\newblock {\em Computer Methods in Applied Mechanics and Engineering}, 317(C),
  12 2016.

\bibitem{gappy}
R.~Everson and L.~Sirovich.
\newblock Karhunen-{L}oeve procedure for gappy data.
\newblock {\em J. Optical Society of America A}, pages 1657--1664, 1995.

\bibitem{FarhatCortial:2014}
Charbel Farhat, Philip Avery, Todd Chapman, and Julien Cortial.
\newblock Dimensional reduction of nonlinear finite element dynamic models with
  finite rotations and energy-based mesh sampling and weighting for
  computational efficiency.
\newblock {\em International Journal for Numerical Methods in Engineering},
  98(9):625--662, 2014.

\bibitem{FarhatChapman:2015}
Charbel Farhat, Todd Chapman, and Philip Avery.
\newblock Structure-preserving, stability, and accuracy properties of the
  energy-conserving sampling and weighting method for the hyper reduction of
  nonlinear finite element dynamic models.
\newblock {\em International Journal for Numerical Methods in Engineering},
  102(5):1077--1110, 2015.

\bibitem{ACE_paper}
Jennifer Frederick, Alejandro Mota, Irina Tezaur, and Diana Bull.
\newblock A thermo-mechanical terrestrial model of arctic coastal erosion.
\newblock {\em Journal of Computational and Applied Mathematics}, 397:113533,
  2021.

\bibitem{fresca}
S.~Fresca, L.~Dede, and A.~Manzoni.
\newblock {A Comprehensive Deep Learning-Based Approach to Reduced Order
  Modeling of Nonlinear Time-Dependent Parametrized PDEs}.
\newblock {\em Journal of Scientific Computing}, 87(61), 2021.

\bibitem{QCAD:2013}
X.~Gao, E.~Nielsen, R.~Muller, R.~Young, A.~Salinger, N.~Bishop, M.~Lilly, and
  M.~Carroll.
\newblock Quantum computer aided design simulation and optimization of
  semiconductor quantum dots.
\newblock {\em Journal of Applied Physics}, 114:1--19, 2013.

\bibitem{gugercin}
S.~Gugercin and A.C. Antoulas.
\newblock A survey of model reduction by balanced truncation and some new
  results.
\newblock {\em Int. J. Control}, 77(8):748--766, 2004.

\bibitem{aztec}
M.~Heroux.
\newblock {AztecOO User Guide}.
\newblock Sandia Report, SAND2014-3796, 2014.

\bibitem{trilinos}
M.~Heroux, R.~Bartlett, V.~Howle, R.~Hoekstra, J.~Hu, T.~Kolda, R.~Lehoucq,
  K.~Long, R.~Pawlowski, E.~Phippsand~A. Salinger, H.~Thornquist, R.~Tuminaro,
  J.~Willenbring, A.~Williams, and K.~Stanley.
\newblock An overview of the trilinos project.
\newblock {\em {ACM} Trans. Math. Softw.}, 31:397--423, 2005.

\bibitem{hoang2021domain}
Chi Hoang, Youngsoo Choi, and Kevin Carlberg.
\newblock Domain-decomposition least-squares {P}etrov--{G}alerkin ({DD-LSPG})
  nonlinear model reduction.
\newblock {\em Computer Methods in Applied Mechanics and Engineering},
  384:113997, 2021.

\bibitem{holmes}
P.~Holmes, J.L. Lumley, and G.~Berkooz.
\newblock {\em Turbulence, Coherent Structures, Dynamical Systems and
  Symmetry}.
\newblock Cambridge University Press, 1996.

\bibitem{Holzapfel}
G.~Holzapfel.
\newblock {\em Nonlinear solid mechanics: a continuum approach for
  engineering}.
\newblock Wiley, 1st edition, 2000.

\bibitem{jiang}
R.~Jiang.
\newblock Pressure preconditioning using proper orthogonal decomposition.
\newblock M.S. Thesis, Stanford University, Stanford, CA (Advisor: H.
  Tchelepi), 2013.

\bibitem{ec_ldrd_sand}
I.~Kalashnikova, S.~Arunajatesan, M.F. Barone, B.G. van Bloemen~Waanders, and
  J.A. Fike.
\newblock Reduced order modeling for prediction and control of large-scale
  systems.
\newblock Sandia National Laboratories Report, SAND No. 2014-4693. Sandia
  National Laboratories, Albuquerque, NM, 2014.

\bibitem{baronekalash}
I.~Kalashnikova and M.F. Barone.
\newblock On the stability and convergence of a {G}alerkin reduced order model
  ({ROM}) of compressible flow with solid wall and far-field boundary
  treatment.
\newblock {\em J. Comput. Phys.Int. J. Numer. Meth. Engng.}, 83:1345--1375,
  2010.

\bibitem{youngsoo}
Y.~Kim, D.~Widemann, Y.~Choi, and T.~Zohdi.
\newblock Efficient nonlinear manifold reduced order model.
\newblock ArXiV pre-print, \url{https://arxiv.org/pdf/2011.07727.pdf}, 2021.

\bibitem{lall}
Sanjay Lall, Petr Krysl, and Jerrold~E Marsden.
\newblock Structure-preserving model reduction for mechanical systems.
\newblock {\em Physica D: Nonlinear Phenomena}, 184(1):304--318, 2003.
\newblock Complexity and Nonlinearity in Physical Systems -- A Special Issue to
  Honor Alan Newell.

\bibitem{kookjin}
Kookjin Lee and Kevin~T. Carlberg.
\newblock Model reduction of dynamical systems on nonlinear manifolds using
  deep convolutional autoencoders.
\newblock {\em Journal of Computational Physics}, 404:108973, 2020.

\bibitem{collocation}
P.~LeGresley.
\newblock Application of proper orthogonal decomposition (pod) to design
  decomposition methods.
\newblock Ph.D. thesis, Stanford University, 2006.

\bibitem{moore}
B.~Moore.
\newblock Principal component analysis in linear systems: Controllability,
  observability, and model reduction.
\newblock {\em IEEE Transactions on Automatic Control}, 26(1), 1981.

\bibitem{Mota:2022}
A.~Mota, I.~Tezaur, and G.~Phlipot.
\newblock The {S}chwarz alternating method for dynamic solid mechanics.
\newblock Int. J. Numer. Meth. Engng (in press), 2022.

\bibitem{schwarz}
Alejandro Mota, Irina Tezaur, and Coleman Alleman.
\newblock The schwarz alternating method in solid mechanics.
\newblock {\em Computer Methods in Applied Mechanics and Engineering}, 319:19
  -- 51, 2017.

\bibitem{bp1}
N.C. Nguyen, A.T. Patera, and J.~Peraire.
\newblock A `best points' interpolation method for efficient approximation of
  parametrized functions.
\newblock {\em Int. J. Numer. Meth. Engng.}, 73:521--543, 2008.

\bibitem{bp2}
N.C. Nguyen and J.~Peraire.
\newblock An efficient reduced-order modeling approach for non-linear
  parametrized partial differential equations.
\newblock {\em Int. J. Numer. Meth. Engng.}, 76:27--55, 2008.

\bibitem{noack_papas_monkewitz_2005}
B.~Noack, R.~Papas, and P.~Monkewitz.
\newblock The need for a pressure-term representation in empirical galerkin
  models of incompressible shear flows.
\newblock {\em Journal of Fluid Mechanics}, 523:339--365, 2005.

\bibitem{Parish:2022}
E.~Parish and F.~Rizzi.
\newblock On the impact of dimensionally-consistent and physics-based inner
  products for pod-galerkin and least-squares model reduction of compressible
  flows.
\newblock ArXiV pre-print, 2022.

\bibitem{parish2021windowed}
Eric~J Parish and Kevin~T Carlberg.
\newblock Windowed least-squares model reduction for dynamical systems.
\newblock {\em Journal of Computational Physics}, 426:109939, 2021.

\bibitem{pasetto}
Damiano Pasetto, Massimiliano Ferronato, and Mario Putti.
\newblock A reduced order model-based preconditioner for the efficient solution
  of transient diffusion equations.
\newblock {\em International Journal for Numerical Methods in Engineering},
  109(8):1159--1179, 2017.

\bibitem{epetraext_me}
R.~Pawlowski, R.~Bartlett, N.~Belcourt, R.~Hooper, and R.~Schmidt.
\newblock A therory manual for multi-physics code coupling in {LIME}.
\newblock Sandia National Laboratories Report, SAND No. 2011-2195, Sandia
  National Laboratories, Albuquerque, NM, 2011.

\bibitem{Rempfer_JFM_1994}
D.~Rempfer and H.F. Fasel.
\newblock Dynamics of three-dimensional coherent structures in a flat-plate
  boundary layer.
\newblock {\em Journal of Fluid Mechanics}, 275:257--283, 1994.

\bibitem{rowley}
C.W. Rowley.
\newblock Model reduction for fluids using balanced proper orthogonal
  decomposition.
\newblock {\em Int. J. Bif. Chaos}, 15(3):997--1013, 2005.

\bibitem{rowley1}
C.W. Rowley, T.~Colonius, and R.M. Murray.
\newblock Model reduction for compressible flows using {POD} and {G}alerkin
  projection.
\newblock {\em Physica D}, 189:115--129, 2004.

\bibitem{rozza}
G.~Rozza.
\newblock Reduced basis approximation and error bounds for potential flows in
  parametrized geometries.
\newblock {\em Commun. Comput. Phys.}, 9(1):1--48, 2011.

\bibitem{ifpack-guide}
M.~Sala and M.~Heroux.
\newblock Robust algebraic preconditioners with {IFPACK} 3.0.
\newblock Technical Report SAND-0662, Sandia National Laboratories, 2005.

\bibitem{albany}
A.~Salinger, R.~Bartlett, A.~Bradley, Q.~Chen, I.~Demeshko, X.~Gao, G.~Hansen,
  A.~Mota, R.~Muller, E.~Nielsen, J.~Ostien, R.~Pawlowski, M.~Perego,
  E.~Phipps, W.~Sun, and I.~Tezaur.
\newblock Albany: Using agile components to develop a flexible, generic
  multiphysics analysis code.
\newblock {\em Int. J. Multiscale Comput. Engng}, 14:415--438, 2016.

\bibitem{AgileComponents}
A.~Salinger, I.~Tezaur, M.~Perego, et~al.
\newblock Component-based application code development, part 1: The agile
  components strategy and albany code.
\newblock Advancing X-cutting Ideas for Computational Climate Science (AXICCS)
  2016, Rockville, MD, 2016.

\bibitem{schein2021preserving}
Alexander Schein, Kevin~T Carlberg, and Matthew~J Zahr.
\newblock Preserving general physical properties in model reduction of
  dynamical systems via constrained-optimization projection.
\newblock {\em International Journal for Numerical Methods in Engineering},
  122(14):3368--3399, 2021.

\bibitem{Schmid_JFM_2010}
P.J. Schmid.
\newblock Dynamic mode decomposition of numerical and experimental data.
\newblock {\em Journal of Fluid Mechanics}, 656:5--28, 2010.

\bibitem{singh1}
N.~P. {Singh} and K.~{Ahuja}.
\newblock Reusing preconditioners in projection based model order reduction
  algorithms.
\newblock {\em IEEE Access}, 8:133233--133247, 2020.

\bibitem{singh2}
Navneet~Pratap Singh and Kapil Ahuja.
\newblock Preconditioned linear solves for parametric model order reduction.
\newblock {\em International Journal of Computer Mathematics},
  97(7):1484--1502, 2020.

\bibitem{sirisup}
S.~Sirisup and G.E. Karniadakis.
\newblock A spectral viscosity method for correcting the long-term behavior of
  {POD} models.
\newblock {\em J. Comput. Phys.}, 194:92--116, 2004.

\bibitem{sirovich}
L.~Sirovich.
\newblock Turbulence and the dynamics of coherent structures, part iii:
  dynamics and scaling.
\newblock {\em Q. Appl. Math.}, 45(3):583--590, 1987.

\bibitem{Aeras:2015}
W.~Spotz, T.~Smith, I.~Demeshko, and J.~Fike.
\newblock Aeras: A next generation global atmosphere model.
\newblock {\em Procedia Computer Science}, 51:2097--2106, 2015.

\bibitem{LCM:2013}
W.~Sun, J.~Ostien, and A.~Salinger.
\newblock A stabilized assumed deformation gradient finite element formulation
  for strongly coupled poromechanical simulations at finite strain.
\newblock {\em International Journal for Numerical and Analytical Methods in
  Geomechanics}, 37:2755--2788, 2013.

\bibitem{fy16_sand}
I.~Tezaur, J.~Fike, K.~Carlberg, M.~Balajewicz, M.~Barone, and E.~Mussoni.
\newblock Model reduction for compressible cavity simulations towards
  uncertainty quantification of structural loading.
\newblock Sandia National Laboratories Report, Sand No. 2016-9463. Sandia
  National Laboratories, Albuquerque, NM, 2016.

\bibitem{fy18_sand}
I.~Tezaur, J.~Fike, K.~Carlberg, and M.~Barone.
\newblock Captive carry reduced order modeling.
\newblock Sandia National Laboratories Report, Sand No. 2018-10824. Sandia
  National Laboratories, Albuquerque, NM, 2018.

\bibitem{fy17_sand}
I.~Tezaur, J.~Fike, K.~Carlberg, M.~Barone, D.~Maddix, E.~Mussoni, and
  M.~Balajewicz.
\newblock Advanced fluid reduced order models for compressible flow.
\newblock Sandia National Laboratories Report, Sand No. 2017-10335. Sandia
  National Laboratories, Albuquerque, NM, 2017.

\bibitem{Tezaur:2015}
I.~Tezaur, M.~Perego, A.~Salinger, R.~Tuminaro, and S.~Price.
\newblock Albany/felix: A parallel, scalable and robust finite element
  higher-order stokes ice sheet solver built for advanced analysis.
\newblock {\em Geoscientific Model Development}, 8:1--24, 2015.

\bibitem{Tiso:2013}
P.~Tiso and D.~Rixen.
\newblock Discrete empirical interpolation method for finite element structural
  dynamics.
\newblock {\em Topics in Nonlinear Dynamics}, 1:203--212, 2013.

\bibitem{veroy}
K.~Veroy and A.T. Patera.
\newblock Certified real-time solution of the parametrized steady
  incompressible {N}avier-{S}tokes equations: rigorous reduced-bases \textit{a
  posteriori} error bounds.
\newblock {\em J. Num. Meth. Fluids}, 47:773--788, 2005.

\bibitem{kyle}
K.~Washabaugh.
\newblock Faster fidelity for better design.
\newblock Ph.D. Thesis, Stanford University, 2016.

\bibitem{willcox}
K.~Willcox and J.~Peraire.
\newblock Balanced model reduction via the proper orthogonal decomposition.
\newblock {\em AIAA Journal}, 40(11):2323--2330, 2002.

\bibitem{zahm}
O.~Zahm and A.~Nouy.
\newblock Interpolation of inverse operators for preconditioning
  parameter-dependent equations.
\newblock {\em SIAM J. Sci. Comput.}, 38:A1044--A1074, 2016.

\end{thebibliography}
\bibliographystyle{plain}

\end{document}